\documentclass{amsart}

\usepackage [mathscr]{eucal}
\usepackage{amsmath,amsthm,amssymb,amscd,color}
\usepackage{amsrefs}
\usepackage{epsf}
\usepackage{amsfonts}
\usepackage{dsfont}
\usepackage{latexsym}
\usepackage{mathrsfs}
\usepackage{layout}
\usepackage{bm}
\usepackage{epsfig}
\usepackage{subfigure}
\usepackage{here}
\usepackage{caption}
\usepackage{graphics}
\usepackage{graphicx}

\newtheorem{theorem}{Theorem}[section]
\newtheorem{lemma}[theorem]{Lemma}
\newtheorem{proposition}[theorem]{Proposition}
\newtheorem{corollary}[theorem]{Corollary}

\newtheorem{remark}[theorem]{Remark}
\newtheorem{definition}[theorem]{Definition}

\newtheorem{assumption}[theorem]{Assumption}
\numberwithin{equation}{section}

\newcommand{\R}{{ \mathbb{R}  }}

\newcommand{\abs}[1]{\lvert#1\rvert}
\newcommand{\bke}[1]{\left( #1 \right)}

\newcommand{\norm}[1]{\left\Vert #1 \right\Vert}

\def\XXint#1#2#3{{\setbox0=\hbox{$#1{#2#3}{\int}$ }
\vcenter{\hbox{$#2#3$ }}\kern-.6\wd0}}

\def\BigRoman{\uppercase\expandafter{\romannumeral\number\count255 }}
\def\Romannumeral{\afterassignment\BigRoman\count255=}

\DeclareMathOperator *{\essosc}{ess\ osc}
\DeclareMathOperator *{\osc}{osc}
\DeclareMathOperator *{\esssup}{ess\ sup}
\DeclareMathOperator *{\essinf}{ess\ inf}
\DeclareMathOperator *{\di}{div} 
\DeclareMathOperator *{\meas}{meas}
\DeclareMathOperator *{\dist}{dist}
\DeclareMathOperator *{\data}{data}
\DeclareMathOperator *{\diam}{diam}
\DeclareMathOperator *{\loc}{loc}
\DeclareMathOperator *{\Lip}{Lip}
\DeclareMathOperator *{\Tr}{Tr}
\DeclareMathOperator *{\BMO}{BMO}
\DeclareMathOperator *{\Proj}{Proj}

\newcommand{\calC}{{\mathcal C}}

\newenvironment{thm1.5}{{\par\noindent\bf
           Proof of Theorem \ref{weak-1}. }}
           {\hfill\fbox{}\par\vspace{.2cm}}

\newenvironment{thm1.6}{{\par\noindent\bf
           Proof of Theorem \ref{weak-2}. }}
           {\hfill\fbox{}\par\vspace{.2cm}}

\newenvironment{thm1.7}{{\par\noindent\bf
           Proof of Theorem \ref{weak-4}. }}
           {\hfill\fbox{}\par\vspace{.2cm}}

\newenvironment{thm1.8}{{\par\noindent\bf
           Proof of Theorem \ref{bdd weak-1}. }}
           {\hfill\fbox{}\par\vspace{.2cm}}

\newenvironment{thm1.9}{{\par\noindent\bf
           Proof of Theorem \ref{bdd weak-2}. }}
           {\hfill\fbox{}\par\vspace{.2cm}}

\newenvironment{thm1.10}{{\par\noindent\bf
           Proof of Theorem \ref{bdd weak-4}. }}
           {\hfill\fbox{}\par\vspace{.2cm}}

\begin{document}
\title[Keller-Segel-fluid model]{H\"older continuity of Keller-Segel equations of porous medium type coupled to fluid equations}

\author{Yun-Sung Chung}
\address{Department of Mathematics, Yonsei University, Seoul, Republic of Korea}
\email{ysjung93@hanmail.net}
\author{Sukjung Hwang}
\address{Center for Mathamatical Analysis and Computation, Yonsei University}
\email{sukjung{\_}hwang@yonsei.ac.kr}
\author{Kyungkeun Kang}
\address{Department of Mathematics, Yonsei University, Seoul, Republic of Korea}
\email{kkang@yonsei.ac.kr}
\author{Jaewoo Kim}
\address{Department of Mathematics, Yonsei University, Seoul, Republic of Korea}
\email{baseballer@skku.edu}

\begin{abstract}
We consider a coupled system consisting of a degenerate porous
medium type of Keller-Segel system and Stokes system modeling the
motion of swimming bacteria living in fluid and consuming oxygen. We
establish the global existence of weak solutions and H\"older
continuous solutions in dimension three, under the assumption that
the power of degeneracy is above a certain number depending on given
parameter values. To show H\"older continuity of weak solutions, we
consider a single degenerate porous medium equation with lower order
terms, and via a unified method of proof, we obtain H\"older regularity,
which is of independent interest.
\end{abstract}

\maketitle

\section{Introduction}

We study a Keller-Segel model coupled to the fluid equations, where
the equation of biological cells is of porous medium type. To be
more precise, we consider
\begin{eqnarray}\label{eq:Chemotaxis}
\mbox{(KS-PME)}\quad \left\{
\begin{array}{cl}
& \partial_{t}n - \Delta n^{1+\alpha} + u\cdot \nabla n = - \nabla
\cdot (\chi (c)n^{q}\nabla c),
\\
& \partial_{t}c - \Delta c + u\cdot \nabla c = -\kappa (c)n,
\\
& \partial_{t}u - \Delta u + \nabla p = -n\nabla \phi,
\\
& \nabla \cdot u=0,
\end{array}
\right.
\end{eqnarray}
where $\alpha > 0$ and $q \geq 1$ are given constants. Here, the
unknowns $n$, $c$, $u$ and $p$ denote the density of bacteria, the
oxygen concentration, the velocity vector of the fluid and the
associated pressure, respectively. In addition, the locally bounded
functions $\chi : \mathbb{R}\rightarrow \mathbb{R}$ and $\kappa :
\mathbb{R}\rightarrow \mathbb{R}$ represent the chemotactic
sensitivity and consumption rate of oxygen. Moreover, $\phi =
\phi(x)$ is a given potential function. It is known that the above
system models the motion of swimming bacteria, so called
$\textit{Bacillus subtilis}$, which live in fluid and consume
oxygen.  This system has been proposed by Tuval \emph{et al.} in
\cite{TCDWKG} for the case $\alpha=0$ and $q=1$, which can be
extended to the case $\alpha>0$ when the diffusion of bacteria is
viewed like movement in a porus medium. In this manuscript, we call
the above system a Keller-Segel porous medium equation(KS-PME),
since fluid equations are restricted to the Stokes system under our
considerations.


The main purpose of this paper is to establish the existence of weak
and H\"older continuous solutions globally in time for the Cauchy
problem of (KS-PME) under general conditions of $\chi$ and $\kappa$
and more extended range of $\alpha$ and $q$ ever known.

We first introduce local H\"{o}lder regularity results for a scalar
equation under proper conditions on the lower order term, that
contributes later obtaining H\"{o}lder continuity of a weak solution
of system (KS-PME).

In the domain $\Omega_{T} \subset \mathbb{R}^{d}\times [0, T]$ for $d\geq 2$,
we consider parabolic porous medium type equations in the form of
\begin{equation}\label{E:PM}
n_{t} = \Delta n^{1+\alpha} + \nabla \cdot \left(B(x,t) n\right)
\end{equation}
for $\alpha \geq 0$ under proper conditions on $B$ where $B: \mathbb{R}^{d} \to \mathbb{R}^{d}$ is a vector field.
Roughly speaking, if we are able to obtain regularity results of
\eqref{E:PM} under the condition on $B$ that is expected from  $u, \chi(c),$ and $\nabla c$ of (KS-PME), then H\"{o}lder continuity of a weak solution $n$ of \eqref{E:PM} yields the same regularity for $n$, a weak solution of (KS-PME).

In fact, our method of showing H\"{o}lder continuity of \eqref{E:PM} works under the conditions on $B$ and $\nabla B$ such that
\begin{equation}\label{B}
B \in L^{2\hat{q}_1, 2\hat{q}_2}_{\text{loc}} (\Omega_{T}) \quad
\text{and } \quad \nabla B \in L^{\hat{q}_1, \hat{q}_2}_{\text{loc}}
(\Omega_{T})
\end{equation}
where positive constants $\hat{q_1}, \hat{q_2} > 1$ satisfy
\begin{equation}\label{hatq-hatr}
\frac{2}{\hat{q}_2} + \frac{d}{\hat{q}_1} = 2 - d \kappa
\end{equation}
for some $\kappa \in (0, 2/d).$ By letting
\[
q_1 = \frac{2\hat{q}_1 (1+\kappa)}{\hat{q}_1-1}, \quad q_2 =
\frac{2\hat{q}_2 (1+\kappa)}{\hat{q}_2 -1},
\]
the admissible range of constants are obtained from
Proposition~\ref{P:admissible} when $p=2$.

There are many papers working on the continuity of weak
solutions to porous medium type equations (refer \cite{DB83},
\cite{Ar}, \cite{CaFr}, et al. for a general porous medium equation
and special classes of equations). Focusing the main term of \eqref{E:PM}, we share some common mathematical approaches.

For the system (KS), we refer recent paper \cite{KL} carrying
H\"{o}lder regularity and uniqueness results (when $\alpha >0$)
relying on technical proofs originated from \cite{ChDB88} and
\cite{DB93}. Compare to similar H\"{o}lder regularity results on
\cite{KL}, we play with a scalar equation \eqref{E:PM} to obtain the
same results under the weaker assumptions on $B$ and $\nabla B$ that
belongs to scaling invariant class. By following
natural behaviour of a solution using a more geometrical approach (refer \cite{HL15a} and \cite{HL15b}), as a separate interest of its own, we provide a unified method of
proof in the sense that the method has no limitation including
$\alpha =0$ (usually it is important to have $\alpha >0$ in
\cite{KL},\cite{ChDB88}, and \cite{DB93} and showing stability when
$\alpha \to 0$ is regarded as an another computational issue). Besides simplicity of
computations in this manuscript, our method of proof carries potentials to provide significant common
elements to the similar proofs for singular type of equations(when
$-1 <\alpha \leq 0$) and even for generalized structured equations (refer Remark~\ref{R:gen} for details).

Here we provide the definition of a weak solution of \eqref{E:PM}.
\begin{definition}\label{D:WS}
Let $\Omega$ be an open set in $\mathbb{R}^{d}$, $B\in
L^{2}((0,T)\times \Omega)$, and $T>0$.
\[
n \in C_{\loc}(0,T; L^{2}_{\loc}(\Omega)), \quad n^{\frac{\alpha
+2}{2}} \in L^{2}_{\loc}(0,T; W^{1,2}_{\loc}(\Omega))
\]
is a local weak solution to \eqref{E:PM} with $\alpha \geq 0$ if for
every compact set $K \subset \Omega$ and every subinterval $[t_1,
t_2] \subset (0,T]$

\begin{equation}\label{WS}
\left. \int_{K} n \varphi \,dx \right|_{t_0}^{t_1} +
\int_{t_0}^{t_1}\int_{K} \left\{ -n\varphi_{t} + \nabla n^{1 +
\alpha}\nabla \varphi + B n \nabla \varphi \right\} \,dx\,dt = 0
\end{equation}
for all nonnegative testing functions
\[
\varphi \in W^{1,2}_{\loc}(0,T; L^2(K)) \cap L^{2}_{\loc}(0, T;
W^{1,2}_{o}(K)).
\]
\end{definition}

From the definition of weak solutions, we compute two types of energy estimates, provided in Propositions~\ref{P:EE} and \ref{P:LEE} which is called local and logarithmic energy estimates, respectively. Due to the difference of the nature of porous medium and $p-$Laplacian equations, we modify the method of proof in \cite{HL15a}, for example, considering the a weak solution directly rather sub or super solutions, also cutting off a weak solution when $u$ may stay near zero for DeGiorgi iteration. Moreover, another technical issue follows
because the lower order term in \eqref{E:PM} does not follow the structure of main term (not given in the form of $n^{1+\alpha}$ but $n$). By imposing conditions of $B$ and $\nabla B$ in scaling invariant class, we can provide simpler proof compare to computation in \cite{KL}. Also conditions on $\nabla B$ does follow global estimates from (KS-PME).

Before we deliver the local H\"{o}lder continuity results, we make comments on intrinsic scaling due to the nonhomogeneity of the equation \eqref{E:PM}. More precisely, the local energy estimate derived from \eqref{E:PM} appears in
Proposition~\ref{P:EE} is nonhomegenous unless $\alpha = 0$. Roughly speaking, in an intrinsically (rescaled with the behaviour of a solution) scaled cylinder, a weak solution behaves like a solution to the heat equation.
That is, more specifically, rescaling the time length
\begin{equation}\label{intrinsic}
T_{\omega, \rho} = \theta \omega^{-\alpha} \rho^2
\end{equation}
for some constant $\theta$ and $\rho$ and
\begin{equation}\label{omega}
\omega := \essosc_{\Omega_{T}} n = \mu_{+} - \mu_{-} : =\esssup_{\Omega_T} n -
\essinf_{\Omega_T} n.
\end{equation}

Since $\Omega_{T}$ is open, there are positive constants $r$ and $s$
such that $K^{x_0}_{r} \times (t_0 - s, t_0) \subset \Omega_{T}$. If
we set
\[
R = \frac{1}{4} \min \left\{ r,\ \frac{\omega^{\alpha/2}
s^{1/2}}{\theta^{1/2}}\right\},
\]
then we conclude that
\[
Q^{x_0, t_0}_{\omega, 4R}(\theta) = K^{x_0}_{4R} \times (t_0 -
\theta \omega^{-\alpha}R^2, t_0 ) \subset \Omega_{T}.
\]
Then for any positive constants $\theta$ and $\omega$, we can fit the cylinder $Q^{x_0,
t_0}_{\omega, 4R}(\theta)$ in $\Omega_{T}$ by selecting $R$
properly. Basically, we are going to work with the cylinder $Q^{x_0,
t_0}_{\omega, 4R}(\theta)$ to find a proper subcylinder where a
solution has less oscillation eventually leading to H\"{o}lder
continuity.

Due to the intrinsic scaling \eqref{intrinsic}, we define a time scale in
terms of the function $n$ and the set $\Omega$ on which $n$ is
defined. For any real number $\tau$, we define
\begin{equation}\label{tI}
|\tau|_{I} = \omega^{\alpha/2}|\tau|^{1/2}.
\end{equation}
With this time scale, we define the parabolic distance between two
sets such $\mathcal{K}_{1}$ and $\mathcal{K}_2$ by
\[
\dist_{p}(\mathcal{K}_1 ; \mathcal{K}_2) = \inf_{\substack {(x,t) \in \mathcal{K}_{1}\\
 (y,s)\in \mathcal{K}_{2}\\
 s\le t} }
 \max\{|x-y|_\infty, \  |t-s|_{I} \}
\]
with $|\cdot|_{\infty}$ (which is defined by $|x-y|_{\infty} =
\max_{1\leq i\leq d} |x^{i} - y^{i}|$). 

Now we state the H\"{o}lder continuity of a bounded weak
solution of \eqref{E:PM}.
\begin{theorem}\label{T:Holder1} (H\"{o}lder continuity of $n$)
Let $n$ be a nonnegative bounded weak solution of \eqref{E:PM} under \eqref{B}
with $\alpha \geq 0$ in $\Omega_{T}$. Then $n$ is locally
continuous. Moreover, there exist positive constant $\beta \in
(0,1)$ and $\gamma$ depending on data(that is, $d, \Omega_{T},
\Omega'_{T}, \alpha, \|B\|_{2\hat{q}_1, 2\hat{q}_2}, \|\nabla
B\|_{\hat{q}_1, \hat{q}_2}$ for some $\hat{q}_1, \hat{q}_2 > 1$ satisfying
\eqref{hatq-hatr}) such that, for any two distinct points $(x_1,
t_1)$ and $(x_2, t_2)$ in any subset $\Omega'_{T}$ of $\Omega_{T}$
with $\dist (\Omega'_{T};
\partial_{p}\Omega_{T})$ positive, we have
\begin{equation}\label{Holder}
\left|n(x_1, t_1) - n(x_2, t_2)\right| \leq \gamma \omega
\left(\frac{|x_1 - x_2| + \omega^{\alpha/2}|t_1 - t_2|^{1/2}}{\dist_{p}
(\Omega'_{T}; \partial_{p} \Omega_{T})}\right)^{\beta}.
\end{equation}
\end{theorem}

The proof of this theorem is given in Section~\ref{S:ProofHolder} considering two alternatives. Then the proofs of two alternatives are shown in Section~\ref{S:ProofAlternatives} as combinations of DeGorgi iterations and the expansion of positivity along the time axis and the spatial axis.


Now we state results on the existence of global-intime weak solution of (KS-PME) and global H\"{o}lder continuity of the Cauchy problems of (KS-PME) as well.
For the notational convenience, we denote
\[
A:=\left\{(q, \alpha)\mid \alpha > 2q-2, q\ge 1 \right\},\quad
B:=\left\{(q, \alpha)\mid \alpha
> \frac{9q-8}{6}, q\ge 1 \right\},
\]
\[
C:=\left\{(q, \alpha)\mid \alpha
> \frac{10q-9}{8}, q\ge 1 \right\}.
\]
We introduce the notions of weak solutions and H\"older continuous
solutions. We start with the definition of weak solutions.
\begin{definition}\label{Def:weak}\,\,(Weak solutions)\,\,
Let $q\geq 1$
and $0<T<\infty$. A triple
$(n,c,u)$ is said to be a \emph{weak solution} of the system
(\ref{eq:Chemotaxis}) if the followings are satisfied:
\begin{itemize}
\item[$({\romannumeral 1})$]$n$ and $c$ are non-negative functions and $u$
is a vector function defined in $\mathbb R^3\times (0,T)$ such that
\[
 n(1+|x|+|\log n|)\in L^\infty(0,T;L^1(\mathbb
R^3)),\quad\nabla n^{\frac{1+\alpha}{2}}\in L^2(0,T;L^2(\mathbb
R^3)),
\]
\[
c\in L^\infty(0,T;H^1(\mathbb R^3))\cap L^2(0,T;H^2(\mathbb
R^3)),\quad c\in L^{\infty}(\mathbb R^3\times [0,T)),
\]
\[
u\in L^\infty(0,T;L^2(\mathbb R^3)),\quad \nabla u\in
L^2(0,T;L^2(\mathbb R^3)),
\]

\item[$({\romannumeral 2})$]$(n,c,u)$ satisfies the system
  (\ref{eq:Chemotaxis}) in the sense of distributions, namely,
    \begin{eqnarray*}
           \int_0^T\int_{\mathbb R^d}\left(n\varphi_t-\nabla
            n^{1+\alpha}\cdot\nabla\varphi+nu\cdot\nabla\varphi+n^{q}\chi(c)\nabla
            c\cdot\nabla\varphi\right) dxdt=-\int_{\mathbb R^d}
            n_0\varphi(\cdot,0) ~dx,\\
            \int_0^T\int_{\mathbb R^d}\left(c\varphi_t-\nabla
            c\cdot\nabla\varphi+cu\cdot\nabla\varphi+n\kappa(c)\varphi\right) dxdt=-\int_{\mathbb R^d}
            c_0\varphi(\cdot,0) ~dx,\\
            \int_0^T\int_{\mathbb R^d}\left( u\cdot\psi_t-\nabla
            u\cdot\nabla\psi+\left(\tau(u\cdot\nabla \psi)\right) u-n\nabla\phi\cdot\psi
            \right)dxdt=-\int_{\mathbb R^d}u_0\cdot\psi(\cdot,0) ~dx\\
    \end{eqnarray*}
for any $\varphi\in C^\infty_0\left(\mathbb R^3\times [0,T)\right)$
and $\psi\in C^\infty_0\left(\mathbb R^3\times [0,T),\mathbb
R^3\right)$ with $\nabla\cdot\psi=0.$
\end{itemize}
\end{definition}

Next we define H\"older continuous solutions. For convenience, we
denote $Q_T:=(0,T)\times \mathbb R^3$.
\begin{definition}\label{Def:bdd-weak}
\,\,(H\"older continuous solutions)\,\, Let $q\geq 1$,
$(q,\alpha)\in (A\cap B)$ and $0<T<\infty.$ A triple $(n,c,u)$ is
said to be a \emph{H\"older continuous solutions} of the system
(\ref{eq:Chemotaxis}) if $(n,c,u)$ is a weak solution in Definition
\ref{Def:weak} and furthermore satisfies the following: there exists
$\beta>0$ such that
\begin{equation}\label{Oct3-CHKK-10}
n,\,\, \partial_{t}c,\,\, \partial_{t}u, \,\,\nabla^2 c,\,\,
\nabla^2 u\in \calC^\beta(Q_T).
\end{equation}
\end{definition}

Before stating our result precisely, we first recall some essential
conditions for $\chi$ and $\kappa.$ To preserve the non-negativity
of the density of bacteria $n(x,t)$ and the oxygen $c(x,t)$ for
 $0<t<T,$ it is necessary to assume that $\kappa(0)=0$.
The condition $\kappa(\cdot)\geq 0$ is also essential since the
bacteria consume the oxygen. Thus, the following hypotheses are
compulsory: $\kappa(\cdot)\geq 0$ and $\kappa(0)=0$. Furthermore, we
suppose that  $\chi' \in L^{\infty}_{\rm loc}$. Summing up,
throughout this thesis, we assume that
\begin{itemize}
  \item [$(P_1)$] $\chi' \in L^{\infty}_{\rm loc}$,\,\, $\kappa \in L^{\infty}_{\rm loc}$, \,\,$\kappa(\cdot)\geq 0$\,\, and\,\,
  $\kappa(0)=0.$
\end{itemize}
To obtain more extended range of $\alpha$, we sometimes make further
assumptions on $\kappa$, which are given by
\begin{itemize}
  \item [$(P_2)$] $\kappa'\in L^\infty_{\rm loc}$ with$\quad\kappa'(\cdot)\geq \kappa_0
     \mbox{ for some constant } \kappa_0>0.$
\end{itemize}
We now present two different types of assumptions on $\chi$,
$\kappa$ together with the range of $\alpha$ and $q$. The first one
is reserved for weak solutions.
\begin{assumption}\label{assume-1}
$\chi,\kappa$ and $\alpha$ satisfy $(P_{1})$ and one of the
following holds:
\begin{itemize}
\item[$(\romannumeral 1)$] $
(q, \alpha)\in B$.
\item[$(\romannumeral 2)$] $(q,\alpha)\in A \cup B$ and $\kappa$ satisfies
$(P_2)$.
\end{itemize}
\end{assumption}
Next assumption is prepared for H\"{o}lder continuous solutions.
\begin{assumption}\label{assume-2}
$\chi,\kappa$ and $\alpha$ satisfy $(P_1)$ and one of the following
holds:
\begin{itemize}
\item[$(\romannumeral 1)$] $
(q, \alpha)\in A\cap B$.
\item[$(\romannumeral 2)$] $(q,\alpha)\in (A\cup B) \cap C$ and $\kappa$ satisfies $(P_2)$.
\end{itemize}
\end{assumption}
We recall some known results related to our concerns. Firstly, we
compare the system (KS-PME) to the classical Keller-Segel model (KS)
of porous medium type, which is given as
\begin{equation}\label{eq:Keller-Segel}
        \mbox{(KS)}\quad\left\{
            \begin{array}{l}
                \partial_t n=\Delta
                n^{1+\alpha}-\nabla\cdot\left(\chi n\nabla
                c\right),
                \\
                \vspace{-3mm}\\
                \tau\partial_t c=\Delta c-c+n,
            \end{array}
         \right.
\end{equation}
where $\chi$ is a positive constant, $q=1$ and $\tau=0$ or $1$(for
example, \cite{KS1, KS2}). 
We remark that the equation of $c$ in \eqref{eq:Keller-Segel} is
modeled by the chemical substance, which is produced by biological
organism, but in our case the equation \eqref{eq:Chemotaxis}$_2$
indicates the dynamics of oxygen, which is consumed by a certain
type of bacteria. That's the reason opposite sign of the right side
of each equation appears, which causes main difference regarding
global existence or blow-up for the value on $\alpha$. In case that
\eqref{eq:Keller-Segel}, the equation of $c$, is of elliptic type,
i.e. $\tau=0$, existence of bounded weak solutions was shown  in
\cite{SK} globally in time, provided that $q \geq 1$ and
$\alpha>q-\frac{2}{d}$. If $0<\alpha\le q-\frac{2}{d}$, blow-up may
occur in a finite time. Later, in \cite{IY}, the result of \cite{SK}
was extended to the case that the equation of $c$ is of parabolic
type, i.e. $\tau>0$.

\indent For the chemotaxis fluid system \eqref{eq:Chemotaxis} with
$q=1$ in two dimensions, it was known that bounded weak solutions
exist globally in time under some assumptions on $\kappa$ and $\chi$
for sufficiently regular data. We remark that results in dimension
two are even valid in replacement with the Navier-Stokes equations
for fluid equations (refer to \cite{CKK} and \cite{TW_1})

In three dimensions, it was shown in \cite{LL} that if $\alpha
=\frac{1}{3}$, then the chemotaxis-Stokes system
\eqref{eq:Chemotaxis} with $q=1$  has global-in-time  bounded weak
solutions. For the special case that $\chi=1$ and $\kappa(x)=x$,
existence of bounded weak solutions was proved in \cite{TW_2} for
(KS-PME) with $q=1$, provided that $\alpha> 1/7$. In
\cite{CKK}, for \eqref{eq:Chemotaxis} with $q=1$, it was proved that
global-in-time existence of weak solutions and bounded weak
solutions under the same conditions as $(P1)$ and $(P2)$, if
$\alpha>1/6$ and $\alpha>1/4$, respectively. The range of $\alpha$
was improved in \cite{CK}. More precisely, if $\alpha>\frac{1}{6}$,
bounded weak solutions exist under only the condition $(P1)$.
Furthermore, it was also proved that if $\chi$ or $\kappa$ satisfy
$(P1)$ and $(P2)$, and if $\alpha > \frac{1}{8}$, then there exists
bounded weak solutions for the system \eqref{eq:Chemotaxis} with
$q=1$.

As mentioned earlier, our main goal is to study more general
Keller-segel-fluid system \eqref{eq:Chemotaxis} with $q\ge 1$ and
obtain global existence of weak and H\"older continuous solutions
for extended range of $\alpha$ and $q$. Our results are summarized
in the Table 1. We remark that in case that $q=1$, our results
recover those of \cite{CK}.

\begin{table}
\caption{Relations between parameters and conditions} \scriptsize
\begin{center}
\begin{tabular}{|c|c|c|}
    \hline
               &weak solutions &H\"older continuous
solutions\\
    \hline
    \hline
    $(P_1)$ &$\alpha > \frac{9q-8}{6}$ &$\alpha > \text{max}\left\{ 2q-2, \ \frac{9q-8}{6} \right\}$\\
     &(Theorem \ref{weak-1}) &(Theorem \ref{bdd weak-1})\\
    \hline
    $(P_1)$ and $(P_2)$ &$\alpha > \text{min}\left\{
2q-2, \frac{9q-8}{6}\right\}$ &$\alpha > \text{max}\left\{
\text{min}\left\{ 2q-2,
\frac{9q-8}{6}\right\}, \ \frac{10q-9}{8} \right\}$\\
    &(Theorem \ref{weak-2}) &(Theorem \ref{bdd weak-2})\\
    \hline
\end{tabular}
\end{center}
\end{table}

Now we are ready to state our main results of the system \eqref{eq:Chemotaxis},
and the first one is about existence of weak solutions, which reads
as follows:

\begin{theorem}\label{weak-1} $\bke{\text{Weak solutions}}$
Let $\alpha$ belong to B i.e., $\alpha > \frac{9q-8}{6}$ and initial
data $(n_{0}, c_{0}, u_{0})$ satisfy
\begin{equation}\label{weak-initial-1}
n_{0}(1+|x|+|\log n_{0}|) \in L^{1}(\mathbb{R}^{3}), \ c_{0} \in
L^{\infty}(\mathbb{R}^{3}) \cap H^{1}(\mathbb{R}^{3})  \ \text{and}
\ u_{0} \in L^{2}(\mathbb{R}^{3}).
\end{equation}
Suppose that $\chi$, $\kappa$ satisfy the hypothesis $(P_1)$. Then,
there exists a weak solution $(n,c,u)$ for the system
\eqref{eq:Chemotaxis}. Furthermore, for any $p$ with $1\leq p\leq
\alpha -q+2$
\[
n\in L^\infty(0,T;L^p(\mathbb R^3)),\qquad \nabla
n^{\frac{p+\alpha}{2}}\in L^2(0,T;L^2(\mathbb R^3)),
\]
and the following inequality is satisfied:
\begin{equation*}
\begin{split}
& \sup_{0\leq t \leq T}\left(\int_{\mathbb R^3} n\left(\abs{\log n}
+ \langle x\rangle \right) + \int_{\mathbb{R}^3}n^{\alpha -q+2}
+\int_{\mathbb{R}^3}\mid \nabla
c \mid^{2}  + \int_{\mathbb{R}^3}\abs{u}^{2}  \right) \\
& \hspace{0.8cm} + C \int_{0}^{T}\left( \norm{\nabla
n^{\frac{1+\alpha}{2}}}_{2}^{2}+ \norm{\nabla
n^{\frac{2\alpha-q+2}{2}}}_{2}^{2}+\norm{\Delta
c}_{2}^{2}+\norm{\nabla
u}_{2}^{2} \right) \\
& \leq C \left( T, \norm{c_0}_{L^{\infty}\cap H^{1}},
\norm{n_{0}(1+|x|+|\log n_{0}|)}_1, \norm{n_0}_{\alpha -q +2},
\norm{u_0}_2 \right),
\end{split}
\end{equation*}
where $\langle x\rangle=(1+|x|^2)^{\frac{1}{2}}$.
\end{theorem}
If the condition $(P2)$ is additionally assumed, the range of
$\alpha$ is a bit expanded, compared to that of Theorem
\ref{weak-1}. More precisely, we have the following:

\begin{theorem} \label{weak-2}$\bke{\text{Weak solutions}}$
Let $\alpha$ belong to $A\cup B$ i.e., $\alpha > \text{min}\left\{
2q-2, \frac{9q-8}{6}\right\} $. Suppose that $\chi$, $\kappa$
satisfy the hypothesis $(P_1)$, $(P_2)$ and initial data $(n_{0},
c_{0}, u_{0})$ satisfies
\begin{equation}\label{weak-initial-2}
n_{0}(1+|x|+|\log n_{0}|) \in L^{1}(\mathbb{R}^{3}), \ c_{0} \in
L^{\infty}(\mathbb{R}^{3}) \cap H^{1}(\mathbb{R}^{3}) \ \text{and} \
u_{0} \in H^{1}(\mathbb{R}^{3}).
\end{equation}
Then, there exists a weak solution $(n,c,u)$ for the system
(KS-PME). Furthermore, for any $p$ with $1\leq p\leq
\alpha -2q+3$
\[
n\in L^\infty(0,T;L^p(\mathbb R^3)),\quad n^{\frac{1}{2}}\nabla c
\in L^2(0,T;L^2(\mathbb R^3)),\quad \nabla n^{\frac{p+\alpha}{2}}\in
L^2(0,T;L^2(\mathbb R^3)),
\]
and the following inequality is satisfied:
\begin{equation*}
\begin{split}
& \sup_{0 \leq t \leq T}\left(\int_{\mathbb R^3} n\left(\abs{\log n}
+ \langle x\rangle \right)  + \int_{\mathbb{R}^3}n^{\alpha -2q+3}
+\int_{\mathbb{R}^3}\abs{ \nabla
c }^{2}  + \int_{\mathbb{R}^3}\abs{u}^{2}  \right) \\
& \hspace{0.3cm} + C \int_{0}^{T}\left( \norm{\nabla
n^{\frac{1+\alpha}{2}}}_{2}^{2}+ \norm{\nabla
n^{\frac{2\alpha-2q+3}{2}}}_{2}^{2}+\norm{\Delta
c}_{2}^{2}+\norm{\nabla
u}_{2}^{2} \right) \\
& \leq C \left( T, \norm{\nabla c_0}_2, \norm{n_{0}\log n_{0}}_1,
\norm{n_0}_{\alpha -2q +3}, \norm{u_0}_2 \right),
\end{split}
\end{equation*}
where $\langle x\rangle=(1+|x|^2)^{\frac{1}{2}}$.
\end{theorem}



Next, if $\alpha$ is greater than a certain value depending on $q$,
we prove existence of H\"older continuous solutions for (KS-PME)
under the condition $(P_1)$. To be more precise, the result reads as
follows:
\begin{theorem}\label{bdd weak-1}$\bke{\text{H\"older continuous
solutions}}$ Let $\alpha$ belongs to $A\cap B$ i.e., \\
$\alpha > \text{max}\left\{ 2q-2, \ \frac{9q-8}{6} \right\}$.
Suppose that $\chi$, $\kappa$ satisfy the hypothesis $(P_1)$ and
initial data $(n_{0}, c_{0}, u_{0})$ satisfies
\eqref{weak-initial-1} as well as
\begin{equation}\label{initial-bound-1}
n_{0} \in L^{\infty}(\mathbb{R}^{3}), \ c_{0} \in
W^{1,m}(\mathbb{R}^{3}), \ u_{0} \in W^{1,m}(\mathbb{R}^{3}), \
\text{for any} \ m < \infty.
\end{equation}
Then, there exists a H\"older continuous solution $(n,c,u)$ for the
system (KS-PME).
\end{theorem}
Furthermore, we assume the condition $(P_2)$ and we then see that
the restriction of $\alpha$ is relaxed for the existence of H\"older
continuous solutions.

\begin{theorem}\label{bdd weak-2}$\bke{\text{H\"older continuous
solutions}}$ Let $\alpha$ belongs to $(A\cup B) \cap C$ i.e.,
$\alpha > \text{max}\left\{ \text{min}\left\{ 2q-2,
\frac{9q-8}{6}\right\}, \ \frac{10q-9}{8} \right\}$. Suppose that
$\chi$, $\kappa$ satisfy the hypothesis $(P_1)$, $(P_2)$ and initial
data $(n_{0}, c_{0}, u_{0})$ satisfies \eqref{weak-initial-2} as
well as
\begin{equation}\label{initial-bound-2}
n_{0} \in L^{\infty}(\mathbb{R}^{3}), \ c_{0} \in
W^{1,m}(\mathbb{R}^{3}), \ u_{0} \in W^{1,m}(\mathbb{R}^{3}), \ \text{for any} \ m <
\infty.
\end{equation}
Then, there exists a H\"older continuous solution $(n,c,u)$ for the
system (KS-PME).
\end{theorem}

\begin{remark}
There are some known results regarding uniqueness of H\"older
continuous solutions for Keller-Segel system of porous medium type
(see e.g. \cite{MS} and \cite{KL}). As for us, uniqueness of
solutions in Theorem \ref{bdd weak-1} and Theorem \ref{bdd weak-2}
doesn't seem to be obvious, in particular, due to presence of the
fluid velocity field. Therefore, we leave it as an open question.
\end{remark}

This paper is organized as follows: In section 2, we introduce some
notations and review known results. Section 3 is devoted for the
proof of Theorem \ref{T:Holder1} with the crucial aid of two
alternatives, whose are clarified in section 4. In section 5, we
present the proofs of existence for weak solutions of (KS-PME)  in
dimension three. We also provide the proofs of Theorem \ref{bdd
weak-1} and Theorem \ref{bdd weak-2} in section 6. In appendix,
proofs of Propositions ~\ref{P:EE} and ~\ref{P:LEE} are given.


\section{Preliminaries}\label{Preliminaries}

\subsection{Notations and useful inequalities} In this subsection,
We introduce the notations throughout this paper and recall some
useful inequalities for our purpose. Let $\Omega$ be an open domain
in $\mathbb{R}^{d}$, $d\ge 1$ and $I$ a finite interval.
\[
L^{p}(\Omega)=\{ f : \Omega \rightarrow \mathbb{R} \mid f \text{ is
Lebesgue measurable }, \norm{f}_{L^{p}(\Omega)}<\infty\},
\]
 where
\[
 \norm{f}_{L^{p}(\Omega)}= \bke{\
 \int_{\Omega}|f|^{p}dx}^{\frac{1}{p}}, \qquad (1 \leq p < \infty).
\]
We will write $\norm{f}_{L^{p}(\Omega)}:=\norm{f}_{p}$, unless there
is any confusion to be expected. For $1\leq p \leq \infty$,
$W^{k,p}(\Omega)$ denotes the usual Sobolev space, i.e.,
\[
W^{k,p}(\Omega)=\{u\in L^{p}(\Omega):D^{\alpha}u\in L^{p}(\Omega), 0
\leq |\alpha|\leq k\}.
\]
We also write the mixed norm of $f$ in spatial and temporal
variables as
\[
\norm{f}_{L^{p,q}_{x,t}(\Omega\times I)}=\norm{f}_{L^{q}_{t}(I;
L^p_x(\Omega))}=\norm{\norm{f}_{L^p_x(\Omega)}}_{L^q_t(I)}.
\]
Let $m$ and $p$ be positive constants greater than $1$ and consider
the Banach spaces
\[
V^{m,p}(\Omega_{T}) := L^{\infty} (0, T; L^{m}(\Omega)) \cap
L^{p}(0,T; W^{1,p}(\Omega))
\]
and
\[
V^{m,p}_{0}(\Omega_{T}) := L^{\infty} (0, T; L^{m}(\Omega)) \cap
L^{p}(0,T; W^{1,p}_{0}(\Omega)),
\]
both equipped with the norm $v \in V^{m,p}(\Omega_{T})$,
\[
\| v\|_{V^{m,p}(\Omega_{T})} : = \esssup_{0<t<T} \|v(\cdot, t)\|_{m,
\Omega} + \| \nabla v\|_{p, \Omega_{T}}.
\]
When $m=p$, we set $V^{p,p}(\Omega_{T}) = V^{p}(\Omega_{T})$. Note
that both spaces are embedded in $L^{q}(\Omega_{T})$ for some $q >
p$. We denote by $C=C(\alpha,\beta,...)$ a constant depending on the
prescribed quantities $\alpha,\beta,...$, which may change from line
to line.

Now we introduce basic embedding inequalities and auxiliary lemmas
for fast geometric convergence. (Refer Chapter~I in \cite{DB93})

\begin{theorem}\label{T:GN} (Gagliardo-Nirenberg multiplicative embedding inequality)
Let $v \in W_{0}^{1,p}(\Omega)$, $p\geq 1$. For every fixed number
$s \geq 1$ there exists a constant $C$ depending only upon $d$, $p$
and $s$ such that
\[
\| v\|_{q, \Omega} \leq C \|\nabla v\|^{\alpha}_{p, \Omega}
\|v\|^{1-\alpha}_{s,\Omega},
\]
where $\alpha \in [0,1]$, $p, q \geq 1$, are linked by
\[
\alpha = \left( \frac{1}{s} - \frac{1}{q}\right) \left(\frac{1}{d} -
\frac{1}{p} + \frac{1}{s}\right)^{-1},
\]
and their admissible range is
\[\begin{cases}
q \in [s, \infty], \ \alpha \in [0, \frac{p}{p + s(p-1)}], & \text{ if } d = 1, \\
q\in [s, \frac{dp}{d-p}], \ \alpha \in [0,1], & \text{ if } 1\leq p < d, \ s \leq \frac{dp}{d-p}, \\
q\in [\frac{dp}{d-p}, s], \ \alpha \in [0,1], & \text{ if } 1\leq p < d, \ s \geq \frac{dp}{d-p}, \\
q \in [s, \infty), \ \alpha\in [0, \frac{dp}{dp + s(p-d)}), & \text{
if } 1<d\leq p .
\end{cases}\]
\end{theorem}

\begin{theorem}\label{T:Sobolev} (Sobolev embedding theorem) There exists a constant $C$ depending only upon $d, p, m$ such that for every $ v\in L^{\infty}(0,T; L^{m}(\Omega)) \cap L^{p}(0, T; W^{1,p}_{0}(\Omega))$,
\[
\iint_{\Omega_{T}} \left| v(x,t) \right|^{q} \,dx\,dt \leq C^{q}\left(\esssup_{0<t<T} \int_{\Omega} \left|
v(x,t)\right|^{m} \,dx \right)^{p/d} \left( \iint_{\Omega_{T}}
\left| \nabla v(x,t)\right|^{p}\,dx\,dt\right)
\]
where $q = \frac{p(d+m)}{d}$. Moreover,
\[
\|v\|_{q, \Omega_{T}}
\leq C \|v\|_{ V^{m,p} (\Omega_{T})} = C \left( \esssup_{0<t<T} \|v(\cdot, t)\|_{m,\Omega} + \|\nabla
v\|_{p, \Omega_{T}}\right).
\]
\end{theorem}

When $p=m$, we apply H\"{o}lder's inequality to obtain the following
corollary.

\begin{corollary}
Let $p>1$. There exists a constant $C$ depending only upon $d$ and
$p$, such that for every $v\in V_{0}^{p}(\Omega_{T})$,
\[
\| v \|^{p}_{p, \Omega_{T}} \leq C \left| \{ |v| > 0
\}\right|^{\frac{p}{d+p}} \|v\|^{p}_{V^{p}(\Omega_{T})}.
\]
\end{corollary}
\begin{proposition}\label{P:admissible}
There exists a constant $C$ depending only upon $d$ and $p$ such
that for every $v\in V^{p}_{0}(\Omega_{T})$,
\[
\| v \|_{q_1,q_2; \Omega_{T}} \leq C \|v\|_{V^{p}(\Omega_{T})}
\]
where the numbers $q_1, q_2 \geq 1$ are linked by
\[
\frac{1}{q_2} + \frac{d}{p q_1} = \frac{d}{p^2},
\]
and their admissible range is
\[\begin{cases}
q_1 \in (p,\infty), \ q_2\in (p^2, \infty); & \text{ for } d = 1, \\
q_1 \in (p, \frac{dp}{d-p}), \ q_2 \in (p, \infty); & \text{ for } 1 < p < d, \\
q_1 \in (p, \infty), \ q_2\in (\frac{p^2}{d}, \infty); & \text{ for } 1
< d \leq p .
\end{cases}\]
\end{proposition}

\begin{proof}
Let $v\in V_{0}^{p}(\Omega_{T})$ and let $r \geq 1$ to be chosen.
From Theorem~\ref{T:GN} with $s=p$ follows that
\[\begin{split}
&\left( \int_{0}^{T} \|v(\cdot, \tau)\|^{q_2}_{q_1, \Omega} \,d\tau\right)^{1/r} \\
&\leq C \left( \int_{0}^{T} \|\nabla v(\cdot, \tau)\|^{\alpha q_2}_{p}
\,d\tau\right)^{1/r} \esssup_{ 0\leq r \leq T} \| v (\cdot,
\tau)\|^{1-\alpha}_{p, \Omega_{T}}.
\end{split}\]
Choose $\alpha$ such that $\alpha q_2 = p$.
\end{proof}

We state a lemma concerning the geometric convergence
of sequences of numbers.


\begin{lemma}\label{L:iter2}
Let $\{Y_n\}$ and $\{Z_n\}$, $n=0,1,2, \ldots,$ be sequences of
positive numbers, satisfying the recursive inequalities
\[\begin{cases}
Y_{n+1} \leq C b^{n} \left( Y_{n}^{1+\alpha} + Z_{n}^{1+\kappa} Y_{n}^{\alpha}\right) & \\
Z_{n+1} \leq C b^{n}\left( Y_{n} + Z_{n}^{1+\kappa} \right) &
\end{cases}\]
where $C,b>1$ and $\kappa, \alpha >0$ are given numbers. If
\[
Y_{0} + Z_{0}^{1+\kappa} \leq
(2C)^{-\frac{1+\kappa}{\sigma}}b^{-\frac{1+\kappa}{\sigma^2}}, \quad
\text{where } \sigma = \min\{\kappa, \alpha\},
\]
then $\{ Y_n \}$ and $\{ Z_n \}$ tend to zero as $n \to \infty$.
\end{lemma}

The following lemma is introduced in \cite{DBGiVe06}; it states that
if the set where $v$ is bounded away from zero occupies a sizable
portion of $K_{\rho}$, then the set where $v$ is positive cluster
about at least one point of $K_{\rho}$. Here we name the inequality
as the isoperimetric inequality.

\begin{lemma}\label{Iso} (Isoperimetric inequality)
Let $v \in W^{1,1}(K_{\rho}^{x_0}) \cap C(K_{\rho}^{x_0})$ for some
$\rho > 0$ and some $x_0 \in \mathbb{R}^{d}$ and let $k$ and $l$ be
any pair of real numbers wuch that $k < l$. Then there exists a
constant $\gamma$ depending upon $N, p$ and independent of $k, l, v,
x_0, \rho$, such that
\[
(l-k) |K_{\rho}^{x_0} \cap \{v > l\}| \leq \gamma
\frac{\rho^{d+1}}{|K_{\rho}^{x_0} \cap \{v \leq k\}|}
\int_{K_{\rho}^{x_0}\cap \{k< v < l\}} |Dv|\,dx .
\]
\end{lemma}

We consider the following heat equation:
\begin{equation}\label{heat-eqn}
v_t-\Delta v =f,\qquad \mbox{in }\, Q_T:=\mathbb{R}^3\times (0,T),
\end{equation}
with initial data $v(x,0)=v_0(x)$. Next, we recall maximal estimates
of the heat equation in terms of mixed norms.
\begin{lemma}\label{lem1}
Let $1<l,m<\infty$. Suppose that $f\in L^{l,m}_{x,t}(Q_T)$ and
$v_0\in {W^{2,l}(\mathbb{R}^3)}$. If $v$ is the solution of the heat
equation \eqref{heat-eqn}, then the following estimate is satisfied:
\begin{equation}\label{heat-estimate}
\norm{v_t}_{L^{l,m}_{x,t}(Q_T)}+\norm{\nabla^2
v}_{L^{l,m}_{x,t}(Q_T)} \leq C
\bke{\norm{f}_{L^{l,m}_{x,t}(Q_T)}+\norm{v_0}_{{W^{2,l}(\mathbb{R}^3)}}}.
\end{equation}
\end{lemma}
We also recall the following Stokes system, which is the linearized
Stokes equations:
\begin{equation}\label{stokes-eqn}
v_t-\Delta v +\nabla p =f,\qquad \text{div} \ v=0 \qquad \mbox{in
}\, Q_T:=\mathbb{R}^3\times (0,T),
\end{equation}
with initial data $v(x,0)=v_0(x)$. The  maximal estimates of the
Stokes system is given as follows.
\begin{lemma}\label{lem2}
Let $1<l,m<\infty$. Suppose that $f\in L^{l,m}_{x,t}(Q_T)$ and
$v_0\in {W^{2,l}(\mathbb{R}^3)}$. If $v$ is the solution of the
Stokes system \eqref{stokes-eqn}, then the following estimate is
satisfied:
\begin{equation}\label{heat-estimate2}
\norm{v_t}_{L^{l,m}_{x,t}(Q_T)}+\norm{\nabla^2
v}_{L^{l,m}_{x,t}(Q_T)} + \norm{\nabla p}_{L^{l,m}_{x,t}(Q_T)} \leq
C
\bke{\norm{f}_{L^{l,m}_{x,t}(Q_T)}+\norm{v_0}_{{W^{2,l}(\mathbb{R}^3)}}}.
\end{equation}
\end{lemma}


\section{Proof of Theorem \ref{T:Holder1}}\label{S:ProofHolder}

For notational convention, we take $\nu_0$ to be the constant from
Proposition~\ref{P:DG} (DeGiorgi type of iteration) corresponding
to $\theta = 1$ and, with $\omega$ and $R$ given positive constants,
we set
\begin{equation}\label{Delta}
\Delta = \left(\frac{\omega}{2}\right)^{-\alpha} \left(2R\right)^2 .
\end{equation}

Our first alternative is that, if a bounded weak solution $n$ stays
close to its maximum on most of one suitable small subcylinder, then
$n$ is away from its minimum on a suitable subcylinder centered at
$(0,0)$, the target point.

Moreover, for given constants $k$ (in analysis it denotes the level
of solution) and $\rho$ (usually it means the spacial radius), we
assume that
\begin{equation}\label{k-rho}
k^{-\alpha - \frac{2\alpha (1+\kappa)}{q_2}} \rho^{d\kappa} < 1 .
\end{equation}
If \eqref{k-rho} fails, we have $k \leq \rho^{\epsilon}$ for some
$\epsilon >0$ which directly implies the H\"{o}lder continuity of a
solution.

\begin{lemma}\label{L:FA} (The first alternative)
Let $\theta_0 > 1$ be a given constant, $\nu_0$ be a constant in
Proposition~\ref{P:DG} (when $\theta = 1$) and $\Delta$ be in
\eqref{Delta}. Suppose $n$ is a nonnegative bounded weak solution of
\eqref{E:PM} in
\begin{equation}\label{Q}
Q = K_{2R} \times (-\theta_0 \Delta, 0)
\end{equation}
with $\alpha \geq 0$. If there is a constant $T_0 \in [-\theta_0
\Delta, -\Delta]$ such that
\begin{equation}\label{FAa}
\left| K_{2R} \times (T_0, T_0 + \Delta) \cap \{ n < \mu_{-} +
\frac{\omega}{2}\}\right| \leq \nu_0 |K_{2R}| \Delta,
\end{equation}
then there is a constant $\delta_1 \in (0,1)$ determined only by
$\theta_0$ and data such that
\[
\essinf_{Q'} \, n(x,t) \geq \mu_{-} + \delta_1 \omega
\]
where
\begin{equation}\label{Q'}
Q' = K_{R/2} \times [-\omega^{-\alpha}\left(\frac{R}{2}\right)^2,
0].
\end{equation}
\end{lemma}

When the assumption \eqref{FAa} fails, which means that $n$ stays
somewhat close to its maximum on a suitable fraction of all suitable
small subcylinders, then eventually $n$ is away from its supremum on
a suitable subcylinder centered at $(0,0)$.

\begin{lemma}\label{L:SA} (The second alternative)
There are constants $\theta_{0}>1$ and $\nu_0 \in (0,1)$ determined
only by data such that, if $n$ is a bounded weak solution of
\eqref{E:PM} in $Q$ (given by \eqref{Q}) with $\alpha \geq 0$ and
\begin{equation}\label{SAa}
\left| K_{2R} \times (T_0, T_0 + \Delta) \cap \{ n > \mu_{+} -
\frac{\omega}{2}\}\right| \leq (1-\nu_0) |K_{2R}| \Delta,
\end{equation}
for all $T_0 \in [-\theta_0 \Delta, -\Delta]$, then there is a
constant $\delta_2 \in (0, 1)$ determined only by data, such that
\[
\esssup_{Q'} \, n \leq \mu_{+} - \delta_2 \omega
\]
where $Q'$ is given in \eqref{Q'}.
\end{lemma}

Also we prove this lemma in Section~\ref{S:ProofAlternatives}. From two alternatives, we infer a decay estimate for the oscillation of a bounded nonnegative weak  solution of \eqref{E:PM}. 

\begin{lemma}\label{L:Main} (Main Lemma)
Let $\alpha,$ $\rho,$ $\mu_{+}$, $\omega$ be given constants with
$\alpha > 0$. Suppose also that $n$ is a nonnegative bounded weak solution of
\eqref{E:PM} in $Q_{\omega, \rho} = K_{\rho} \times
[-\omega^{-\alpha}\rho^2, 0]$ with $ \essosc_{Q_{\omega, \rho}} \, n
\leq \omega.$ Then there are positive constants $\eta$ and
$\lambda$, both less than one and determined only by data such that
\[
\essosc_{Q_{\eta\omega, \lambda\rho}} \, n \leq \eta\omega
\]
where
\[
Q_{\eta\omega, \lambda\rho} = K_{\lambda \rho} \times
[-(\eta\omega)^{-\alpha}(\lambda\rho)^2, 0].
\]
\end{lemma}

\begin{proof}
Let us call $Q$ and $Q'$ where Lemma~\ref{L:FA} and Lemma~\ref{L:SA}
hold. Then, here, our goal is to choose proper $R$ and $\lambda$
such that
\[
Q_{\eta \omega, \lambda \rho} \subseteq Q' \subset Q \subseteq
Q_{\omega, \rho} .
\]
Then the proper relationship of two essential oscillation is
following rather directly from two alternatives, Lemma~\ref{L:FA}
and Lemma~\ref{L:SA}.

For $\theta_0 >1$ given in both Lemma~\ref{L:FA} and
Lemma~\ref{L:SA}, we choose
\[
R = \theta_{0}^{-1/2} 2^{-\alpha/2 -1} \rho < \rho/2.
\]
Then $Q \subseteq Q_{\omega, \rho}$. For $\eta = \max\{1-\delta_1,
1-\delta_2\} \in (0,1)$ where $\delta_1$ and $\delta_2$ are from
Lemma~\ref{L:FA} and Lemma~\ref{L:SA}, we choose
\begin{equation}\label{lambda}
\lambda = \eta^{\alpha/2} \theta_{0}^{-1/2} 2^{-\alpha/2 - 2}
\end{equation}
so that $Q_{\eta \omega, \lambda \rho} \subseteq Q'$.

If there is a $T_0 \in (-\theta_0 \Delta, -\Delta)$ such that
\begin{equation}\label{M01}
\left| K_{2R} \times (T_0, T_0 + \Delta) \cap \{ n > \mu_{+} -
\frac{\omega}{2}\}\right| \leq \nu_0 |K_{2R}| \Delta,
\end{equation}
then by Lemma~\ref{L:FA} applied to $n$ implies that
\[
\esssup_{Q'} n \leq \mu_{+} - \delta_1 \omega.
\]
Hence, by Lemma~\ref{L:FA} and $Q_{\eta\omega, \lambda\rho}
\subseteq Q'$, it follows that
\[
\essosc_{Q_{\eta\omega, \lambda\rho}} n  \leq \essosc_{Q'} n \leq
\eta\omega .
\]

When \eqref{M01} fails, then it holds
\[
\left| K_{2R} \times (T_0, T_0 + \Delta) \cap \{ n < \mu_{-} +
\frac{\omega}{2}\}\right| \leq (1-\nu_0) |K_{2R}| \Delta.
\]
By Lemma~\ref{L:SA}, we have
\[
\essinf_{Q'} n \geq \mu_{-} + \delta_2 \omega
\]
which implies that
\[
\essosc_{Q_{\eta\omega, \lambda\rho}} n \leq \eta\omega.
\]
This completes the proof.
\end{proof}

Now we provide the proof of Theorem~\ref{T:Holder1}.
\smallskip

\textbf{Proof of Theorem~\ref{T:Holder1}}: If $\omega=0$, then this
result is true for any choice of $\gamma$ and $\beta$, so we assume
that $\omega > 0$ and set $\omega_0 = \omega$. We also set
\[
\rho_0 = \dist_{p} \left( \{(x_1, t_1)\},
\partial_{p}\Omega_{T}\right).
\]
We define
\[
\omega_{i} = \eta^{i} \omega_0, \quad \rho_{i} = \lambda^{i}\rho_0,
\]
where $\lambda$ and $\eta$ are the constants from
Lemma~\ref{L:Main}. Also define a sequence of cylinders $\{Q_{n}\}$
by
\[
Q_{i} = K^{x_1}_{\rho_{i}} \times [ t_1 - \omega_{i}^{-\alpha}
\rho_{i}^{2}, t_1].
\]
It is easy to check that $Q_0 \subset \Omega$ and that $Q_{i+1}
\subset Q_{i}$ for any $i$. Combining with Lemma~\ref{L:Main} with
an induction argument, we find that $\essosc_{Q_i} n \leq \omega_i$
for any $i$.

For $(x_2, t_2) \in Q_0$ with $x_1 \neq x_2$ and $t_1 \neq t_2$,
then there are nonnegative integers $k$ and $l$ such that
\begin{equation}\label{Holx}
\rho_{k+1} < |x_1 - x_2| \leq \rho_{k},
\end{equation}
and
\begin{equation}\label{Holt}
\omega_{l+1}^{-\alpha}\rho_{l+1}^{2} < |t_1 - t_2| \le
\omega_{l}^{-\alpha}\rho_{l}^{2}.
\end{equation}
As a result, we obtain that
\[
|n(x_1, t_1) - n(x_2, t_2)| \leq \max\{ \omega_{k}, \omega_{l}\}.
\]
From \eqref{Holx} , we derive, for $\beta_1 = \log_{\eta} \lambda$,
\[
\frac{|x_1 - x_2|}{\rho_0} > \lambda^{k+1} =
\left(\eta^{\beta_1}\right)^{k+1}
\]
which implies
\[
\omega_{k}= \eta^{k}\omega_0 < \eta \omega_0 \left(\frac{|x_1 -
x_2|}{\rho_0}\right)^{\beta_1}.
\]
On the other hand, the inequality \eqref{Holt} implies that
\[
|t_1 - t_2|_{I} > \eta^{-\alpha(l+1)/2} \rho_{l+1} =
\left(\eta^{-\alpha/2}\lambda\right)^{l+1}\rho_0.
\]
Because of the choice of $\lambda$ from \eqref{lambda}, let us
denote that
\[
\tilde{\lambda} = \eta^{-\alpha/2}\lambda < 1.
\]
Then we have
\[
\omega_{l} = {\tilde{\lambda}}^{\beta_2 l}\omega_{0} \leq \left(
\frac{|t_1 - t_2|_{I}}{\tilde{\lambda}\rho_0}\right)^{\beta_2}
\omega_{0}
\]
for $\beta_2 = \log_{\tilde{\lambda}}\eta$.

Therefore, for some $\gamma > 0$,
\[
|n(x_1, t_1) - n(x_2, t_2)| \leq \gamma \omega \left[
\left(\frac{|x_1 - x_2|}{\rho_0}\right)^{\beta_1}+ \left( \frac{|t_1
- t_2|_{I}}{\rho_0}\right)^{\beta_2} \right].
\]
Then this implies \eqref{Holder} with $\beta = \min\{\beta_1,
\beta_2\}$ and the definition of $|\cdot|_{I}$ \eqref{tI} because
$\rho_0 \geq \dist_{p}(\Omega'; \partial_{p}\Omega)$.

If $x_1 = x_2$ or $t_1 = t_2$, then a similar and simpler arguments
yields the same result.


\begin{remark}\label{R:gen}
Here we make comments that our method of analysis to show the local H\"{o}lder continuity of \eqref{E:PM} can easily modified to expalin the same regularity for a generalized structured equation. Now consider porous medium equation in the form of
\begin{equation}\label{E:GPM}
n_{t} = \Delta \Phi(n) + \nabla \cdot \left( B(x,t) n\right)
\end{equation}
where $\Phi \in C^{1}[0,\infty)$. Let $\Phi'(s) = \phi(s)$ where $\phi$ is a nonnegative increasing function with $\phi(0)=0$ and we assume that there are two constants $\alpha_0$ and $\alpha_1$ satisfying $0\leq \alpha_0 \leq \alpha_1 < \infty$ such that
\[
(1+\alpha_0) \Phi(s) \leq s \phi(s) \leq (1+\alpha_1) \Phi(s)
\]
for all $s>0$. The two inequalities are essentially the $\nabla_2$ and $\Delta_2$ conditions in Orlicz space theory. If $\alpha_0 = \alpha_1$, then \eqref{E:GPM} is exactly \eqref{E:PM}.

The local and logarithmic energy estimates are obtained for \eqref{E:GPM} by replcing $n^{\alpha}$ to $\phi(n)$ in the estimates from Propositions~\ref{P:EE} and \ref{P:LEE}. Hence corresponding intrinsic scaling for \eqref{E:GPM} follows immediately as $ T_{\omega, \rho} = \theta \rho^2 / \phi(\omega)$ (cf. \eqref{intrinsic}). Moreover, we assume
$$ \phi(k)^{-1 - \frac{2(1+\kappa)}{q_2}} \rho^{d\kappa} < 1 $$
instead of \eqref{k-rho}. Then the same proofs in Sections~\ref{S:ProofHolder} and \ref{S:ProofAlternatives} hold for \eqref{E:GPM} as well.
\end{remark}



\section{Proofs of the two alternatives}\label{S:ProofAlternatives}

In this section, we deliver the proofs of two alternatives,
Lemmas~\ref{L:FA} and \ref{L:SA}. The proofs in
Section~\ref{SS:ProofAlternatives} are basically composed with two
parts, DeGiorgi type iteration and the expansion of positivities
along the time and space variables.

\subsection{Local energy estimates}

In this section, we provide two types of local energy estimates that
are key prove the modulus of continuity of $n$. We make two remarks.
First, to carry calculations directly with weak solutions rather
sub(super-)solutions, we take advantage that both $n-\mu_{-}$ and
$\mu_{+} - n$ are nonnegative that leads the positiveness of level
$k$. Second, because the lower order term does not given in the form
of $n^{1+\alpha}$ in \eqref{E:PM}, we assume proper conditions on
$B$ and $\nabla B$ (cf. Chapter~{3.6} of \cite{DBGiVe12}).

For given positive constants $k$ and $\rho$, we denote the set
\begin{equation}\label{Akrho}
A_{k,\rho}^{\pm}(\tau) = \{ x\in K_{\rho}: \, \left(n(x,\tau) -
\mu_{\pm} \pm k \right)_{\pm} > 0 \},
\end{equation}
that indicates a level set (either $n < \mu_{-} + k$ or $n > \mu_{+}
- k$) at a fixed time $\tau$.

\begin{proposition}\label{P:EE}
Suppose that $\zeta$ is a cutoff function on the parabolic cylinder
$Q_{\rho} = K_{\rho} \times [t_0, t_1]$, vanishing on the parabolic
boundary of $Q_{\rho}$ with $0 \leq \zeta \leq 1$. For a nonnegative
bounded weak solution $n$ of \eqref{E:PM} under \eqref{B} , it
follows, for any $k >0$ and for some positive constants
\begin{equation}\begin{split}\label{EE}
  & \sup_{t_0 \leq t \leq t_1} \, \int_{K_{\rho}\times\{t\}} (n- \mu_{\pm} \pm k)_{\pm}^2 \zeta^2 \,dx
  + (1+\alpha)\iint_{Q_{\rho}} \, n^{\alpha} |\nabla (n-\mu_{\pm} \pm k)_{\pm}|^2 \, \zeta^2 \,dx\,dt \\
    &\quad\leq \int_{K_{\rho}\times\{t_0\}} (n- \mu_{\pm} \pm k)_{\pm}^2 \zeta^2 \,dx + 2 \iint_{Q_\rho}(n- \mu_{\pm} \pm k)_{\pm}^2 \zeta \zeta_{t} \,dx\,dt\\
    &\qquad + 16(1+\alpha) \iint_{Q_{\rho}}  \, (n^{\alpha}+ k^{\alpha}) (n-\mu_{\pm}\pm k)^{2}_{\pm} \, |\nabla \zeta|^2 \,dx\,dt \\
  &\qquad + \left(k^{2-\alpha} \|B\|^{2}_{2 \hat{q}_1, 2\hat{q}_2} + k^2 \|\nabla B\|_{\hat{q}_1, \hat{q}_2}\right)
  \left[ \int_{t_0}^{t_1} \left[A^{\pm}_{k,\rho}(t)\right]^{\frac{q_2}{q_1}} \,dt\right]^{\frac{2(1+\kappa)}{q_2}}.
\end{split}\end{equation}
\end{proposition}

We deliver the proof in Section~\ref{SS:ProofEnergy}.


Now we provide a logarithmic energy estimate that is crucial to
capture the expansion of positivity along the time axis,
Proposition~\ref{P:EPT}.

For given positive constants $k$ and $\delta$, let us define
\begin{equation}\label{Psi}
\Psi_{\pm}(n) = \ln^{+} \left[ \frac{k}{(1+\delta)k - (n - \mu_{\pm}
\pm k)_{\pm}} \right].
\end{equation}
By assuming proper conditions on $B$ and $\nabla B$ given in
\eqref{B}, we succeed to calculate logarithmic estimates (refer
Section~B.7 on \cite{DBGiVe12}).

\begin{proposition}\label{P:LEE}
For $n$ a nonnegative weak solution of \eqref{E:PM} with $\alpha
\geq 0$ under \eqref{B}, suppose that $\zeta$ is a cutoff function
in $K_{\rho} \times [t_0, t_1]$ which is vanishing on the lateral
boundary of $K_{\rho} \times [t_0, t_1]$ (independent of the time
variable) with $0 \leq \zeta \leq 1$. For given positive constants
$k$ and $\delta$, it follows that
\begin{equation}\label{LEE}\begin{split}
& \int_{K_{\rho} \times \{t_1\}} \Psi^{2}_{\pm}(n) \zeta^2 \,dx + 2(1+\alpha)\int_{t_0}^{t_1}\int_{K_{\rho}} n^{\alpha} |\nabla n|^2  \Psi_{\pm}(n) (\Psi'_{\pm}(n))^2 \zeta^2  \,dx\,dt \\
&\leq \int_{K_{\rho} \times \{t_0\}} \Psi^{2}_{\pm}(n) \zeta^2 \,dx + 16 (1+\alpha)\int_{t_0}^{t_1}\int_{K_{\rho}} (n^{\alpha}+k^{\alpha}) \Psi^{2}_{\pm}(n) |\nabla \zeta|^2  \,dx\,dt \\
&\quad +  2 k^{-\alpha}\left(\ln \frac{1}{\delta}\right)^2 \|B\|^{2}_{2 \hat{q}_1, 2\hat{q}_2}\left[ \int_{t_0}^{t_1} \left[A^{\pm}_{k,\rho}(t)\right]^{\frac{q_2}{q_1}} \,dt\right]^{\frac{2(1+\kappa)}{q_2}} \\
&\quad + \left( (\ln \frac{1}{\delta})^2 + \frac{\mu_{+}\ln
\frac{1}{\delta}}{\delta k}\right) \|\nabla B\|_{q_1, q_2} \left[
\int_{t_0}^{t_1} \left[A^{\pm}_{k,\rho}(t)\right]^{\frac{q_2}{q_1}}
\,dt\right]^{\frac{2(1+\kappa)}{q_2}}.
\end{split}\end{equation}
\end{proposition}
We deliver the proof in Section~\ref{SS:ProofEnergy}.

\subsection{DeGiorgi type of estimates}\label{SS:DG}

In this section, we modify DeGiorgi iteration for the porous medium
equations given in \eqref{E:PM} with $\alpha \geq 0$ that starts from
the local energy estimates given in Proposition~\ref{P:EE}. (refer Lemma~{7.1} of
\cite{DBGiVe12})

\begin{proposition}\label{P:DG} (DeGiorgi iteration)
Let $n$ be a bounded nonnegative weak solution of \eqref{E:PM} under
\eqref{B} with $\alpha \geq 0$. For given positive constants $\rho$,
$k$, and $\theta$ satisfying \eqref{k-rho}, let
\[
Q_{\rho} = Q (\rho, \theta k^{-\alpha} \rho^2) = K_{\rho} \times
[-\theta k^{-\alpha} \rho^2, 0].
\]
\begin{enumerate}

\item[(i)] There exists $\nu_0 = \nu_0 (\theta, \text{data}) \in (0,1)$ such that, if
\begin{equation}\label{DGa_sub}
\left|Q_{\rho} \cap \{ n > \mu_{+} - k\}\right| \leq \nu_0 \left|
Q_{\rho}\right|,
\end{equation}
then it holds
\begin{equation}\label{DGb_sub}
\esssup_{Q_{\rho/2}} n(x,t) \leq \mu_{+} - \frac{k}{2}.
\end{equation}

\item[(ii)] There exists $\nu_0 = \nu_0 (\theta, \text{data}) \in (0,1)$ such that, if
\begin{equation}\label{DGa_super}
\left|Q_{\rho} \cap \{ n < \mu_{-} + k\}\right| \leq \nu_0 \left|
Q_{\rho}\right|,
\end{equation}
then it holds
\begin{equation}\label{DGb_super}
\essinf_{Q_{\rho/2}} n(x,t) \geq \mu_{-} + k/2 .
\end{equation}

\end{enumerate}
\end{proposition}

\begin{proof}
First, we consider a level set $(n - \mu_{+} + k)_{+}$ for a
nonnegative bounded weak solution $n$ of \eqref{E:PM}. We constructs
sequences $\{\rho_i\}$, $\{k_i\}$, $\{K_{i}\}$, and $\{Q_i\}$ such
that
\begin{gather}\label{DGrho}
\rho_{i} = \frac{\rho}{2} + \frac{\rho}{2^{i+1}}, \quad \rho_0 = \rho, \ \rho_{\infty}= \frac{\rho}{2},\\
k_{i} =\frac{k}{2} + \frac{k}{2^{i+1}}, \quad k_{0}=k, \ k_{\infty} =\frac{k}{2},\\
K_{i}= K_{\rho_i}, \quad Q_{i} =K_{i}\times [-\theta k^{-\alpha}
\rho_{i}^{2}, 0].
\end{gather}
Moreover, we take a sequence of piecewise linear cutoff functions
$\{ \zeta_{i}\}_{i=0}^{\infty}$ such that
\[
\zeta_{i} = \begin{cases}
            1 & \text{ in } Q_{i+1} \\
            0 & \text{ on the parabolic boundary of } Q_{i},
 \end{cases}
\]
satisfying
\[
|\nabla \zeta_{i}| \leq \frac{2^{i+2}}{\rho}, \quad
\partial_{t} \zeta_{i} \leq \frac{2^{i+2} k^{\alpha}}{\theta \rho^{2}}.
\]

It is easy to observe $(n - \mu_{+} + k_i)_{+} \leq k_i \leq k $ and
$ n > \mu_{+} - k_{i} > \gamma k $ by choosing $k > (\gamma +
1/2)^{-1}\mu_{+}$ for some constant $\gamma >0$.

Then the energy estimate given in Proposition~\ref{P:EE} provides
\begin{equation}\begin{split}\label{EE01}
  & \sup_{-\theta k^{-\alpha}\rho_{i}^2 < t < 0} \, \int_{K_i \times\{t\}} (n- \mu_{+} + k_{i})_{+}^2 \zeta_{i}^{2} \,dx \\
  &\quad  + C(\alpha, \gamma) k^{\alpha}\iint_{Q_i} \, |\nabla (n- \mu_{+} + k_{i})_{+}|^2 \, \zeta_{i}^{2} \,dx\,dt \\
  &\leq \left[ \frac{2^{i+3}}{\theta} + 16(1+\alpha)\gamma^{\alpha}\right] \frac{k^{2+\alpha}}{\rho^2} \left| Q_i \cap \{ (n-\mu_{+} + k_{i})_{+} > 0\} \right| \\
  &\quad + C(k^{2-\alpha} + k^2) \left[ \int_{-\theta k^{-\alpha}\rho_{i}^2}^{0} \left[A^{+}_{k_{i},\rho_{i}}(t)\right]^{\frac{q_2}{q_1}} \,dt\right]^{\frac{2(1+\kappa)}{q_2}}.
  \end{split} \end{equation}

Now we take the change of variable that $ \bar{t} = k^{\alpha} t \in
[-\theta \rho^2, 0] .$ Also denote $\bar{n} = n(\cdot, \bar{t})$,
$\bar{\zeta} = \zeta(\cdot, \bar{t})$, and $\overline{Q}_{i} = K_{i}
\times [-\theta \rho_{i}^{2}, 0]$. Then \eqref{EE01} gives
\begin{equation}\begin{split}\label{EE02}
  & \sup_{-\theta \rho_{i}^2 < \bar{t} < 0} \, \int_{K_i \times\{\bar{t}\}} (\bar{n}- \mu_{+} + k_{i})_{+}^2 \bar{\zeta}_{i}^{2} \,dx \\
  &\quad  + C(\alpha, \gamma) \iint_{\overline{Q}_i} \, |\nabla (\bar{n}- \mu_{+} + k_{i})_{+}|^2 \, \bar{\zeta}_{i}^{2} \,dx\,d\bar{t} \\
  &\leq \left[ \frac{2^{i+3}}{\theta} + 16(1+\alpha)\gamma^{\alpha}\right] \frac{k^{2}}{\rho^2} \left| \overline{Q}_i \cap \{ (\bar{n}-\mu_{+} + k_{i})_{+} > 0\} \right| \\
  &\quad + C(k^{2-\alpha} + k^2)k^{-\frac{2\alpha (1+\kappa)}{q_2}} \left[ \int_{-\theta \rho_{i}^2}^{0} \left[A^{+}_{k_{i},\rho_{i}}(\bar{t})\right]^{\frac{q_2}{q_1}} \,d\bar{t}\right]^{\frac{2(1+\kappa)}{q_2}}.
  \end{split} \end{equation}

For simplicity, denote two sets
\begin{align*}
A_{i} &= Q_{i} \cap \{ (n-\mu_{+} + k_{i})_{+} > 0\}, \\
\overline{A}_{i} &= \overline{Q}_{i} \cap \{ (\bar{n}-\mu_{+} +
k_{i})_{+} > 0\}.
\end{align*}

To handle the left hand side of \eqref{EE02}, we apply Sobolev
embedding theorem, Theorem~\ref{T:Sobolev}, from which we calculate
\begin{equation}\label{EE03}\begin{split}
& \iint_{\overline{Q}_{i+1}}  (\bar{n}- \mu_{+} + k_{i})_{+}^2  \,dx \,d\bar{t}
\leq \iint_{\overline{Q}_{i}}  (\bar{n}- \mu_{+} + k_{i})_{+}^2 \bar{\zeta}_{i}^{2} \,dx \,d\bar{t} \\
&\leq |\overline{A}_{i}|^{\frac{2}{d+2}}\left[
\sup_{-\theta \rho_{i}^2 < \bar{t} < 0} \, \int_{K_i \times\{\bar{t}\}} (\bar{n}- \mu_{+} + k_{i})_{+}^2 \bar{\zeta}_{i}^{2} \,dx \right.\\
&\quad  \left. + C(\alpha, \gamma) \iint_{\overline{Q}_i} \, \nabla
\left[(\bar{n}- \mu_{+} + k_{i})_{+} \, \bar{\zeta}_{i} \right]^{2}
\,dx\,d\bar{t} \right].
\end{split}\end{equation}

In the set $\{ (\bar{n} - \mu_{+} + k_{i+1})_{+} > 0 \}$, we observe
that
$$ (\bar{n} - \mu_{+} + k_{i})_{+} \geq k_{i} - k_{i+1} = \frac{k}{2^{i+2}}. $$
Then the combination of \eqref{EE02} and \eqref{EE03}, carrying
cancellation on $k^2$, we are able to say that
\begin{equation}\label{EE04}\begin{split}
|\overline{A}_{i+1}|
&\leq C(\alpha, d, \gamma, \theta)\frac{2^{3i}}{\rho^2} |\overline{A}_{i}|^{1+\frac{2}{d+2}}  \\
  &\quad + C(k^{-\alpha} + 1)k^{-\frac{2 \alpha  (1+\kappa)}{q_2}} \left[ \int_{-\theta \rho_{i}^2}^{0} \left[A^{+}_{k_{i},\rho_{i}}(\bar{t})\right]^{\frac{q_2}{q_1}} \,d\bar{t}\right]^{\frac{2(1+\kappa)}{q_2}}|\overline{A}_{i}|^{\frac{2}{d+2}} .
\end{split}\end{equation}

We note that
$$ \rho^2 \sim |\overline{Q}_{i}|^{\frac{2}{d+2}}, \quad |K_{i}| \sim |\overline{Q}_{i}|^{\frac{d}{d+2}}. $$
Let
$$ \overline{Z}_{i} = \frac{1}{|K_{i}|} \left[ \int_{-\theta \rho_{i}^2}^{0} \left[A^{+}_{k_{i},\rho_{i}}(\bar{t})\right]^{\frac{q_2}{q_1}} \,d\bar{t}\right]^{\frac{2}{q_2}}.$$

By dividing \eqref{EE04} by
$|\overline{Q}_{i}|=c(d)|\overline{Q}_{i+1}|$, we obtain (by
applying \eqref{k-rho})
\begin{equation}\label{EE05}
\frac{|\overline{A}_{i+1}|}{|\overline{Q}_{i+1}|} \leq C(\alpha, d,
\gamma, \theta)2^{3i}
\left\{\frac{|\overline{A}_{i}|}{|\overline{Q}_{i}|}\right\}^{1+\frac{2}{d+2}}
+ C \bar{Z}_{i}^{1+\kappa}
\left\{\frac{|\overline{A}_{i}|}{|\overline{Q}_{i}|}\right\}^{\frac{2}{d+2}}.
\end{equation}
Therefore, we take change variable back to $t$ from $\bar{t}$ from
the dimensionless inequality \eqref{EE05}, by letting $ Y_{i} =
|A_i|/|Q_i|, $ we obtain the following inequality
\begin{equation}
Y_{i+1} \leq C_0 2^{3i} Y_{i}^{1+ \frac{2}{d+2}} + C_1 2^{2i}
Z_{i}^{1+\kappa}Y_{i}^{\frac{2}{d+2}}.
\end{equation}
Then we are able to apply Lemma~\ref{L:iter2} that there exists
$$ \nu_0 \leq (2C)^{-\frac{1+\kappa}{\sigma}}2^{-\frac{3(1+\kappa)}{\sigma^2}}$$
where $C = \max\{C_0, C_1\}$ and $\sigma = \min\{\kappa,
\frac{2}{d+2}\}$ such that $Y_{i}$ and $Z_i$ converge to $0$ as $i
\to \infty$. Hence, $n$ is greater than $\mu_{+}- \frac{k}{2}$ in
almost everywhere of the set $Q_{\rho/2}$.

\bigskip

Second, we now carry the DeGiorgi iteration with the level sets $(n
- \mu_{-} - k)_{-}$. It is easy to see that $ 0 \leq (n - \mu_{-} -
k_{i})_{-} < k_{i} \leq k .  $ We wish to avoid when $\mu_{-}$ is
near zero so it is hard to estimate the lower bound of the local
energy estimate \eqref{EE}. Therefore, we introduce
\begin{equation}\label{Ntil}
  m = \max\left\{ n, \ \mu_{-} + \frac{k}{4} \right\}.
\end{equation}
which provides that
\begin{equation}\label{mn01}
\iint_{Q_i} \, m^{\alpha}|\nabla (m- \mu_{-} - k_{i})_{-}|^2 \,
\zeta_{i}^{2} \,dx\,dt \leq \iint_{Q_i} \, n^{\alpha}|\nabla (n-
\mu_{-} - k_{i})_{-}|^2 \, \zeta_{i}^{2} \,dx\,dt
\end{equation}
considering two sets where $\{Q_{i}: \, m= n\}$ and $\{ Q_{i}: \, m
= \mu_{-} + \frac{k}{4}\}$ where $|\nabla (m- \mu_{-} - k_{i})_{-}|
= 0$ in the latter set. Moreover, we compute that
\begin{equation}\label{mn02}\begin{split}
& \sup_{-\theta k^{-\alpha}\rho_{i}^2 < t < 0} \, \int_{K_i \times\{t\}} (m- \mu_{-} - k_{i})_{-}^2 \zeta_{i}^{2} \,dx \\
&\leq \sup_{-\theta k^{-\alpha}\rho_{i}^2 < t < 0} \, \int_{K_i \times\{t\}} (n- \mu_{-} - k_{i})_{-}^2 \zeta_{i}^{2} \,dx
 + \frac{\theta}{4^2} \theta k^{2+\alpha}\rho_{i}^{-2} \left|
Q_i \cap \{ (n-\mu_{-}-k_i)_{-} > 0 \} \right|.
\end{split}\end{equation}

Then the combination of two inequalities \eqref{mn01} and
\eqref{mn02} provides the energy estimates \eqref{EE01} in terms of
$m$ on the left-hand-side. By taking DeGiorgi iteration, it provides
the existence of $\nu_0$ such that if $|Q_{\rho} \cap \{ n < \mu_{-}
+ k\}| \leq \nu_0 |Q_{\rho}|$, the it holds
\[
\essinf_{Q_{\rho/2}} m \geq \mu_{-} + \frac{k}{2},
\]
which leads our conclusion.
\end{proof}


Next proposition is a variant of DeGiorgi iteration using the
information at a certain fixed time level to obtain the same
conclusion as in Proposition~\ref{P:DG}(in this poroposition $\nu_0$
is depending on both data and $\theta$) where $\nu^{*}$ depending
only on data.

\begin{proposition}\label{P:DG02}
Let $n$ be a bounded nonnegative weak solution of \eqref{E:PM} under
\eqref{B} with $\alpha \geq 0$. For given positive constants $\rho$,
$k$, and $\theta$ satisfying \eqref{k-rho}, let $Q_{\rho}$ and $Q_{\rho/2}$ be given in Proposition~\ref{P:DG}.
\begin{enumerate}

\item[(i)] There exists $\nu^{*} \in (0,1)$ determined only by data, such that, if
$$ n(x, -\theta k^{-\alpha} \rho^2) \leq \mu_{+} - k, $$ and if
\begin{equation}\label{DGa_sub2}
\left|Q_{\rho} \cap \{ n > \mu_{+} - k\}\right| \leq
\frac{\nu^{*}}{\theta} \left| Q_{\rho}\right|,
\end{equation}
then it holds
\begin{equation}\label{DGb_sub2}
\esssup_{Q_{\rho/2}} n(x,t) \leq \mu_{+} - \frac{k}{2}.
\end{equation}

\item[(ii)] There exists $\nu^{*} \in (0,1)$ determined only by data, such that, if
$$ n(x, -\theta k^{-\alpha} \rho^2) \geq \mu_{-} + k ,$$ and if
\begin{equation}\label{DGa_super2}
\left|Q_{\rho} \cap \{ n < \mu_{-} + k\}\right| \leq \nu_0 \left|
Q_{\rho}\right|,
\end{equation}
then it holds
\begin{equation}\label{DGb_super2}
\essinf_{Q_{\rho/2}} n(x,t) \geq \mu_{-} + k/2 .
\end{equation}

\end{enumerate}
\end{proposition}

The proof is done by repeating the same proof for
Proposition~\ref{P:DG} with $\frac{\nu^{*}}{\theta}$ instead of
$\nu_0$. We refer Proposition~4.5 from \cite{HL15a}.

\subsection{The expansion of positive data}\label{SS:Expansion}


This section is to understand the behavior of a nonnegative bounded
weak solution of \eqref{E:PM} explaining positive data's flows in time and spatial axis sperately in measure sense (so called the expansion of positivities).
The following proposition shows that if a nonnegative function is
large on part of a cylinder, then it keeps largeness on part of a suitable
time slice same as Proposition~{4.1} in \cite{HL15a}.

\begin{proposition}\label{P:prop1}
Let $k$, $\rho$, and $T$ be positive constants. If $n$ is a
measurable nonnegative function defined on $Q= K_{\rho} \times (-T,
0)$ and if there is a constant $\nu_1 \in (0,1)$ such that
\[
\left| Q \cap \{ n > \mu_{+} - k \}\right| \leq (1-\nu_1) |Q|
\]
there there is a number $ \tau_1 \in (-T, - \frac{\nu_1}{2-\nu_1} T)
$ for which
\[
\left| \{ x\in K_{\rho}: n(x, \tau_1) > \mu_{+} - k \}\right| \leq
(1-\frac{\nu_1}{2})|K_{\rho}|.
\]
\end{proposition}

With the aid of logarithmic energy estimate in Proposition~\ref{P:LEE}, we are able to control the measure where a weak solution keeps its largeness in a certain later time when we have according measure information at a fixed time level.
\begin{proposition}\label{P:EPT}
Let $n$ be a nonnegative weak solution of \eqref{E:PM} with $\alpha
\geq 0$ under \eqref{B}. Suppose that positive constants $\rho$,
$k$, and $\mu \in (0,1)$ are given satisfying \eqref{k-rho}.
\begin{enumerate}
 \item[(i)]  Assume that
  \begin{equation}\label{EPTa}
  \left| \left\{ x\in K_{\rho} : n(x, t_0) > \mu_{+} - k \right\} \right| \leq (1-\nu) |K_{\rho}|.
  \end{equation}
  Then for any $\epsilon \in (0,1)$, if
  $ t_1 - t_0 \leq \theta k^{-\alpha} \rho^2, $
  there exist $\delta$ depending on the data, $\nu$, and $\epsilon$ such that
  \begin{equation}\label{EPTc}
  \left| \left\{ x\in K_{\rho} : n(x,t) > \mu_{+} - \delta k \right\}\right| < (1 - (1-\epsilon)\nu)|K_{\rho}|
  \end{equation}
  for any $ t_0 \leq t \leq t_1$.

  \item[(ii)] Assume that
  \begin{equation}\label{EPTaa}
  \left| \left\{ x\in K_{\rho} : n(x, t_0) < \mu_{-} + k \right\} \right| \leq (1-\nu) |K_{\rho}|.
  \end{equation}
  Then for any $\epsilon \in (0,1)$, if
  $ t_1 - t_0 \leq \theta k^{-\alpha} \rho^2, $
  there exist $\delta$ depending on the data, $\nu$, and $\epsilon$ such that
  \begin{equation}\label{EPTcc}
  \left| \left\{ x\in K_{\rho} : n(x,t) < \mu_{-} + \delta k \right\}\right| < (1 - (1-\epsilon)\nu)|K_{\rho}|
  \end{equation}
  for any $ t_0 \leq t \leq t_1$.

\end{enumerate}
\end{proposition}

\begin{proof}
 We apply the logarithmic energy estimates from Proposition~\ref{P:LEE}. Let $\zeta$ be a linear cutoff function independent of the time variable such that
 \[
 \zeta =\begin{cases}
 1 & \text{ in } K_{(1-\sigma)\rho} \times [t_0, t], \\
 0 & \text{ on the lateral boundary of } K_{\rho} \times [t_0, t]
 \end{cases}\]
 for any $t \in (t_0, t_1)$ satisfying
 $|\nabla \zeta| \leq \frac{1}{\sigma \rho}$ and $\zeta_{t}= 0$
 where $\sigma \in (0,1)$ is to be determined later.
Let $\delta = 2^{-j}$ for $j$ is a positive integer which will be
chosen later large enough.
For $\Psi_{\pm}(n)$ in \eqref{Psi}, we observe first, $ (n -
\mu_{\pm} \pm k)_{\pm} \leq k \quad $ which provides
$$ \Psi_{\pm}(n) \leq \ln \frac{1}{\delta} = j \ln 2 \quad
\text{and} \quad |\Psi'_{\pm}(n)| \leq \frac{1}{\delta k}. $$
Moreover, in the set $\{ n > \mu_{+} - \delta k \}$ and $\{ n <
\mu_{-} + \delta k\}$ respectively, we have
$$ (n - \mu_{\pm} \pm k)_{\pm} > (1-\delta)k . $$
In both cases, it gives that
$$ \Psi_{\pm}(n) \geq \ln \frac{1}{2\delta} = (j-1)\ln 2 .$$

From \eqref{LEE} and \eqref{k-rho}, it follows
\begin{equation}\label{EPT01}\begin{split}
& (j-1)^2 (\ln 2)^2 \left| \{ K_{(1-\sigma)\rho}: (n(\cdot, t)-\mu_{\pm}\pm \delta k)_{\pm} > 0 \} \right| \\
&\leq j^2 (\ln 2)^2 \left| \{ K_{\rho} : (n(\cdot, t_0)-\mu_{\pm}\pm k)_{\pm} > 0 \} \right| \\
&\quad + C(\alpha, \gamma) \frac{k^{\alpha} j \ln 2}{\sigma^2\rho^2} |t-t_0| |K_{\rho}| + C(\theta,\|B\|^{2}_{2\hat{q}_1,2\hat{q}_2},\|\nabla
B\|_{\hat{q}_1, \hat{q}_2})(j\ln 2)^2 |K_{\rho}|.
\end{split}\end{equation}

Hence all estimates above yields the following inequality:
\begin{equation}\label{EPT03}\begin{split}
& \left| \{ K_{(1-\sigma)\rho}: (n(\cdot, t)-\mu_{\pm}\pm \delta k)_{\pm} > 0 \} \right| \\
&\leq \left[ \left(\frac{j}{j-1}\right)^2 (1-\nu) + \frac{C(c_0,
\alpha) j}{\sigma^2 (j-1)^2} + C(\delta) \right] |K_{\rho}|
\end{split}\end{equation}
by assumptions. Therefore, it leads that
\begin{equation}\label{EPT05}\begin{split}
& \left| \{ K_{\rho}: (n(\cdot, t)-\mu_{\pm}\pm \delta k)_{\pm} > 0 \} \right| \\
&\leq \left[ \left(\frac{j}{j-1}\right)^2 (1-\nu) + \frac{C(c_0,
\alpha) j}{\sigma^2 (j-1)^2} + C(\delta) + d \sigma\right]
|K_{\rho}|.
\end{split}\end{equation}
Then next we make a choice of $j$ (so $\delta = 2^{-j}$) satisfying
the following inequalities:
\[
\left(\frac{j}{j-1}\right)^2 \leq 1+\epsilon \nu, \quad
\frac{C(c_0, \alpha) j}{\sigma^2 (j-1)^2} \leq \frac{\epsilon \nu^2}{4}, \quad
C(\delta) \leq \frac{\epsilon \nu^2}{4}, \quad  d \sigma \leq
\frac{\epsilon \nu^2}{4},
\]
which provides \eqref{EPTc} and \eqref{EPTcc}. Let us choose
\[
\sigma = \frac{\epsilon \nu^2}{4 d}, \ j = \max\left\{ 1+
\frac{1}{\sqrt{1+\epsilon\nu} - 1}, \ \frac{C 4^3 d^2 (1+\epsilon
\nu)}{\epsilon^3 \nu^6} \right\}.
\]
The inequality for $C(\delta)$ is trivially following.
\end{proof}


The following proposition is obtaining an arbitrary control over the measure of a
cylinder where a weak solution is lager than some constant, if at
each time level we know the measure of the set where a weak solution is somewhat large.

\begin{proposition}\label{P:EPX}
  Let $k$, $\rho$, and $\theta$ be positive constants satisfying \eqref{k-rho}.
  Suppose that $n$ is a nonnegative bounded weak solution of \eqref{E:PM} under \eqref{B} with $\alpha \geq 0$ in $K_{2\rho}\times [-2\tau, 0]$. Then for any $\beta$ and $\nu$ in $(0,1)$, there exists $\delta^{*} = \delta^{*}(\beta, \nu, \theta, \data) \in (0,1)$ depending on data such that, if
  \begin{equation}\label{tau1}
  \tau \geq \theta (\delta^{*} k)^{-\alpha} \rho^2.
  \end{equation}
  and if
   \begin{equation}\label{EPXa}
  \left| \{ x\in K_{2\rho}: n(x,t) > \mu_{+} - k \} \right| \leq (1-\beta)|K_{2\rho}|
  \end{equation}
  for all $t\in (-2\tau, 0]$, then we have
  \begin{equation}\label{EPXb}
  \left| \{ (x,t)\in K_{\rho}\times[-\tau, 0]: n > \mu_{+} - \delta^{*}k \} \right| \leq \nu|K_{\rho}\times [-\tau, 0]|.
  \end{equation}
\end{proposition}

\begin{proof}
Let $k_{j}=2^{-j}k$ for $j=0,1,2,\ldots,j^{*}$ with $j^{*}$ to be
determined later. Denote that $\delta^{*}=2^{-j^{*}}$. For
simplicity, denote that
  \[
  A_{j} = \{(x,t)\in K_{\rho}\times [-\tau, 0]: n(x,t) < k_{j}\}.
  \]
Let $\zeta$ be a poecewise linear cutoff function
  \[
  \zeta = \begin{cases}
            1, & \mbox{ in } K_{\rho}\times[-\tau, 0] \\
            0, & \mbox{ on the parabolic boundary of } K_{2\rho}\times [-2\tau, 0].
          \end{cases}
  \]
  satisfying, for all $j=1, \ldots, j^{*}$,
  \[
  |\nabla \zeta| \leq \frac{1}{\rho},
  \text{ and }
  \zeta_{t} \leq \frac{1}{\tau} = \frac{1}{\theta} (\delta^{*}k)^{\alpha}\rho^{-2} \leq \frac{1}{\theta} k_{j}^{\alpha}\rho^{-2}.
  \]
Then the local energy estimate Proposition~\ref{P:EE} provides the
following:
\begin{equation}\label{EPX01}\begin{split}
& (1+\alpha)\int_{-2\tau}^{0}\int_{K_{2\rho}} n^{\alpha} \left|\nabla (n-\mu_{+}+k_{j})_{+}\zeta]\right|^2 \,dx\,dt \\
&\leq 2(1+\alpha) \int_{-2\tau}^{0}\int_{K_{2\rho}} k_{j}^{\alpha} (n-\mu_{+}+k_{j})_{+}^{2} \zeta \partial_{t}\zeta \,dx\,dt \\
&\quad + 16(1+\alpha)\int_{-2\tau}^{0}\int_{K_{2\rho}} (n^{\alpha} + k_{j}^{\alpha})(n-\mu_{+}+k_{j})_{+}^{2} |\nabla \zeta|^2 \,dx\,dt \\
&\quad + 16 k^{2}(k^{-\alpha}\|B\|^{2}_{2\hat{q}_1, 2\hat{q}_2} +
\|\nabla B\|_{\hat{q}_1, \hat{q}_2}) |2\tau|^{\frac{2(1+\kappa)}{q_2}}
|K_{2\rho}|^{\frac{2(1+\kappa)}{q_1}}
\end{split}\end{equation}

In the set $\{(n-\mu_{+}+k_{j})_{+}>0 \},$ we observe that $n >
\mu_{+}-k_{j}$. Because of \eqref{k-rho}, the last term is
simplifies as
\[
16 \frac{k^{2-\alpha}}{\rho^2}(k^{-\alpha}\|B\|^{2}_{2\hat{q}_1, 2\hat{q}_2} + \|\nabla B\|_{\hat{q}_1, \hat{q}_2}) |2\tau|^{\frac{2(1+\kappa)}{q_2}} |K_{2\rho}|^{\frac{2(1+\kappa)}{q_1}}
\leq C\frac{k^2}{\rho^2}|K_{2\rho}\times[-2\tau, 0]|.
\]
Then the inequality \eqref{EPX01} yields
\begin{equation}\label{EPX02}
\int_{-2\tau}^{0}\int_{K_{2\rho}} \left|\nabla (n-\mu_{+}+k_{j})_{+}
\zeta]\right|^2 \,dx\,dt \leq C_0  \frac{k_{j}^{2}}{\rho^2}
\left|K_{2\rho} \times [-2\tau, 0]\right|
\end{equation}
where $C_0 = C_0 (\alpha, \theta, \gamma,\|B\|^{2}_{2\hat{q}_1,
2\hat{q}_2}, \|\nabla B\|_{\hat{q}_1, \hat{q}_2} )$.

By taking integration of the obtained inequality obtained from Lemma~\ref{Iso} for all $t \in [-\tau, 0]$, we have
\begin{equation}\label{EPX04}
 2^{-2}k_{j}^{2} |K_{\rho} \times [-\tau, 0] \cap \{ n > \mu_{-} - k_{j} \}|
\leq \frac{\gamma \rho}{\beta}\iint_{A_{j}\setminus A_{j+1}} |\nabla
(n-\mu_{+}+k_{j})_{+}| \,dx\,dt.
\end{equation}
The inequality \eqref{EPXa} yields the following inequalty for all $t \in (-2\tau, 0)$
\[
\left| K_{\rho} \cap \{\mu_{+} - k_{j} < n < \mu_{+} - k_{j+1}\}
\right| \geq \beta \left| K_{2\rho}\right|.
\]

For the simplicity, let $\Omega_{\tau} = K_{\rho} \times [-\tau, 0]$.
Let us divide \eqref{EPX04} by $|A_{j}\setminus A_{j+1}|$ and apply
Jensen's inequality to obtain the following inequality:
\begin{equation}\label{EPX05}\begin{split}
 k_{j}^{2+2\alpha} \left(\frac{|A_{j}|}{|A_{j}\setminus A_{j+1}|}\right)^{2}
&\leq \left(\frac{\gamma \rho}{\beta}\right)^2 \frac{1}{|A_{j}\setminus A_{j+1}|}\iint_{A_{j}\setminus A_{j+1}} |\nabla (n-\mu_{+}+k_{j})_{+}|^2 \,dx\,dt \\
&\leq C_1 (\alpha, \theta, \gamma, \beta)
k_{j}^{2+2\alpha}\frac{|\Omega_{\tau}|}{|A_j \setminus A_{j+1}|}.
\end{split}\end{equation}

Hence, it yields constant $C_2$ independent of $j$ satisfying
\begin{equation}\label{EPX06}
\left(\frac{|A_{j}|}{|A_{j}\setminus A_{j+1}|}\right)^{2} \leq C_2
\frac{|\Omega_{\tau}|}{|A_j \setminus A_{j+1}|}.
\end{equation}
Therefore, \eqref{EPX06} provides that
\begin{equation}\label{EPX07}
\left(\frac{|A_{j+1}|}{|\Omega_{\tau}|}\right)^{2} =
\left(\frac{|A_{j+1}|}{|A_{j}\setminus A_{j+1}|}\right)^{2}
   \left(\frac{|A_{j}\setminus A_{j+1}|}{|\Omega_{\tau}|}\right)^{2}
\leq C_2 \frac{|A_{j}\setminus A_{j+1}|}{|\Omega_{\tau}|}.
\end{equation}
Then by taking the sum over $j=0, 1, \ldots, j^{*}-1$, we obtain
\[
j^{*} \left(\frac{|A_{j^{*}}|}{|\Omega_{\tau}|}\right)^{2} \leq C_2
\frac{|A_0 \setminus A_{j^{*}}|}{|\Omega_{\tau}|} \leq C_2.
\]
By choosing $j^{*}$ larege enough such that $ j^{*} \geq C_2
\nu^{-2}, $ it leads that $ |A_{j^{*}}| \leq \nu |\Omega_{\tau}|. $
\end{proof}

\subsection{Proof of two alternatives}\label{SS:ProofAlternatives}

First, we provide the proof of Lemma~\ref{L:FA}.

\textbf{Proof of Lemma~\ref{L:FA}}: Because of the assumption
\eqref{FAa}, we apply Proposition~\ref{P:DG} with $\theta =1$, $k =
\omega/2$, and $\rho = R$ to infer that
\[
n \geq \mu_{-} + \frac{\omega}{4} \quad \text{ on } K_{R} \times
\{T_1\}
\]
where
\[
T_1 = T_0 + \Delta - \left(\frac{\omega}{2}\right)^{-\alpha} R^2.
\]

Therefore, it holds, for $\nu =1$,
\[
\left|\{K_{R}: n(x,T_1) < \mu_{-} + \frac{\omega}{4}\} \right| \leq
(1-\nu)|K_{R}|.
\]
By applying Proposition~\ref{P:EPT} with $\rho = R$, $\nu = 1$, $t_0
= T_1$, $\theta = 1$ and $k = \omega/4$, for any $\epsilon \in
(0,1)$, there exists $\delta = \delta (\data, \epsilon, \theta_1)$
satisfying
\begin{equation}\label{FA01}
\left|\{K_{R}: n(x,t) < \mu_{-} + \frac{\delta \omega}{4}\} \right|
\leq (1-(1-\epsilon)\nu)|K_{R}|
\end{equation}
for all
\[
t \geq T_1 + \left(\frac{\omega}{4}\right)^{-\alpha} R^2.
\]
Let us choose
\[
\theta_1 = \theta_0 2^{2-\alpha} - 3 \cdot 2^{-\alpha}
\]
which provides that \eqref{FA01} holds for all time $[-T_1, 0]$.

Let us choose $\epsilon = \nu^{*}/(2^{2-\alpha}\theta_0)$ from
Proposition~\ref{P:DG02} and set
\[
\theta = \frac{-T_1}{\left(\delta \omega / 4\right)^{-\alpha} R^2}.
\]
We observe that $\theta \leq 2^{2-\alpha}\theta_0$. Hence,
Proposition~\ref{P:DG02} provides that
\[
\essinf_{Q} n(x,t) \geq \mu_{-} + \frac{\delta \omega}{8}
\]
where
\[
Q = K_{R/2} \times [- \theta \left(\frac{\delta
\omega}{4}\right)^{-\alpha} \left(\frac{R}{2}\right)^2, 0] .
\]
We make conclusion by choosing $\theta = 4^{-\alpha}\omega^{\alpha}$
and $\theta_0 = 2^{-2-\alpha}\delta^{\alpha}+1$.
\smallskip


Now, we provide the proof of Lemma~\ref{L:SA}.

\textbf{Proof of Lemma~\ref{L:SA}}: From the assumption \eqref{SAa},
we apply Proposition~\ref{P:prop1} with $k=\omega/2$, $\nu_1 =
\nu_0$, and $\rho = 2R$, that there exists
\[
\tau_1 \in \left( T_0, T_0 + \frac{\nu_0}{2-\nu_0} \Delta \right)
\]
such that
\[
\left| K_{2R} \cap \{ n(x, \tau_1) > \mu_{+}
-\frac{\omega}{2}\}\right| \leq \left( 1-
\frac{\nu_0}{2}\right)|K_{2R}|.
\]

We apply Proposition~\ref{P:EPT} with $\rho = 2R$, $k=\omega/2$,
$\epsilon = 1/2$ and $\theta = \theta_0$. Then there exists $\delta$
depending on data such that
\begin{equation}\label{SA01}
\left| K_{2R} \cap \{  n(x, t) > \mu_{+} -
\frac{\delta\omega}{2}\}\right| \leq \left( 1-
\frac{\nu_0}{4}\right)|K_{2R}| \text{ for all } t \in [\tau_1, 0 ]
\end{equation}
because $\tau_1 + \theta_0 \Delta \geq 0 $.

We are ready to apply Proposition~\ref{P:EPX} with $k=\delta \omega
/2$, $\theta = \theta_0 -1$, $\beta = \nu_0 /4$ and $\rho = R$. Then
for any $\nu \in (0,1)$, there there exists $\delta^{*}=
\delta^{*}(\data, \nu)$ such that
\[
\left|K_{R}\times[T, 0] \cap \{n > \mu_{+}- \frac{\delta^{*}\delta
\omega}{2}\}\right| \leq \nu^{*} |K_{R} \times [T, 0]|
\]
where
\[
T = -(\theta_0 -1)\left(\frac{\delta^{*}\delta
\omega}{2}\right)^{-\alpha} R^2,
\]
by determining $\theta_0 = 1+ (\delta^{*}\delta)^{-\alpha}$. Let us
choose $\nu = \nu_0$, the constant from Proposition~\ref{P:DG02}
which yields
\[
\esssup_{Q} n(x,t) \leq \mu_{+} - \frac{\delta^{*}\delta \omega}{4}
\]
where
\[
Q = K_{R/2} \times [-\left( \frac{\omega}{2}\right)^{-\alpha}
\left(\frac{R}{2}\right)^2, 0].
\]


\section{Proofs of Theorem \ref{weak-1} and Theorem \ref{weak-2}}

We introduce the approximate system of \eqref{eq:Chemotaxis}, which
is given by
\begin{eqnarray}\label{eq:regu-chemo}
\left\{
\begin{array}{cl}
\begin{split}
& \partial_{t}n_{\epsilon} - \Delta
(n_{\epsilon}+\epsilon)^{1+\alpha} + u_{\epsilon}\cdot \nabla
n_{\epsilon} = - \nabla \cdot (\chi
(c_{\epsilon})(n_{\epsilon}+\epsilon)^{q}\nabla c_{\epsilon}),
\\
& \partial_{t}c_{\epsilon} - \Delta c_{\epsilon} + u_{\epsilon}\cdot
\nabla c_{\epsilon} = -\kappa (c_{\epsilon})n_{\epsilon},
\\
& \partial_{t}u_{\epsilon} - \Delta u_{\epsilon} + \nabla
p_{\epsilon} = -n_{\epsilon}\nabla \phi,
\\
& \nabla \cdot u_{\epsilon}=0,
\end{split}
\end{array}
\right.
\end{eqnarray}
in $\mathbb{R}^{3}\times (0,T)$ with smooth initial data
$(n_{0_{\epsilon}}, c_{0_{\epsilon}}, u_{0_{\epsilon}})$ defined by
\begin{equation*}
n_{0_{\epsilon}}=\psi_{\epsilon}*n_{0}, \quad
c_{0_{\epsilon}}=\psi_{\epsilon}*c_{0} \quad \text{and} \quad
u_{0_{\epsilon}}=\psi_{\epsilon}*u_{0},
\end{equation*}
where $\psi_{\epsilon}$ denotes the usual mollifier with $\epsilon
\in (0, 1)$ and $*$ denotes the space convolution. \\
It is known that, due to the standard theory of existence and
regularity as done in \cite{FLM} and \cite{TW_2}, there exists a
classical solution of the equation \eqref{eq:regu-chemo} locally in
time for each $\epsilon \in (0, 1)$. In the sequel, although we have
to obtain estimates from the approximate system
\eqref{eq:regu-chemo}, we compute a priori estimates only for
simplicity, since it turns out that all computations are independent
of $\epsilon$.
\\
\\
Now we present the proof of Theorem \ref{weak-1}.

\begin{thm1.5}

We consider first  the case that  as $\frac{9q-8}{6}<\alpha \leq
\frac{3q-2}{2}$ and the other case $\alpha \geq \frac{3q-2}{2}$ will
be treated later.

$\bullet$\,\, (Case; $\frac{9q-8}{6}<\alpha \leq
\frac{3q-2}{2}$)\,\, Multiplying equation $\eqref{eq:Chemotaxis}_1$
with $\log n$ and using the integration by parts, we have
\[
\frac{d}{dt}\int_{\mathbb{R}^3}n\log n  +
\frac{4}{1+\alpha}\norm{\nabla n^{\frac{1+\alpha}{2}}}_{2}^{2} =
\frac{1}{q}\int_{\mathbb{R}^3}\nabla n^{q}\chi(c)\nabla c
\]
\begin{equation}\label{CKK-sept12-100}
\leq
\frac{2\overline\chi}{2\alpha-q+2}\bke{\int_{\mathbb{R}^3}\abs{\nabla
n^{\frac{2\alpha -q+2}{2}}}n^{\frac{3q-2\alpha-2}{2}}\abs{\nabla
c}}:=J_1,
\end{equation}
where $\overline\chi:=\underset{\mathbb
R^3_T}{\max}\abs{\chi(c(\cdot))}$. Here we remark that the
restriction that $\alpha \leq \frac{3q-2}{2}$ is due to the
requirement that $3q-2\alpha-2 \geq 0$ in \eqref{CKK-sept12-100}.
Applying Young's inequality,
\begin{equation}\label{Kim-Kang-sept2-10}
J_1\leq \varepsilon_1\norm{\nabla n^{\frac{2\alpha
-q+2}{2}}}_{2}^{2}+C_{\varepsilon_1}\int_{\mathbb{R}^3}n^{3q-2\alpha-2}\abs{\nabla
c}^{2},
\end{equation}
where $\varepsilon_1$ is a sufficiently small constant, which will
be chosen later. In the sequel, we indicate $\varepsilon_i$,
$i=1,2,\cdots$ as a small constant, which will be decided later.
Combining the above estimates, we obtain
\begin{equation}\label{eq:first}
\begin{split}
& \frac{d}{dt}\int_{\mathbb{R}^3}n\log n  +
\frac{4}{1+\alpha}\norm{\nabla
n^{\frac{1+\alpha}{2}}}_{2}^{2} \\
& \leq \varepsilon_{1}\norm{\nabla n^{\frac{2\alpha
-q+2}{2}}}_{2}^{2}+C_{\varepsilon_{1}}\int_{\mathbb{R}^3}n^{3q-2\alpha-2}\abs{\nabla
c}^{2}.
\end{split}
\end{equation}
Next,testing $n^{\alpha -q+1}$ to $\eqref{eq:Chemotaxis}_1$ and
using H\"{o}lder and Young's inequalities,
\begin{equation*}
\frac{1}{\alpha -q+2}\frac{d}{dt}\int_{\mathbb{R}^3}n^{\alpha-q+2}
+\frac{4(1+\alpha)(\alpha-q+1)}{(2\alpha-q+2)^2}\norm{\nabla
n^{\frac{2\alpha-q+2}{2}}}_{2}^{2}
\end{equation*}
\begin{equation*}
\leq \frac{\alpha-q+1}{1+\alpha}\overline\chi
\int_{\mathbb{R}^3}\abs{\nabla
n^{\frac{2\alpha-q+2}{2}}}n^{\frac{q}{2}}\abs{\nabla c}
\end{equation*}
\begin{equation}\label{Kim-Kang-sept2-20}
\leq \varepsilon_{2}\norm{\nabla n^{\frac{2\alpha-q+2}{2}}}_{2}^{2}+
C_{\varepsilon_{2}} \int_{\mathbb{R}^3} n^{q}\abs{\nabla c}^2.
\end{equation}
We estimate second term in righthand side of
\eqref{Kim-Kang-sept2-20} via integration by parts.
\begin{equation*}
\int_{\mathbb{R}^3}n^{q}\abs{\nabla c}^{2} =
\int_{\mathbb{R}^3}n^{q}\nabla c\cdot \nabla c\leq
\int_{\mathbb{R}^3}\abs{\nabla n^{q}} c\abs{\nabla c} +
\int_{\mathbb{R}^3} n^{q} c\abs{\nabla^2 c}.
\end{equation*}
\begin{equation*}
\leq C\int_{\mathbb{R}^3} \abs{ \nabla n^{\frac{2\alpha-q+2}{2}} }
n^{\frac{3q-2\alpha-2}{2}} \abs{\nabla c} +
C\int_{\mathbb{R}^3}n^{q}\abs{ \nabla^2 c }
\end{equation*}
\begin{equation}\label{K-Nov25-10}
\leq \varepsilon_{3}\norm{\nabla
n^{\frac{2\alpha-q+2}{2}}}_{2}^{2}+C_{\varepsilon_{3}}\int_{\mathbb{R}^3}n^{3q-2\alpha-2}\abs{\nabla
c}^{2} + C\int_{\mathbb{R}^3}n^{q}\abs{ \nabla^2 c },
\end{equation}
where 
we used that $\norm{c(t)}_{L^{\infty}(\mathbb{R}^3}\le
\norm{c_0}_{L^{\infty}(\mathbb{R}^3)}$. Combining the estimates
above and taking $\epsilon_i$, $i=1,2,3$ sufficiently small, we
conclude that
\begin{equation}\label{eq:second}
\frac{d}{dt}\int_{\mathbb{R}^3}n^{\alpha-q+2} +\norm{\nabla
n^{\frac{2\alpha-q+2}{2}}}_{2}^{2}\leq
C\bke{\int_{\mathbb{R}^3}n^{3q-2\alpha-2}\abs{\nabla c}^{2} +
\int_{\mathbb{R}^3}n^{q}\abs{ \nabla^2 c }}.
\end{equation}
Multiplying equation $\eqref{eq:Chemotaxis}_2$ with $-\Delta c$ and
integrating by parts, we have
\begin{equation*}
\frac{1}{2}\frac{d}{dt}\int_{\mathbb{R}^3}\abs{\nabla c}^{2}
+\norm{\nabla^2 c}_2^2\leq\int_{\mathbb R^3}\bke{u\cdot\nabla
c}\nabla^2 c+\overline\kappa\int_{\mathbb R^3}n\abs{\nabla^2 c},
\end{equation*}
where $\overline\kappa:=\underset{\mathbb
R^3_T}{\max}\abs{\kappa(c(\cdot))}$. Due to $\nabla \cdot u=0$ and
uniform bound of $c$, we observe that
\begin{align*}
        \int_{\mathbb R^3}\bke{u\cdot\nabla c}\Delta c&=\sum_{1\leq i,j\leq
        3}\int_{\mathbb R^3}~u_ic_{x_i}c_{x_jx_j}=-\sum_{1\leq i,j\leq
        3}\int_{\mathbb R^3}u_{i,x_j}c_{x_i}c_{x_j}\\
        &=\sum_{1\leq i,j\leq 3}\int_{\mathbb
        R^3}u_{i,x_j}cc_{x_jx_i}\leq C\norm{\nabla u}_2\norm{\nabla^2
        c}_2.
    \end{align*}
Summing up estimates, we obtain
    \begin{equation}\label{eq:third}
        \frac{1}{2}\frac{d}{dt}\int_{\mathbb{R}^3}\abs{\nabla c}^{2} +\norm{\nabla^2c}_2^2\leq
        C\bke{\norm{\nabla u}_2\norm{\nabla^2c}_2+\int_{\mathbb R^3}n\abs{\nabla^2c}}.
    \end{equation}
Let $M$ be a sufficiently large positive constant, which will be
specified later. Multiplying equation $\eqref{eq:Chemotaxis}_3$ with
$Mu$ and using the integration by parts, we note that
    \begin{equation}\label{eq:fourth}
        \frac{M}{2}\frac{d}{dt}\int_{\mathbb{R}^3}\abs{u}^{2} +2M\norm{\nabla
        u}_2^2\leq M\norm{\nabla \phi}_{L^{\infty}(\R^3)}\int_{\mathbb
        R^3}n\abs{u}.
    \end{equation}

We recall that $\int n \abs{\log n}$, we recall that (see e.g., (18)
of \cite{CKK})
\begin{equation}\label{nlogn-negative}
\int_{\mathbb R^3} n|\log n| \leq \int_{\mathbb R^3} n\log n
+2\int_{\mathbb R^3}\langle x\rangle n +C
\end{equation}
and we note that (see e.g., (19) of \cite{CKK})
\begin{equation}\label{nlogn-10}
\frac{d}{dt}\int_{\mathbb R^3} \langle x\rangle n  \leq
C\left(1+\norm{\nabla c}_2^2 + \norm{ \nabla
u}_{2}^{2}\right)+\varepsilon\norm{\nabla
n^{\frac{1+\alpha}{2}}}_{2}^{2}+\varepsilon\norm{\nabla
n^\frac{2\alpha-q+2}{2}}_2^2.
\end{equation}

Summing up \eqref{eq:first}$-$\eqref{eq:fourth} and
\eqref{nlogn-10}, and taking sufficiently small $\varepsilon_{i}$,
$i=1,2,3$, we have
\begin{equation}\label{eq:sum-1}
\begin{split}
& \frac{d}{dt}\left(\int_{\mathbb R^3} n\log n +\int_{\mathbb R^3}
\langle x\rangle n  +\frac{1}{\alpha -q+2}
\int_{\mathbb{R}^3}n^{\alpha -q+2}
+\frac{1}{2}\int_{\mathbb{R}^3}\mid \nabla
c \mid^{2} + \frac{M}{2}\int_{\mathbb{R}^3}\abs{u}^{2}  \right) \\
& \qquad + C \left( \norm{\nabla n^{\frac{1+\alpha}{2}}}_{2}^{2}+
\norm{\nabla^2 c}_{2}^{2}+M\norm{\nabla u}_{2}^{2} + \norm{\nabla
n^{\frac{2\alpha-q+2}{2}}}_{2}^{2}
\right) \\
& \leq  C \biggl(\int_{\mathbb{R}^3}n^{3q-2\alpha-2}\mid\nabla
c\mid^{2} + \int_{\mathbb{R}^3}n^{q}\mid \nabla^2 c \mid
\\
& \qquad +  \norm{\nabla u}_2\norm{\nabla^2 c}_2 +\int_{\mathbb
R^3}n\abs{\Delta c} + M\int_{\mathbb R^3}n\abs{u} \biggl) +C\bke{ 1
+ \norm{\nabla c}_2^2}
 \\
& = C\left(\text{{\Romannumeral 1}} + \text{{\Romannumeral 2}} +
\text{{\Romannumeral 3}} + \text{{\Romannumeral 4}} +
\text{{\Romannumeral 5}}\right) +C\bke{ 1 + \norm{\nabla c}_2^2}.
\end{split}
\end{equation}
Due to $0\le 3q-2\alpha-2<\frac{2}{3}$ via $\frac{9q-8}{6}<\alpha\le
\frac{3q-2}{2}$, we estimate $I$ as follows:
\begin{equation}\label{A1}
\text{{\Romannumeral 1}}\leq\left\{
\begin{array}{cl}
\norm{\nabla c}_2^2, \quad & \text{if} \quad \alpha=\frac{3q-2}{2},
\\
C_{\varepsilon_4}\norm{\nabla
c}_2^2+\varepsilon_{4}\norm{n_0}_1^\frac{2}{3}\norm{\nabla^2 c}_2^2,
\quad & \text{if} \quad  \frac{9q-8}{6}<\alpha<\frac{3q-2}{2},
\end{array}
\right.
\end{equation}
where H\"{o}lder inequality and Sobolev embedding are used.

To estimate the term $\text{\Romannumeral 2}$, applying H\"{o}lder,
Young's and interpolation inequalities, we have
\[
\text{{\Romannumeral 2}} = \int_{\mathbb{R}^3}n^{q}\abs{ \nabla^2 c
}
 \leq C_{\varepsilon_{5}}\norm{n}_{2q}^{2q}+
\varepsilon_{5}\norm{\nabla^2 c}_{2}^{2} \leq
C_{\varepsilon_{5}}\norm{n}_{1}^{2q\theta}\norm{n}_{3(2\alpha-q+2)}^{2q(1-\theta)}
+ \varepsilon_{5}\norm{\nabla^2 c}_{2}^{2}
\]
\[
\leq C_{\varepsilon_{5}}\norm{n}_{1}^{2q\theta}\norm{\nabla
n^{\frac{2\alpha-q+2}{2}}}_{2}^{2q(1-\theta)\frac{2}{2\alpha-q+2}}+\varepsilon_{5}\norm{\nabla^2
c}_{2}^{2},
\]
where $\theta=\frac{6\alpha-5q+6}{2q(6\alpha-3q+5)}$. Since $\alpha
> \frac{9q-8}{6}$, we observe that
\[
2q(1-\theta)\frac{2}{2\alpha-q+2}=\frac{6(2q-1)}{6\alpha-3q+5}<2.
\]
Therefore, we have
\begin{equation}\label{A2}
\text{{\Romannumeral 2}} \leq
C_{\varepsilon_{5}}C_{\varepsilon_{6}}+\varepsilon_{6}\norm{\nabla
n^{\frac{2\alpha-q+2}{2}}}_{2}^{2}+\varepsilon_{5}\norm{\nabla^2
c}_{2}^{2}.
\end{equation}

The term $\text{{\Romannumeral 3}}$ is easily estimated as follows:
\begin{equation}\label{A3}
\text{{\Romannumeral 3}}=\norm{\nabla u}_2\norm{\nabla^2 c}_2 \leq
C_{\varepsilon_{7}}\norm{\nabla
u}_{2}^{2}+\varepsilon_{7}\norm{\nabla^2 c}_{2}^{2}.
\end{equation}

Next, we estimate the term $\text{\Romannumeral 4}$. H\"{o}lder and
Young's and interpolation inequalities yield
\[
\text{{\Romannumeral 4}}=\int_{\mathbb R^3}n\abs{\nabla^2 c}  \leq
C_{\varepsilon_{8}}\norm{n}_{1}^{2\theta}\norm{n}_{3(2\alpha-q+2)}^{2(1-\theta)}
+ \varepsilon_{8}\norm{\nabla^2 c}_{2}^{2}
\]
\[
\leq C_{\varepsilon_{8}}\norm{n}_{1}^{2\theta}\norm{\nabla
n^{\frac{2\alpha-q+2}{2}}}_{2}^{2(1-\theta)\frac{2}{2\alpha-q+2}}+\varepsilon_{8}\norm{\nabla^2
c}_{2}^{2},
\]
where $\theta=\frac{6\alpha-3q+4}{2(6\alpha-3q+5)}$. Similarly as
above, we note, due to $\alpha
> \frac{9q-8}{6}$, that
\[
2(1-\theta)\frac{2}{2\alpha-q+2}=\frac{6}{6\alpha-3q+5}<2.
\]
Therefore, we have
\begin{equation}\label{A4}
\text{{\Romannumeral 4}} \leq
C_{\varepsilon_{8}}C_{\varepsilon_{9}}+\varepsilon_{9}\norm{\nabla
n^{\frac{2\alpha-q+2}{2}}}_{2}^{2}+\varepsilon_{8}\norm{\nabla^2
c}_{2}^{2}.
\end{equation}
Finally, the term $\text{\Romannumeral 5}$ is estimated via
H\"{o}lder, Young's inequalities and Sobolev embedding. Indeed,
\begin{equation*}
\text{{\Romannumeral 5}}=M\int_{\mathbb R^3}n\abs{u}~dx  \leq M
\bke{\norm{n}_{\frac{6}{5}}\norm{u}_{6}} \leq
M \bke{C_{\varepsilon_{10}}\norm{n}_{\frac{6}{5}}^{2}+\varepsilon_{10}\norm{\nabla u}_{2}^{2}}. \\
\end{equation*}
\begin{equation*}
\begin{split}
\leq M \bke{
C_{\varepsilon_{10}}\norm{n}_{1}^{2\theta}\norm{n}_{3(2\alpha-q+2)}^{2(1-\theta)}
+\varepsilon_{10}\norm{\nabla u}_{2}^{2}}
\end{split}
\end{equation*}
\begin{equation*}
\begin{split}
\leq M \bke{ C_{\varepsilon_{10}}\norm{n}_{1}^{2\theta}\norm{\nabla
n^{\frac{2\alpha-q+2}{2}}}_{2}^{2(1-\theta)\frac{2}{2\alpha-q+2}}+\varepsilon_{10}\norm{\nabla
u}_{2}^{2}},
\end{split}
\end{equation*}
where $\theta=\frac{10\alpha-5q+8}{2(6\alpha-3q+1)}$. As before, via
$\alpha> \frac{9q-8}{6}$, we can see that
\[
2(1-\theta)\frac{2}{2\alpha-q+2}=\frac{2}{6\alpha-3q+5}<2.
\]
Therefore, we obtain
\begin{equation}\label{A5}
\text{{\Romannumeral 5}} \leq M \bke{
C_{\varepsilon_{10}}C_{\varepsilon_{11}}+\varepsilon_{11}\norm{\nabla
n^{\frac{2\alpha-q+2}{2}}}_{2}^{2}+\varepsilon_{10}\norm{\nabla
u}_{2}^{2}}.
\end{equation}
Summing up \eqref{A1}$-$\eqref{A5}, we have for sufficiently small
$\varepsilon_{i}$, $i=4, \cdots , 11$
\begin{equation}\label{eq:sum-2}
\begin{split}
& \sup_{0\leq t \leq T}\left(\int_{\mathbb R^3}  n\log
n+\int_{\mathbb R^3} \langle x\rangle n +\frac{1}{\alpha -q+2}
\int_{\mathbb{R}^3}n^{\alpha -q+2} +\int_{\mathbb{R}^3}\mid \nabla
c \mid^{2} + \frac{1}{2}\int_{\mathbb{R}^3}\abs{u}^{2}  \right) \\
& + C \int_{0}^{T}\left( \norm{\nabla
n^{\frac{1+\alpha}{2}}}_{2}^{2} +  \norm{\nabla^2
c}_{2}^{2}+M\norm{\nabla u}_{2}^{2} + \norm{\nabla
n^{\frac{2\alpha-q+2}{2}}}_{2}^{2} \right) \\
& \leq C\bke{1+\int_{0}^{T} \left( \norm{\nabla c}_2^2 +1
\right)}\leq C,
\end{split}
\end{equation}
where $C=C\bke{T, \norm{c_0}_{L^{\infty}\cap H^{1}},
\int_{\mathbb{R}^3}n_{0}\log n_{0}, \norm{n_0}_{\alpha -q +2},
\norm{n_0(1+\langle x\rangle)}_{1}, \norm{u_0}_2}$. Combining
estimates \eqref{nlogn-negative} and \eqref{eq:sum-2}
\begin{equation}\label{eq:sum-2-1}
\begin{split}
& \sup_{0\leq t \leq T}\left(\int_{\mathbb R^3}  n\abs{\log n} +
\langle x\rangle n + \int_{\mathbb{R}^3}n^{\alpha -q+2}
+\int_{\mathbb{R}^3}\mid \nabla
c \mid^{2} + \int_{\mathbb{R}^3}\abs{u}^{2}  \right) \\
& + \int_{0}^{T}\left( \norm{\nabla n^{\frac{1+\alpha}{2}}}_{2}^{2}
+  \norm{\nabla^2 c}_{2}^{2}+M\norm{\nabla u}_{2}^{2} + \norm{\nabla
n^{\frac{2\alpha-q+2}{2}}}_{2}^{2} \right) \\
& \leq C\left(1+ \int_{0}^{T} \left( \norm{\nabla c}_2^2 +1
\right)\right) \leq C,
\end{split}
\end{equation}
where $C=C \left( T, \norm{c_0}_{L^{\infty}\cap H^{1}},
\norm{n_{0}(1+|x|+|\log n_{0}|)}_1, \norm{n_0}_{\alpha -q +2},
\norm{u_0}_2 \right)$.

 $\bullet$\,\,$\bke{\text{Case} \ \alpha > \frac{3q-2}{2}}$
Multiplying equation $\eqref{eq:Chemotaxis}_1$ with $\log n$ and
integrating it by parts, we have
\begin{equation*}\label{eq:first-1}
\begin{split}
\frac{d}{dt}\int_{\mathbb{R}^3}n\log n  +
\frac{4}{1+\alpha}\norm{\nabla n^{\frac{1+\alpha}{2}}}_{2}^{2} =
\frac{1}{q}\int_{\mathbb{R}^3}\chi(c)\nabla n^{q}\cdot\nabla c,
\end{split}
\end{equation*}
\begin{equation}\label{eq:first-1-2}
\begin{split}
& =-\frac{1}{q}\int_{\mathbb{R}^3} \bke{\chi{'}(c) n^{q} \abs{
\nabla
c }^{2} + \chi(c)n^{q}\nabla^2 c } \\
& \leq \frac{1}{q}\int_{\mathbb{R}^3} \bke{\chi{'}(c) n^{q} \abs{
\nabla c }^{2} + \chi(c)n^{q}\abs{ \nabla^2 c }}.
\end{split}
\end{equation}
Using estimates \eqref{Kim-Kang-sept2-20}, \eqref{eq:third},
\eqref{eq:fourth}, \eqref{nlogn-10} together with \eqref{eq:first-1}
\begin{equation}\label{eq:sum-3}
\begin{split}
& \frac{d}{dt}\left(\int_{\mathbb R^3} n\left(\log n + 2\langle
x\rangle \right) +\frac{1}{\alpha -q+2} \int_{\mathbb{R}^3}n^{\alpha
-q+2} +\frac{1}{2}\int_{\mathbb{R}^3}\mid \nabla
c \mid^{2} + \frac{1}{2}\int_{\mathbb{R}^3}\abs{u}^{2}  \right) \\
& \hspace{0.8cm} + C \left( \norm{\nabla
n^{\frac{1+\alpha}{2}}}_{2}^{2}+ \norm{\nabla^2
c}_{2}^{2}+M\norm{\nabla u}_{2}^{2} + \norm{\nabla
n^{\frac{2\alpha-q+2}{2}}}_{2}^{2}
\right) \\
& \leq C \biggl(\int_{\mathbb{R}^3}n^{q}\abs{\nabla c}^{2} +
n^{q}\abs{ \nabla^2 c } + \norm{\nabla u}_2\norm{\nabla^2 c}_2
+\int_{\mathbb R^3}n\abs{\nabla^2 c} + \int_{\mathbb
R^3}n\abs{u} \biggr) +C\bke{ 1 + \norm{\nabla c}_2^2}\\
& = C \left( \text{{\Romannumeral 1}} + \text{{\Romannumeral 2}} +
\text{{\Romannumeral 3}} + \text{{\Romannumeral 4}} +
\text{{\Romannumeral 5}} \right)+C\bke{ 1 + \norm{\nabla c}_2^2}.
\end{split}
\end{equation}
We estimate $\text{{\Romannumeral 2}}$, $\text{{\Romannumeral 3}}$,
$\text{{\Romannumeral 4}}$ and $\text{{\Romannumeral 5}}$ exactly
the same ways as \eqref{A2}, \eqref{A3}, \eqref{A4} and \eqref{A5},
respectively. It remains to estimate $\text{{\Romannumeral 1}}$.
Due to $0<q<\frac{2\alpha -q+2}{2}$ via $\alpha > \frac{3q-2}{2}$,
we have
\begin{equation*}
\text{{\Romannumeral 1}} \leq \int_{\mathbb{R}^3}\left(
C_{\varepsilon_{12}}+\varepsilon_{12} n^{
\frac{2\alpha-q+2}{2}}\right)\abs{\nabla c}^{2} \leq
C_{\varepsilon_{12}}\norm{\nabla c}_{2}^{2} +
\varepsilon_{12}\int_{\mathbb{R}^3}n^{
\frac{2\alpha-q+2}{2}}\abs{\nabla c}^{2}.
\end{equation*}
We note that
\begin{equation}\label{first term}
\begin{split}
& \int_{\mathbb{R}^3}n^{ \frac{2\alpha-q+2}{2}}\abs{\nabla c}^{2}
\leq \int_{\mathbb{R}^3}\abs{\nabla n^{
\frac{2\alpha-q+2}{2}}}\abs{\nabla c}  + \int_{\mathbb{R}^3}n^{
\frac{2\alpha-q+2}{2}}\abs{ \nabla^2 c } \\
& \leq \frac{1}{2}\bke{\norm{\nabla
n^{\frac{2\alpha-q+2}{2}}}_{2}^{2}+\norm{\nabla^2
c}_{2}^{2}+\norm{\nabla c}_{2}^{2}+ \norm{
n^{\frac{2\alpha-q+2}{2}}}_{2}^{2}}.
\end{split}
\end{equation}
The last term of the right hand side in \eqref{first term} is
estimated as follows:
\begin{equation*}
\begin{split}
\norm{ n^{\frac{2\alpha-q+2}{2}}}_{2}^{2} =
\norm{n}_{2\alpha-q+2}^{2\alpha-q+2} \leq
\norm{n}_{1}^{(2\alpha-q+2)\theta}\norm{n}_{3(2\alpha-q+2)}^{(2\alpha-q+2)(1-\theta)}
\end{split}
\end{equation*}
\begin{equation*}
\begin{split}
\le C\norm{n}_{1}^{(2\alpha-q+2)\theta}\norm{\nabla
n^{\frac{2\alpha-q+2}{2}}}_{2}^{(2\alpha-q+2)(1-\theta)\frac{2}{2\alpha-q+2}},
\end{split}
\end{equation*}
where $\theta=\frac{2}{6\alpha-3q+5}$. We then note that
\[
(2\alpha-q+2)(1-\theta)\frac{2}{2\alpha-q+2}=\frac{6(2\alpha-q+1)}{6\alpha-3q+5}<2.
\]
Hence we obtain
\begin{equation}\label{A1-1}
\text{{\Romannumeral 1}} \leq C\norm{\nabla c}_{2}^{2} +
C\varepsilon_{12}\bke{\norm{\nabla^2 c}_{2}^{2} + \norm{ \nabla
n^{\frac{2\alpha-q+2}{2}}}_{2}^{2}}.
\end{equation}
Adding up estimates, we conclude that
\begin{equation}\label{eq:sum-4}
\begin{split}
& \sup_{0\leq t \leq T}\left(\int_{\mathbb R^3}  n\log
n+\int_{\mathbb R^3} \langle x\rangle n +\frac{1}{\alpha -q+2}
\int_{\mathbb{R}^3}n^{\alpha -q+2} +\int_{\mathbb{R}^3}\mid \nabla
c \mid^{2} + \frac{1}{2}\int_{\mathbb{R}^3}\abs{u}^{2}  \right) \\
& + C \int_{0}^{T}\left( \norm{\nabla
n^{\frac{1+\alpha}{2}}}_{2}^{2} +  \norm{\nabla^2
c}_{2}^{2}+M\norm{\nabla u}_{2}^{2} + \norm{\nabla
n^{\frac{2\alpha-q+2}{2}}}_{2}^{2} \right) \\
& \leq C\bke{1+\int_{0}^{T} \left( \norm{\nabla c}_2^2 +1
\right)}\leq C,
\end{split}
\end{equation}
where $C=C\bke{T, \norm{c_0}_{L^{\infty}\cap H^{1}},
\int_{\mathbb{R}^3}n_{0}\log n_{0}, \norm{n_0}_{\alpha -q +2},
\norm{n_0(1+\langle x\rangle)}_{1}, \norm{u_0}_2}$. Combining
estimates \eqref{nlogn-negative} and \eqref{eq:sum-4}
\begin{equation}
\begin{split}
& \sup_{0\leq t \leq T}\left(\int_{\mathbb R^3}  n\abs{\log n} +
\langle x\rangle n + \int_{\mathbb{R}^3}n^{\alpha -q+2}
+\int_{\mathbb{R}^3}\mid \nabla
c \mid^{2} + \int_{\mathbb{R}^3}\abs{u}^{2}  \right) \\
& + \int_{0}^{T}\left( \norm{\nabla n^{\frac{1+\alpha}{2}}}_{2}^{2}
+  \norm{\nabla^2 c}_{2}^{2}+M\norm{\nabla u}_{2}^{2} + \norm{\nabla
n^{\frac{2\alpha-q+2}{2}}}_{2}^{2} \right) \\
& \leq C\left(1+ \int_{0}^{T} \left( \norm{\nabla c}_2^2 +1
\right)\right) \leq C,
\end{split}
\end{equation}
where $C=C \left( T, \norm{c_0}_{L^{\infty}\cap H^{1}},
\norm{n_{0}(1+|x|+|\log n_{0}|)}_1, \norm{n_0}_{\alpha -q +2},
\norm{u_0}_2 \right)$. This completes the proof.
\end{thm1.5}

We present proof of Theorem \ref{weak-2}.
\begin{thm1.6}
Since we showed already the case $\alpha
> \frac{9q-8}{6}$ in Theorem \ref{weak-1}, it suffices to prove
the case $2q-2 < \alpha \leq \frac{9q-8}{6}$, $1 \leq q \leq
\frac{4}{3}$.

Multiplying equation $\eqref{eq:Chemotaxis}_1$ with $\log n$ and
using the integration by parts, we have
\begin{equation*}
\frac{d}{dt}\int_{\mathbb{R}^3}n\log n  +
\frac{4}{1+\alpha}\norm{\nabla n^{\frac{1+\alpha}{2}}}_{2}^{2} =
\frac{1}{q}\int_{\mathbb{R}^3}\nabla n^{q}\chi(c)\nabla c
\end{equation*}
\begin{equation*}
\begin{split}
\leq \frac{2\overline{\chi}}{1+\alpha}\int_{\mathbb{R}^3}\abs{\nabla
c}\abs{\nabla n^{\frac{1+\alpha}{2}}}n^{\frac{2q-\alpha-1}{2}} .
\end{split}
\end{equation*}
where $\overline{\chi}$ denote $\underset{\mathbb
R^3_T}{\max}\abs{\chi(c(\cdot))}$. Applying H\"{o}lder and Young's
inequalities, we have
\\
\begin{equation*}
\frac{2\overline{\chi}}{1+\alpha}\int_{\mathbb{R}^3}\abs{\nabla
c}\abs{\nabla n^{\frac{1+\alpha}{2}}}n^{\frac{2q-\alpha-1}{2}}  \leq
\int_{\mathbb{R}^3}\abs{\nabla c}\abs{\nabla
n^{\frac{1+\alpha}{2}}}\bke{C_{\varepsilon_{1}}+\varepsilon_{1}
n^{\frac{1}{2}}}
\end{equation*}
\begin{equation*}
\leq C_{\varepsilon_{1}}\bke{C_{\varepsilon_{2}}\norm{\nabla
c}_{2}^{2} + \varepsilon _{2}\norm{\nabla
n^{\frac{1+\alpha}{2}}}_{2}^{2}} + \varepsilon_{1}
\bke{C_{\varepsilon_{3}}\norm{\nabla n^{\frac{1+\alpha}{2}}}_{2}^{2}
+ \varepsilon_{3}\norm{n^{\frac{1}{2}}\nabla c}_{2}^{2}}.
\end{equation*}
Hence we have
\begin{equation}\label{eq:first-2-1}
\frac{d}{dt}\int_{\mathbb{R}^3}n\log n  + \bke{
\frac{4}{1+\alpha}-C_{\varepsilon_{1}}\varepsilon_{2}-C_{\varepsilon_{3}}\varepsilon_{1}
} \norm{\nabla n^{\frac{1+\alpha}{2}}}_{2}^{2} \leq
C_{\varepsilon_{1}}C_{\varepsilon_{2}}\norm{\nabla c}_{2}^{2} +
\varepsilon_{1}\varepsilon_{3}\norm{n^{\frac{1}{2}}\nabla
c}_{2}^{2},
\end{equation}
where $2q-\alpha-1 > 0$ and $\frac{2q-\alpha-1}{2} < \frac{1}{2}$,
which is equivalent to $2q-2 < \alpha \leq
\frac{9q-8}{6}$, $1 \leq q \leq \frac{4}{3}$.\\
Multiplying equation $\eqref{eq:Chemotaxis}_1$ with $n^{\alpha
-2q+2}$ and using the integration by parts, we have
\[
\frac{1}{\alpha-2q+3}\frac{d}{dt}\int_{\mathbb{R}^3}\abs{n}^{\alpha-2q+3}+\frac{4(1+\alpha)(\alpha-2q+2)}{(2\alpha-2q+3)^2}\norm{\nabla
n^{\frac{2\alpha-2q+3}{2}}}_{2}^{2}
\]
\[
\leq
\frac{\alpha-2q+2}{\alpha-q+2}\overline\chi\int_{\mathbb{R}^3}\abs{\nabla
n^{\alpha - q+2}}\cdot\abs{\nabla c}.
\]
Applying H\"{o}lder and Young's inequalities, we have
\begin{equation}\label{eq:second-2}
\frac{1}{\alpha-2q+3}\frac{d}{dt}\int_{\mathbb{R}^3}\abs{n}^{\alpha-2q+3}+\bke{\frac{4(1+\alpha)(\alpha-2q+2)}{(2\alpha-2q+3)^2}-\varepsilon_{7}}\norm{\nabla
n^{\frac{2\alpha-2q+3}{2}}}_{2}^{2} \leq
C_{\varepsilon_{7}}\int_{\mathbb{R}^3}n\abs{\nabla c}^{2} ,
\end{equation}
where $\alpha-2q+3 > 1$, which is equivalent to $\alpha > 2q-2$.
Multiplying equation $\eqref{eq:Chemotaxis}_2$ with $-\nabla^2 c$
and using the integration by parts, we have
\begin{equation*}
\begin{split}
\frac{1}{2}\frac{d}{dt}\int_{\mathbb{R}^3}|\nabla
c|^{2}  +\norm{\nabla^2 c}^{2}_{2} & \leq \int_{\mathbb{R}^3}(u \cdot \nabla c)\nabla^2 c  + \int_{\mathbb{R}^3}\kappa (c)n\abs{\nabla^2 c}  \\
& = \text{{\Romannumeral 1}} + \text{{\Romannumeral 2}}.
\end{split}
\end{equation*}
The term $\text{\Romannumeral 1}$ is estimated as follows. Applying
H\"{o}lder and Young's inequalities, we have
\begin{equation*}
\text{{\Romannumeral 1}} = C_{\varepsilon_{8}}\norm{\nabla
u}_{2}^{2}+\varepsilon_{8}\norm{\nabla^2 c}_{2}^{2}.
\end{equation*}
Using the integration by parts, we have
\begin{equation*}
\begin{split}
\text{{\Romannumeral 2}} = \int_{\mathbb{R}^3}\kappa
(c)n\abs{\nabla^2 c}  & =
-\int_{\mathbb{R}^3}\kappa^{\prime}(c)n\abs{\nabla c}^{2}
 +
\int_{\mathbb{R}^3}\kappa (c) \abs{\nabla n}\cdot \abs{\nabla c}. \\
\end{split}
\end{equation*}
And using H\"{o}lder and Young's inequalities, we have
\[
\text{{\Romannumeral 2}} \leq
-\kappa_{0}\int_{\mathbb{R}^3}n\abs{\nabla c}^{2}  +
\int_{\mathbb{R}^3}\kappa (c)\abs{\nabla
n^{\frac{1+\alpha}{2}}}(C_{\varepsilon_{9}}+\varepsilon_{9}
n^{\frac{1}{2}})\abs{\nabla c}
\]
\[
\leq -\kappa_{0}\int_{\mathbb{R}^3}n\abs{\nabla c}^{2}  +
\bke{\varepsilon_{10}+\varepsilon_{9}\varepsilon_{11}}\norm{\nabla
n^{\frac{1+\alpha}{2}}}_{2}^{2} +
C_{\varepsilon_{8}}C_{\varepsilon_{9}}\norm{\nabla c}_{2}^{2} +
\varepsilon_{9}C_{\varepsilon_{11}}\norm{n^{\frac{1}{2}}\nabla
c}_{2}^{2}.
\]
Hence we have
\begin{equation}\label{eq:third-2}
\begin{split}
\frac{1}{2}\frac{d}{dt} & \int_{\mathbb{R}^3}\abs{\nabla c}^{2}
+(1-\varepsilon_{8})\norm{\nabla^2
c}^{2}+ \bke{\kappa_{0}-\varepsilon_{9}C_{\varepsilon_{11}}}\int_{\mathbb{R}^3}n\abs{\nabla c}^{2} \\
& \leq
\bke{\varepsilon_{10}+\varepsilon_{9}\varepsilon_{11}}\norm{\nabla
n^{\frac{1+\alpha}{2}}}_{2}^{2} + C_{\varepsilon_{8}}\norm{ \nabla c
}_{2}^{2} + C_{\varepsilon_{8}}\norm{\nabla u}_{2}^{2},
\end{split}
\end{equation}
where $\frac{1}{1-\alpha}>1$, which is equivalent to $0<\alpha <1$.
Let $M$ be a sufficiently large positive constant, which will be
decided later. Multiplying equation $\eqref{eq:Chemotaxis}_3$ with
$Mu$ and using the integration by parts, we have
\begin{equation*}
\frac{M}{2}\frac{d}{dt}\int_{\mathbb{R}^3}\abs{u}^{2}
+2M\norm{\nabla
        u}_2^2\leq M\norm{\nabla \phi}_{L^{\infty}(\R^3)}\int_{\mathbb R^3}n\abs{u}.
\end{equation*}
Using interpolation and Young's inequalities, we have
\begin{equation*}
\frac{M}{2}\frac{d}{dt}\int_{\mathbb{R}^3}\abs{u}^{2}
+2M\norm{\nabla u}_{2}^{2} \leq
C_{\varepsilon_{12}}\norm{n}_{\frac{6}{5}}^{2}+\varepsilon_{12}\norm{\nabla
u}_{2}^{2}
\end{equation*}
\begin{equation*}
\leq
C_{\varepsilon_{12}}\norm{n}_{1}^{2\theta}\norm{n}_{3(1+\alpha)}^{2(1-\theta)}+\varepsilon_{12}\norm{\nabla
u}_{2}^{2} \leq
C_{\varepsilon_{12}}C_{\varepsilon_{13}}+\varepsilon_{13}\norm{\nabla
n^{\frac{1+\alpha}{2}}}_{2}^{2}+\varepsilon_{13}\norm{\nabla
u}_{2}^{2},
\end{equation*}
where $\alpha > 0$. Hence we have
\begin{equation}\label{eq:fourth-2}
\frac{M}{2}\frac{d}{dt}\int_{\mathbb{R}^3}\abs{u}^{2}
+(2M-\varepsilon_{12})\norm{\nabla u}_{2}^{2} \leq
C_{\varepsilon_{12}}C_{\varepsilon_{13}}+\varepsilon_{13}\norm{\nabla
n^{\frac{1+\alpha}{2}}}_{2}^{2}.
\end{equation}
\begin{equation}\label{nlogn-100}
\frac{d}{dt}\int_{\mathbb R^3} \langle x\rangle n  \leq
C\left(1+\norm{\nabla c}_2^2 + \norm{\nabla
u}_{2}^{2}\right)+\varepsilon\norm{\nabla n^\frac{1+\alpha}{2}}_2^2
+ \varepsilon\norm{\nabla n^\frac{2\alpha-2q+3}{2}}_2^2.
\end{equation}
Summing up \eqref{eq:first-2-1}$-$\eqref{nlogn-100}, we have we have
for sufficiently small $\varepsilon_{i}$, $i=1, \cdots , 13$
\begin{equation}\label{eq:sum_5}
\begin{split}
& \sup_{0 \leq t \leq T}\left(\int_{\mathbb R^3} n\left(\log n +
2\langle x\rangle \right)  + \int_{\mathbb{R}^3}n^{\alpha -2q+3}
+\int_{\mathbb{R}^3}\abs{ \nabla
c }^{2} + \int_{\mathbb{R}^3}\abs{u}^{2} \right) \\
& \hspace{0.3cm} + C \int_{0}^{T}\left( \norm{\nabla
n^{\frac{1+\alpha}{2}}}_{2}^{2}+ \norm{\nabla
n^{\frac{2\alpha-2q+3}{2}}}_{2}^{2}+\norm{\nabla^2
c}_{2}^{2}+\norm{\nabla
u}_{2}^{2} + \int_{\mathbb{R}^3}n\abs{\nabla c}^{2}\right) \\
& \leq C\bke{1+\int_{0}^{T} \left( \norm{\nabla c}_2^2 +1
\right)}\leq C,
\end{split}
\end{equation}
where $C=C\bke{T, \norm{c_0}_{L^{\infty}\cap H^{1}},
\int_{\mathbb{R}^3}n_{0}\log n_{0}, \norm{n_0}_{\alpha -2q +3},
\norm{n_0(1+\langle x\rangle)}_{1}, \norm{u_0}_2}$. Combining
estimates \eqref{nlogn-negative} and \eqref{eq:sum_5}
\begin{equation}
\begin{split}
& \sup_{0\leq t \leq T}\left(\int_{\mathbb R^3}  n\abs{\log n} +
\langle x\rangle n + \int_{\mathbb{R}^3}n^{\alpha -q+2}
+\int_{\mathbb{R}^3}\mid \nabla
c \mid^{2} + \int_{\mathbb{R}^3}\abs{u}^{2}  \right) \\
& + \int_{0}^{T}\left( \norm{\nabla n^{\frac{1+\alpha}{2}}}_{2}^{2}
+  \norm{\nabla^2 c}_{2}^{2}+M\norm{\nabla u}_{2}^{2} + \norm{\nabla
n^{\frac{2\alpha-q+2}{2}}}_{2}^{2} \right) \\
& \leq C\left(1+ \int_{0}^{T} \left( \norm{\nabla c}_2^2 +1
\right)\right) \leq C,
\end{split}
\end{equation}
where $C=C \left( T, \norm{c_0}_{L^{\infty}\cap H^{1}},
\norm{n_{0}(1+|x|+|\log n_{0}|)}_1, \norm{n_0}_{\alpha -2q +3},
\norm{u_0}_2 \right)$. This completes the proof.
\end{thm1.6}


\section{Proofs of Theorem \ref{bdd weak-1} and Theorem \ref{bdd weak-2}}

\begin{lemma}\label{velocity-sept12-10}
Let $u$ be a solution of $\eqref{eq:Chemotaxis}_3$ constructed in
Theorem \ref{weak-1}. If $\alpha > \text{max} \left\{2q-2, \
\frac{9q-8}{6}
 \right\}$, then $u\in
L^{\infty}_tL^6_x$. Similarly, Suppose that $u$ is a solution of
$\eqref{eq:Chemotaxis}_3$ constructed in Theorem \ref{weak-2}. If
$\alpha> \text{max}\left\{ \text{min}\left\{ 2q-2,
\frac{9q-8}{6}\right\}, \ \frac{10q-9}{8} \right\}$, then $u\in
L^{\infty}_tL^6_x$.
\end{lemma}
\begin{proof}

$\bullet$\,(Case : $\alpha > \text{max} \left\{2q-2, \
\frac{9q-8}{6}
 \right\}$).

Since $u$ is a solution of $\eqref{eq:Chemotaxis}_3$ constructed in
Theorem \ref{weak-1}, we remind that $\int_{0}^{T}\norm{\nabla
n^{\frac{2\alpha -q+2}{2}}}_{2}^{2}<\infty$. We consider the
vorticity equation of $\eqref{eq:Chemotaxis}_3$
\begin{equation}\label{E:omega}
\omega_{t}-\Delta \omega = - \nabla \times (n\nabla \phi),
\end{equation}
where $\omega=\nabla \times u$. The energy estimate yileds
\begin{equation*}
\frac{d}{dt}\int_{\mathbb{R}^3}\abs{\omega}^{2}  + \norm{\nabla
\omega}_{2}^{2} \leq C\int_{\mathbb{R}^3}\abs{\nabla \omega} \cdot n
 \leq \frac{1}{2}\norm{\nabla \omega}_{2}^{2}+C\norm{n}_{2}^{2}.
\end{equation*}
Therefore, integrating time, we obtain
\begin{equation}\label{vol-1}
\int_{\mathbb{R}^3}\abs{\omega}^{2}  + \int_{0}^{T}\norm{\nabla
\omega}_{2}^{2} \leq C\int_{0}^{T}\norm{n}_{2}^{2}.
\end{equation}
Applying H\"{o}lder, interpolation inequalities  and Sobolev
embedding, we note
\begin{equation}\label{vol-2}
\begin{split}
\int_{0}^{T}\norm{n}_{2}^{2} & \leq
\int_{0}^{T}\norm{n}_{1}^{2\theta}\norm{n}_{3(2\alpha-q+2)}^{2(1-\theta)}
\leq \int_{0}^{T}\norm{n}_{1}^{2\theta}\norm{\nabla n^{\frac{2\alpha
-q+2}{2}}}_{2}^{2(1-\theta)\frac{2}{2\alpha -q+2}} \\
& \leq C\int_{0}^{T}\norm{\nabla n^{\frac{2\alpha
-q+2}{2}}}_{2}^{2(1-\theta)\frac{2}{2\alpha -q+2}},
\end{split}
\end{equation}
where $ \theta=(6\alpha -3q+4)/(12\alpha -6q+10)). $ Since $\alpha
> \frac{9q-8}{6}$, we note that $4(1-\theta)/(2\alpha
-q+2)<2$, which implies that the righthand side of \eqref{vol-2} is
finite. Therefore, $\omega\in L_{x,t}^{2,\infty}$, which immediately
yields $u \in L_{x,t}^{6,\infty}$.
\\
\\
$\bullet$\,(Case : $\alpha > \text{max}\left\{ \text{min}\left\{
2q-2, \frac{9q-8}{6}\right\}, \ \frac{10q-9}{8} \right\})$. Since
$\text{max}\left\{ \text{min}\left\{ 2q-2, \frac{9q-8}{6}\right\}, \
\frac{10q-9}{8} \right\})\ge \frac{10q-9}{8}$, it is enough to
consider the case $\frac{10q-9}{8}< \alpha$. We first treat the case
that $\frac{10q-9}{8}< \alpha < 2q-1$. We note, due to H\"{o}lder,
interpolation inequalities  and Sobolev embedding, that
\begin{equation}\label{vol-3}
\int_{0}^{T}  \norm{n}_{2}^{2} \leq \int_{0}^{T}\norm{n}_{\alpha
-2q+3}^{2\theta}\norm{n}_{3(2\alpha-2q+3)}^{2(1-\theta)}  \leq
C\int_{0}^{T}\norm{\nabla n^{\frac{2\alpha
-2q+3}{2}}}_{2}^{2(1-\theta)\frac{2}{2\alpha -2q+3}},
\end{equation}
where
\[
\theta=\frac{(\alpha -2q+3)(6\alpha-6q -7)}{2(5\alpha -4q+6)}
\]
and we used that $\norm{n}_{L^{\infty}_tL^{\alpha -2q+3}_x}<C$ in
Theorem \ref{weak-2}. Observing that $4(1-\theta)/2\alpha -2q+3<2$,
we can see that the righthand side of \eqref{vol-3} is finite. For
the case that $\alpha \geq 2q-1$, using a different interpolation
inequality, we estimate
\begin{equation}\label{vol-4}
\int_{0}^{T}\norm{n}_{2}^{2} \leq
\int_{0}^{T}\norm{n}_{1}^{2\theta}\norm{n}_{3(2\alpha-2q+3)}^{2(1-\theta)}
\leq C\int_{0}^{T}\norm{\nabla n^{\frac{2\alpha
-2q+3}{2}}}_{2}^{2(1-\theta)\frac{2}{2\alpha -2q+3}},
\end{equation}
where
\[
\theta=\frac{6\alpha -6q+7}{2(6\alpha -6q+8)}.
\]
Since $\alpha \geq 2q-1$, we note that $4(1-\theta)/2\alpha
-2q+3<2$, which implies that the righthand side of \eqref{vol-4} is
finite. Due to estimates \eqref{vol-2}, \eqref{vol-3} and
\eqref{vol-4}, we deduce the lemma.
\end{proof}

We present proof of Theorem \ref{bdd weak-1}.

\begin{thm1.8}
$\bullet$ $\bke{\text{Case} \ \text{max} \left\{2q-2, \
\frac{9q-8}{6}
 \right\} <\alpha \ \text{and} \ \frac{3q-1}{6}\le\alpha}$
\\
Multiplying equation $\eqref{eq:Chemotaxis}_1$ with $n^{p-1}$ and
using the integration by parts, we have
\\
\[
\frac{1}{p}\frac{d}{dt}\int_{\mathbb{R}^3}\abs{n}^{p}+\frac{4(p-1)(1+\alpha)}{(p+\alpha)^2}\norm{\nabla
n^{\frac{p+\alpha}{2}}}_{2}^{2} \leq
\frac{p-1}{p+q-1}\int_{\mathbb{R}^3}\chi(c)\abs{\nabla c}\abs{
\nabla n^{p+q-1}}
\]
\[
\leq \frac{p-1}{p+q-1}\int_{\mathbb{R}^3}\chi(c)\abs{\nabla c}
\abs{\nabla n^{\frac{p+\alpha}{2}}}n^{\frac{p+2q-\alpha -2}{2}}.
\]
Applying H\"{o}lder and Young's inequalities, we have
\begin{equation}\label{oct01-10}
\frac{1}{p}\frac{d}{dt}\int_{\mathbb{R}^3}\abs{n}^{p}+\frac{4(p-1)(1+\alpha)}{(p+\alpha)^2}\norm{\nabla
n^{\frac{p+\alpha}{2}}}_{2}^{2} \leq \varepsilon_{1}\norm{\nabla
n^{\frac{p+\alpha}{2}}}_{2}^{2}+C_{\varepsilon_{1}}\int_{\mathbb{R}^3}n^{p+2q-\alpha
-2}\abs{\nabla c}^2.
\end{equation}
The right hand side of the above is estimated as
\\
\[
\int_{\mathbb{R}^3}n^{p+2q-\alpha -2}\abs{\nabla c}^2\leq
\norm{n^{p+2q-\alpha -2}}_{\frac{p}{p+2q-\alpha -2}}\norm{\mid
\nabla c \mid^{2}}_{\frac{p}{\alpha +2-2q}}
\]
\[
\leq \norm{n}_{p}^{p+2q-\alpha -2}\norm{\nabla^2
c}_{\frac{6p}{2p+3\alpha+6-6q}}^{2} \leq C \bke{1+
\norm{n}_{p}^{p}}\norm{\nabla^2 c}_{\frac{6p}{2p+3\alpha+6-6q}}^{2}.
\]
Therefore, taking
$\varepsilon_1\le\frac{2(p-1)(1+\alpha)}{(p+\alpha)^2}$, we obtain
\[
\frac{d}{dt}\norm{n}_{p}^{p} \leq
Cp^2\bke{1+\norm{n}_{p}^{p}}\norm{\nabla^2
c}_{\frac{6p}{2p+3\alpha+6-6q}}^{2}.
\]
Gronwall inequality implies that
\begin{equation}\label{sept27-20}
\sup_{0 \leq t \leq T}\norm{n}_{p}^{p} \leq \exp \left\{
Cp^2\int_{0}^{T}\norm{\nabla^2 c}_{\frac{6p}{2p+3\alpha+6-6q}}^{2}
\right\}\int_{0}^{T}\norm{\nabla^2
c}_{\frac{6p}{2p+3\alpha+6-6q}}^{2} + \norm{n_0}_{p}^{p}.
\end{equation}
We will show that $\int_{0}^{T}\norm{\nabla^2
c}_{\frac{6p}{2p+3\alpha+6-6q}}^{2}<\infty$ for any $p>\alpha -q+2$.
Indeed, from maximal regularity theorem for heat equation, we have
\begin{equation}\label{sept27-30}
\begin{split}
\int_{0}^{T}\norm{\nabla^2 c}_{\frac{6p}{2p+3\alpha+6-6q}}^{2} &
\leq
C\int_{0}^{T}\bke{\norm{n}_{\frac{6p}{2p+3\alpha+6-6q}}^{2}+\norm{u\cdot
\nabla
c}_{\frac{6p}{2p+3\alpha+6-6q}}^{2}}+C\norm{\nabla^2 c_0}_{\frac{6p}{2p+3\alpha+6-6q}}^{2} \\
& = C(\text{{\Romannumeral 1}}+\text{{\Romannumeral
2}})+C\norm{\nabla^2 c_0}_{\frac{6p}{2p+3\alpha+6-6q}}^{2} .
\end{split}
\end{equation}
The term $\text{{\Romannumeral 1}}$ is estimated as follows. Since $
\frac{3q-1}{6}\le\alpha$, we have via interpolation inequality and
Sobolev embedding
\begin{equation}\label{sept26-10}
\text{{\Romannumeral 1}}  =\int_{0}^{T}
\norm{n}_{\frac{6p}{2p+3\alpha+6-6q}}^{2} \leq
\int_{0}^{T}\norm{n}_{1}^{2\theta}\norm{n}_{3(2\alpha-q+2)}^{2(1-\theta)}
\leq \int_{0}^{T}\norm{\nabla n^{\frac{2\alpha
-q+2}{2}}}_{2}^{2-\delta_{p}},
\end{equation}
where
\[
2(1-\theta)=\frac{(2\alpha -q+2)(4p-3\alpha -6+6q)}{p(6\alpha-3q+5)}
\ \ \text{and}
\]
\[
2(1-\theta)\frac{2}{2\alpha-q+2}=\frac{8}{6\alpha-3q+5}-\frac{6\alpha
+12-12q}{p(6\alpha-3q+5)}=2-\delta_{p}<2.
\]
Thus, it is direct that the term $\text{{\Romannumeral 1}}$ is
finite. On the other hand, applying H\"{o}lder inequality and Lemma
\ref{velocity-sept12-10}, we estimate the term $\text{{\Romannumeral
2}}$ as follows.
\[
\text{{\Romannumeral 2}} = \int_{0}^{T}\norm{u\cdot \nabla
c}_{\frac{6p}{2p+3\alpha+6-6q}}^{2} \leq \int_{0}^{T}\norm{u}_{6}^{2}\norm{\nabla c}_{\frac{6p}{p+3\alpha +6-6q}}^{2} \\
\leq C\int_{0}^{T}\norm{\nabla^2 c}_{\frac{2p}{p+\alpha +2-2q}}^{2}.
\]
Using the maximal regularity for heat equation, interpolation
inequality, Sobolve embedding and $\frac{2p}{p+\alpha
+2-2q}<\frac{6p}{2p+3\alpha+6-6q}$, we have
\[
\int_{0}^{T}\norm{\nabla^2 c}_{\frac{2p}{p+\alpha +2-2q}}^{2} \leq
C\int_{0}^{T} \bke{ \norm{n}_{\frac{2p}{p+\alpha
+2-2q}}^{2}+\norm{u\cdot \nabla c}_{\frac{2p}{p+\alpha
+2-2q}}^{2}}+\norm{\nabla^2 c_0}_{\frac{2p}{p+\alpha +2-2q}}^{2}.
\]
\[
\leq
C\int_{0}^{T}\bke{\norm{n}_{1}^{2}+\norm{n}_{\frac{6p}{2p+3\alpha+6-6q}}^{2}
+\norm{u}_{6}^{2}\norm{\nabla
c}_{\frac{6p}{2p+3\alpha+6-6q}}^{2}}+\norm{\nabla^2
c_0}_{\frac{2p}{p+\alpha +2-2q}}^{2}.
\]
\begin{equation}\label{sept26-20}
\leq C\int_{0}^{T}\bke{\norm{n}_{\frac{6p}{2p+3\alpha+6-6q}}^{2} +
\norm{\nabla c}_{\frac{6p}{2p+3\alpha+6-6q}}^{2}}+CT+\norm{\nabla^2
c_0}_{\frac{2p}{p+\alpha +2-2q}}^{2}.
\end{equation}
We note that the first term in \eqref{sept26-20} is the same as
$\text{{\Romannumeral 1}}$ in \eqref{sept26-10} and thus it is
finite. It is straightforward that $\int_{0}^{T}\norm{\nabla
c}_{\frac{6p}{2p+3\alpha+6-6q}}^{2}<\infty$, since
$2<\frac{6p}{2p+3\alpha+6-6q}<3$. Hence, the second term
$\text{{\Romannumeral 2}}$ is also finite, which deduces the
boundedness of $L^p$-norm of $n$, for any $p>\alpha-q+2$.
\\
\\
$\bullet$ $\bke{\text{Case} \ \text{max}\left\{\frac{9q-8}{6}, \
2q-2 \right\} <\alpha < \frac{3q-1}{6}}$\quad We first note that
this case is equivalent to the case that
$\frac{9q-8}{6}<\alpha<\frac{3q-1}{6}$ with $1\le q<\frac{7}{6}$. We
set $p$ with
\begin{equation}\label{condition p}
\text{max}\left\{ \alpha - q + 2, \ 3\alpha-4q+4 \right\} < p <
4\alpha -5q+6.
\end{equation}
Testing $n^{p-1}$  to equation $\eqref{eq:Chemotaxis}_1$ and
following similar computations as in previous case, we have
\begin{equation}\label{L-p rightterm}
\frac{1}{p}\frac{d}{dt}\int_{\mathbb{R}^3}\abs{n}^{p}+\frac{4(p-1)(1+\alpha)}{(p+\alpha)^2}\norm{\nabla
n^{\frac{p+\alpha}{2}}}_{2}^{2} \leq \varepsilon_{2} \norm{\nabla
n^{\frac{p+\alpha}{2}}}_{2}^{2}+C_{\varepsilon_{2}}\int_{\mathbb{R}^3}n^{p+2q-\alpha
-2}\abs{\nabla c}^2.
\end{equation}
Noting that
\begin{equation*}
\int_{\mathbb{R}^3}n^{p+2q-\alpha -2}\abs{\nabla
c}^2=\int_{\mathbb{R}^3}n^{p+2q-\alpha -2}\nabla c\cdot\nabla c
\end{equation*}
\begin{equation*}
\leq \int_{\mathbb{R}^3}\abs{\nabla n^{p+2q-\alpha -2}}c\abs{ \nabla
c}  +\int_{\mathbb{R}^3}\abs{ n^{p+2q-\alpha -2}}c\abs{ \nabla^2 c}
\end{equation*}
and integrating in time, we estimate \eqref{L-p rightterm} as
\[
\sup_{0 \leq t \leq
T}\norm{n}_{p}^{p}+\bke{\frac{4p(p-1)(1+\alpha)}{(p+\alpha)^2}-\varepsilon_{2}}\int_{0}^{T}\norm{\nabla
n^{\frac{p+\alpha}{2}}}_{2}^{2}
\]
\begin{equation}\label{sept29-20}
\leq Cp^2\int_{0}^{T}\int_{\mathbb{R}^3}\abs{\nabla n^{p+2q-\alpha
-2}}\abs{ \nabla c} +Cp^2\int_{0}^{T}\int_{\mathbb{R}^3}\abs{
n^{p+2q-\alpha -2}}\abs{ \nabla^2 c} :=
Cp^2\bke{\text{{\Romannumeral 1}} + \text{{\Romannumeral 2}}}.
\end{equation}
We first estimate $\text{{\Romannumeral 1}}$. Applying H\"{o}lder
and Young's inequalities, we observe that
\begin{equation}\label{sept29-30}
\text{{\Romannumeral 1}}\leq \int_{0}^{T}\int_{\mathbb{R}^3}\abs{
\nabla n^{\frac{p+\alpha}{2}} } n^{\frac{p+4q-3\alpha -4}{2}}\abs{
\nabla c} \leq \int_{0}^{T}\bke{\varepsilon_{3}\norm{\nabla
n^{\frac{p+\alpha}{2}}}_{2}^{2}+C_{\varepsilon_{3}}\int_{\mathbb{R}^3}
n^{p+4q-3\alpha -4}\abs{ \nabla c}^{2} }.
\end{equation}
Let $r_{1}=\frac{6(\alpha-q+2)}{14\alpha-17q+22-3p}$. Due to
\eqref{condition p}, H\"{o}lder inequality and Sobolev embedding, we
have
\[
\int_{0}^{T}\int_{\mathbb{R}^3} n^{p+4q-3\alpha -4}\abs{ \nabla
c}^{2}
\]
\[
\leq \int_{0}^{T}\norm{n^{p+4q-3\alpha
-4}}_{\frac{\alpha-q+2}{^{p+4q-3\alpha -4}}}\norm{\abs{ \nabla c
}^2}_{\frac{\alpha-q+2}{^{4\alpha -5q+6-p}}}
\]
\[
\leq \int_{0}^{T}\norm{n}_{\alpha-q+2}^{p+4q-3\alpha
-4}\norm{\nabla^2 c}_{r_1}^{2} \le C\int_{0}^{T}\norm{\nabla^2
c}_{r_1}^{2},
\]
where we used that $\norm{n}_{L^{\infty}_tL^{\alpha-q+2}_x}<C$
proved in Theorem \ref{weak-1}. From maximal regularity for heat
equation and results in Lemma \ref{velocity-sept12-10}, we note
\[
\int_{0}^{T}\norm{\nabla^2 c}_{r_1}^{2} \leq
C\int_{0}^{T}\bke{\norm{n}_{r_1}^{2}+\norm{u\cdot \nabla
c}_{r_1}^2}+ C\norm{\nabla^2 c_0}_{r_1}^{2}
\]
\[
\leq C\int_{0}^{T}\bke{\norm{n}_{r_1}^{2}+
\norm{u}_{6}^{2}\norm{\nabla c}_{\frac{6r_1}{6-r_1}}^{2}} +
C\norm{\nabla^2 c_0}_{r_1}^{2}
\leq C\int_{0}^{T}\bke{\norm{n}_{r_1}^{2}+\norm{\nabla^2
c}_{\frac{6r_1}{6+r_1}}^{2}} + C\norm{\nabla^2 c_0}_{r_1}^{2}
\]
\[
\leq
C\int_{0}^{T}\bke{\norm{n}_{r_1}^{2}+\norm{n}_{\frac{6r_1}{6+r_1}}^{2}+\norm{u}_{6}^{2}\norm{\nabla
c}_{r_1}^{2}}+ C\bke{\norm{\nabla^2 c_0}_{r_1}^{2}+\norm{\nabla^2
c_0}_{\frac{6r_1}{6+r_1}}^{2}}
\]
\begin{equation}\label{sept23-20}
\leq
C\int_{0}^{T}\bke{\norm{n}_{r_1}^{2}+\norm{n}_{1}^{2}+\norm{u}_{6}^{2}\norm{\nabla
c}_{r_1}^{2}}+ C:=\text{{\Romannumeral 1}}_{1}
\end{equation}
If $r_1\le \alpha - q + 2$, then \eqref{sept23-20} is bounded, due
to the result of Theorem \ref{weak-1}, by
\begin{equation}\label{sept23-30}
\text{{\Romannumeral 1}}_{1}\le C\int_{0}^{T}\norm{\nabla
c}_{r_1}^{2}+ C(1+T).
\end{equation}
On the other hand, in case that $r_1>\alpha - q + 2$, we can see
that $r_{1} < 3(p+\alpha)$, since $\text{max}\left\{ \alpha - q + 2,
\ 3\alpha-4q+4 \right\} < p < 4\alpha -5q+6$, and thus
\eqref{sept23-20} is estimated as
\[
\text{{\Romannumeral 1}}_1 \leq C\int_{0}^{T}\bke{\norm{n}_{\alpha
-q+2}^{2\theta_{1}}\norm{n}_{3(p+\alpha)}^{2(1-\theta_{1})}
+\norm{\nabla c}_{r_1}^{2}}+C(1+T)
\]
\begin{equation}\label{sept23-10}
\leq C\int_{0}^{T}\bke{\norm{\nabla
n^{\frac{p+\alpha}{2}}}_{2}^{\delta_{1}}+ \norm{\nabla
c}_{r_1}^{2}}+C(1+T),
\end{equation}
where
\[
2(1-\theta_{1})=\frac{6(p+\alpha)(r_{1}-(\alpha-q+2))}{r_{1}(3p+2\alpha+q-2)},\quad
\delta_1=\frac{12(r_1-\alpha+q-2)}{r_1(3p+2\alpha+q-2)}.
\]
Here we used that $\alpha - q + 2 < r_{1} < 3(p+\alpha)$ and $r_1\ge
\frac{6r_1}{r_1+6}$. We note that $\delta_1<2$, since
$\text{max}\left\{ 2q-2, \ \frac{9q-8}{6} \right\} <\alpha <
\frac{3q-1}{6}$. Next we estimate the term $\text{{\Romannumeral
2}}$. Let $r_{2}:=p+2q-\alpha-1$. Using H\"older, Young's and
maximal regularity for heat equation, and following similar
computations as in \eqref{sept23-20}, $\text{{\Romannumeral 2}}$ is
estimated as follows:
\[
\text{{\Romannumeral 2}}
\leq C\int_{0}^{T}\bke{\norm{n}_{r_{2}}^{r_{2}}+\norm{\nabla^2
c}_{r_{2}}^{r_{2}}}\leq C\int_{0}^{T}\bke{\norm{n}_{r_{2}}^{r_{2}} +
\norm{u\cdot \nabla c}_{r_{2}}^{r_{2}}}+C\norm{\nabla^2
c_0}_{r_2}^{r_2}
\]
\[
\leq
C\int_{0}^{T}\left(\|n\|_{r_2}^{r_2}+\|n\|_{\frac{6r_2}{6+r_2}}^{r_2}+\|\nabla
c\|_{r_2}^{r_2}\right) +C\bke{\norm{\nabla^2
c_0}_{r_2}^{r_2}+\norm{\nabla^2 c_0}_{\frac{6r_1}{6+r_1}}^{2}}.
\]
We note that
\[
\int_{0}^{T}\left(\|n\|_{r_2}^{r_2}+\|n\|_{\frac{6r_2}{6+r_2}}^{r_2}+\|\nabla
c\|_{r_2}^{r_2}\right)\le
\int_{0}^{T}\left(C\|n\|_{r_2}^{r_2}+C\|n\|_{1}^{r_2}+\|\nabla
c\|_{r_2}^{r_2}\right)
\]
\[
\leq C\int_{0}^{T}\left(\|n\|_{\alpha
-q+2}^{r_2\theta_{3}}\|n\|_{3(p+\alpha)}^{r_2(1-\theta_{3})}
+\|\nabla c\|_{2}^{r_{2}\theta_{4}}\|\nabla
c\|_{6}^{r_{2}(1-\theta_{4})}\right)+CT
\]
\[
\leq C\int_{0}^{T}\left(\|n\|_{\alpha
-q+2}^{r_2\theta_{3}}\|n\|_{3(p+\alpha)}^{r_2(1-\theta_{3})}
+\|\nabla c\|_{2}^{r_{2}\theta_{4}}\|\nabla^2
c\|_{2}^{r_{2}(1-\theta_{4})}\right)+CT
\]
\[
\leq C\int_{0}^{T}\left(\|\nabla
n^{\frac{p+\alpha}{2}}\|_{2}^{\delta_{3}}+ \|\nabla
c\|_{2}^{r_{2}\theta_{4}}\|\nabla^2 c\|_{2}^{\delta_{4}}\right),
\]
where
\[
\delta_{3}=\frac{2r_{2}(1-\theta_{3})}{p+\alpha}=\frac{6(p-2\alpha+3q-3)}{3p+2\alpha+q-2},\qquad
\delta_4=r_{2}(1-\theta_{4})=\frac{3}{2}(p+2q-\alpha-3).
\]
Here we used that $\alpha -q+2 < r_{2}< 3(p+ \alpha)$ and
$2<r_{2}<6$ and we observe that $\delta_3<2$ and $\delta_4<2$, since
$\frac{9q-8}{6} <\alpha < \frac{3q-1}{6}$ and $\text{max}\left\{
\alpha - q + 2, \ 3\alpha-4q+4 \right\} < p < 4\alpha -5q+6$.

Combining estimates of $\text{{\Romannumeral 1}}$ and
$\text{{\Romannumeral 2}}$, we obtain
\begin{equation*}
\sup_{0 \leq t \leq T}\norm{n(t)}_{p}^{p}+C\int_{0}^{T}\norm{\nabla
n^{\frac{p+\alpha}{2}}}_{2}^{2}  \leq C
\bke{\int_{0}^{T}\norm{\nabla c}_{r_{1}}^{2} +1}.
\end{equation*}
We can see also that  $2<r_{1}<6$ as long as  $\frac{11\alpha
-14q+16}{3} < p < 4\alpha -5q+6$, which is valid, since
$\frac{11\alpha -14q+16}{3}> 3\alpha-4q+4$, in case that
$\frac{9q-8}{6} <\alpha < \frac{3q-1}{6}$. Therefore, for any $p$
with $\text{max}\left\{ \alpha - q + 2, \ 3\alpha-4q+4 \right\} < p
< 4\alpha -5q+6$ we obtain
\begin{equation}\label{sept27-10}
n\in L^{\infty}(0, T;L^{p}(\mathbb{R}^3)),\qquad  \nabla
n^{\frac{p+\alpha}{2}} \in L^{2}(0, T;L^{2}(\mathbb{R}^3))
\end{equation}

Let $p_{0}=\frac{6-3\alpha}{4}$. We then see that $1 < p_{0} <
4\alpha -5q+6$ via $\frac{9q-8}{6} <\alpha < \frac{3q-1}{6}$, and
thus it is evident from \eqref{sept27-10} that
\begin{equation}\label{Lp0}
n \in L^{\infty}(0, T;L^{p_{0}}(\mathbb{R}^{3})).
\end{equation}

Next, we will show that $n \in L^{\infty}(0,
T;L^{p}(\mathbb{R}^{3}))$ for any $p_0<p<\infty$. Similarly as
before, multiplying equation $\eqref{eq:Chemotaxis}_1$ with
$n^{p-1}$ and using Gronwall inequality, we get \eqref{sept27-20},
namely,
\[
\sup_{0 \leq t \leq T}\norm{n}_{p}^{p} \leq \exp \left\{
Cp^2\int_{0}^{T}\norm{\nabla^2 c}_{\frac{6p}{2p+3\alpha+6-6q}}^{2}
\right\}\int_{0}^{T}\norm{\nabla^2
c}_{\frac{6p}{2p+3\alpha+6-6q}}^{2} + \norm{n_0}_{p}^{p}.
\]
We recall \eqref{sept27-30} via  maximal regularity for heat
equation, i.e.
\[
\int_{0}^{T}\norm{\nabla^2 c}_{\frac{6p}{2p+3\alpha+6-6q}}^{2} \leq
C\int_{0}^{T}\bke{\norm{n}_{\frac{6p}{2p+3\alpha+6-6q}}^{2}+\norm{u\cdot
\nabla c}_{\frac{6p}{2p+3\alpha+6-6q}}^{2}} + C\norm{\nabla^2
c_0}_{\frac{6p}{2p+3\alpha+6-6q}}^{2}
\]
\[
= C(\text{{\Romannumeral 3}}+\text{{\Romannumeral 4}}) +
C\norm{\nabla^2 c_0}_{\frac{6p}{2p+3\alpha+6-6q}}^{2}.
\]
The term $\text{{\Romannumeral 3}}$ is estimated as follows. For
$p>\alpha -q+2$, we have
$p_{0}<\frac{6p}{2p+3\alpha+6-6q}<3p_{0}+3\alpha$ and hence it
follows from \eqref{Lp0} that
\begin{equation}\label{sept29-10}
\text{{\Romannumeral 3}} = \norm{n}_{\frac{6p}{2p+3\alpha+6-6q}}^{2}
\leq
\norm{n}_{p_{0}}^{2\theta}\norm{n}_{3(p_{0}+\alpha)}^{2(1-\theta)}
\leq \norm{n}_{p_{0}}^{2\theta}\norm{\nabla
n^{\frac{p_{0}+\alpha}{2}}}_{2}^{2-\delta_{p}},
\end{equation}
where
\[
2(1-\theta)\frac{p_{0}+\alpha}{2}=\frac{12-4p_{0}}{2p_{0}+3\alpha}-\frac{2p_{0}(3\alpha
+6-6q)}{p(2p_{0}+3\alpha)}=2-\delta_{p}.
\]
The term $\text{{\Romannumeral 4}}$ is estimated exactly in the same
way as $\text{{\Romannumeral 2}}$ in ......, and thus we obtain
\begin{equation}\label{sept27-40}
\text{{\Romannumeral 4}} \leq
C\int_{0}^{T}\bke{\norm{n}_{\frac{6p}{2p+3\alpha+6-6q}}^{2} +
\norm{\nabla c}_{\frac{6p}{2p+3\alpha+6-6q}}^{2}}+CT+\norm{\nabla^2
c_0}_{\frac{2p}{p+\alpha +2-2q}}^{2}.
\end{equation}

The first term in \eqref{sept27-40} can be treated as the case
$\text{{\Romannumeral 3}}$ and the second term is bounded, since
$2<\frac{6p}{2p+3\alpha+6-6q}<3$ and $\int_0^T\norm{\nabla c}_m^2
ds<\infty$ for $2\leq m \leq 3$. We finally conclude the boundedness
of $L^\infty$-norm of $n$. Indeed, since $n \in
L^{\infty}(0,T;L^{p}(\mathbb{R}^3))$ for all $1 \leq p < \infty$, we
can see that $c_{t}, \ \nabla^c, \ u_{t}$ and $\nabla^{2}u$ belong
to $L^{p}((0,T)\times \mathbb{R}^3)$ for all $p<\infty$ and
therefore, we also note, due to parabolic embedding, that $\nabla c
\in L^{\infty}((0,T)\times \mathbb{R}^3)$. Using estimate
\eqref{oct01-10} and $\nabla c \in L^{\infty}((0,T)\times
\mathbb{R}^3)$, we obtain
\begin{equation}\label{L-infiny-1}
\frac{d}{dt}\int_{\mathbb{R}^3}\abs{n}^{p} \leq Cp^2
\int_{\mathbb{R}^3}n^{p(1-\delta)}, \qquad
\delta:=\frac{\alpha+2-2q}{p} .
\end{equation}
Multiplying equation $\eqref{L-infiny-1}$ with $p$, we have
\begin{equation*}
\frac{d}{dt}\int_{\mathbb{R}^3}\abs{n}^{p}~dx \leq
Cp^2\norm{n}_{p(1-\delta)}^{p(1-\delta)}.
\end{equation*}
Using interpolation inequality, we have
\begin{equation*}
\frac{d}{dt}\int_{\mathbb{R}^3}\abs{n}^{p}~dx \leq
Cp^2\norm{n}_{p}^{p(1-\delta)(1-\theta)},
\end{equation*}
where $p(1-\delta)(1-\theta)=\frac{p(p(1-\delta)-1)}{p-1}$. Let
$y(t):=\norm{n(t)}_{p}^{p}$ and $\beta := 1-p(1-\delta)(1-\theta)$.
Then we have
\[
y(t)^\prime \le Cp^{2}y(t)^{{1-\beta}}.
\]
Via Gronwall inequality, we observe that
\begin{equation}\label{oct03-10}
\norm{n(t)}_{p} \leq \bke{Cp^2 \beta t }^{\frac{1}{\beta p}} +
\norm{n(0)}_{p}, \qquad t \leq T.
\end{equation}
Passing $p$ to the limit, we obtain for all $t \leq T$.
\begin{equation*}
\norm{n(t)}_{\infty} \leq 1+\norm{n(0)}_{\infty}.
\end{equation*}
H\"older continuity is a direct consequence of Theorem
\ref{T:Holder1}. Indeed, since $n\in L^{\infty}_t(L^1_x\cap
L^{\infty}_x)$, due to Lemma \ref{lem1} and Lemma \ref{lem2}, we
obtain
\[
v_t, \,\,\nabla^2 v, \,\,c_t, \,\,\nabla^2 c\in
L^{l}_{x,t}\qquad\mbox{ for all}\quad 1<l<\infty.
\]
In our case, $B=u+\chi(c) \nabla c$ and we note, due to parabolic
embedding, that $B$ satisfies the condition \eqref{B}. Therefore, we
conclude that $n\in \calC^{\beta}_{x,t}$ for some $\beta>0$. Due to
classical Schauder estimates, $u$ and $c$ are also in the class
$\calC^{2+\beta, 1+\frac{\beta}{2}}_{x,t}$. This completes the
proof.
\end{thm1.8}
We present proof of Theorem \ref{bdd weak-2}. \\
\begin{thm1.9}\\
$\bullet$\quad $\bke{\text{Case} \ : \ \text{max}\left\{
\text{min}\left\{ 2q-2, \frac{9q-8}{6}\right\}, \ \frac{10q-9}{8}
\right\}<\alpha \ \text{and} \ \frac{3q-1}{6}\le \alpha}$

We note first that this case is reduced to the case that $
\text{min}\left\{ 2q-2, \frac{9q-8}{6}\right\}<\alpha$,
$\frac{3q-1}{6}\le \alpha$ and $q\ge \frac{7}{6}$. We also observe
that, in the case that $\alpha\ge 2q$, the result is already
obtained in Theorem \ref{bdd weak-1}, and therefore, it suffices to
treat the case $\alpha<2q$, that is $ \text{min}\left\{ 2q-2,
\frac{9q-8}{6}\right\}<\alpha$, $\frac{3q-1}{6}\le \alpha<2q$ and
$q\ge \frac{7}{6}$.

As in the previous case \eqref{sept27-30}, following similar
computations, we have
\begin{equation}\label{sept28-20}
\sup_{0 \leq t \leq T}\norm{n}_{p}^{p} \leq \exp \left\{
Cp^2\int_{0}^{T}\norm{\nabla^2 c}_{\frac{6p}{2p+3\alpha+6-6q}}^{2}
\right\}\int_{0}^{T}\norm{\nabla^2
c}_{\frac{6p}{2p+3\alpha+6-6q}}^{2} + \norm{n_0}_{p}^{p}.
\end{equation}
We will show that $\int_{0}^{T}\norm{\nabla^2
c}_{\frac{6p}{2p+3\alpha+6-6q}}^{2}<\infty$ for any $p>\alpha
-2q+3$. Indeed, we recall \eqref{sept27-30} via maximal regularity
for heat equation
\begin{equation}\label{sept28-30}
\begin{split}
\int_{0}^{T}\norm{\nabla^2 c}_{\frac{6p}{2p+3\alpha+6-6q}}^{2} &
\leq
C\int_{0}^{T}\bke{\norm{n}_{\frac{6p}{2p+3\alpha+6-6q}}^{2}+\norm{u\cdot
\nabla
c}_{\frac{6p}{2p+3\alpha+6-6q}}^{2}}+C\norm{\nabla^2 c_0}_{\frac{6p}{2p+3\alpha+6-6q}}^{2} \\
& = C(\text{{\Romannumeral 1}}+\text{{\Romannumeral
2}})+C\norm{\nabla^2 c_0}_{\frac{6p}{2p+3\alpha+6-6q}}^{2} .
\end{split}
\end{equation}
The term $\text{{\Romannumeral 1}}$ is estimated as follows. Since $
\text{min}\left\{ 2q-2, \frac{9q-8}{6}\right\}<\alpha$ and
$\frac{3q-1}{6}\le \alpha<2q$, we have via interpolation inequality
and Sobolev embedding
\begin{equation}\label{sept28-15}
\text{{\Romannumeral 1}}  =\int_{0}^{T}
\norm{n}_{\frac{6p}{2p+3\alpha+6-6q}}^{2} \leq
\int_{0}^{T}\norm{n}_{\alpha-2q+3}^{2\theta}\norm{n}_{3(2\alpha-2q+3)}^{2(1-\theta)}
\leq \int_{0}^{T}\norm{\nabla n^{\frac{2\alpha
-2q+3}{2}}}_{2}^{2-\delta_{p}},
\end{equation}
where
\[
2(1-\theta)=\frac{(2\alpha
-2q+3)\left\{p(-2\alpha+4q)-(\alpha-2q+3)(2p+3\alpha+6-6q)\right\}}{p(5\alpha-4q+6)}
\ \ \text{and}
\]
\[
2(1-\theta)\frac{2}{2\alpha-2q+3}=\frac{4q-2\alpha}{5\alpha-4q+6}-\frac{(\alpha
-2q+3)(2p+3\alpha+6-6q)}{p(5\alpha-4q+6)}=2-\delta_{p}<2.
\]
Thus, it is direct that the term $\text{{\Romannumeral 1}}$ is
finite.

On the other hand, the second term $\text{{\Romannumeral 2}}$ can be
computed exactly as the same way as that of \eqref{sept26-20} in
Theorem \ref{bdd weak-1}, and thus the details are omitted.
\\
\\
$\bullet$\quad $\bke{\text{Case} \ : \ \text{max}\left\{
\text{min}\left\{ 2q-2, \frac{9q-8}{6}\right\}, \ \frac{10q-9}{8}
\right\}<\alpha < \frac{3q-1}{6}}$

We note that this is reduced to the case that
$\frac{10q-9}{8}<\alpha < \frac{3q-1}{6}$ and $1\le q <
\frac{7}{6}$. We set $p$ with
\begin{equation}\label{sept27-50}
\text{max}\left\{ \alpha - 2q + 3, \ 3\alpha-4q+4 \right\} < p <
4\alpha -6q+7.
\end{equation}
Similarly as before, testing $n^{p-1}$ to $\eqref{eq:Chemotaxis}_1$,
we obtain \eqref{sept29-20}, namely
\[
\sup_{0 \leq t \leq T}\norm{n}_{p}^{p}+\int_{0}^{T}\norm{\nabla
n^{\frac{p+\alpha}{2}}}_{2}^{2}
\]
\[
\leq Cp^2\int_{0}^{T}\int_{\mathbb{R}^3}\abs{\nabla n^{p+2q-\alpha
-2}}\abs{ \nabla c} +Cp^2\int_{0}^{T}\int_{\mathbb{R}^3}\abs{
n^{p+2q-\alpha -2}}\abs{ \nabla^2 c} :=
Cp^2\bke{\text{{\Romannumeral 1}} + \text{{\Romannumeral 2}}}.
\]
The first term $\text{{\Romannumeral 1}}$ is the same as
\eqref{sept29-30}.
\[
\text{{\Romannumeral 1}} \leq
\int_{0}^{T}\bke{\varepsilon\norm{\nabla
n^{\frac{p+\alpha}{2}}}_{2}^{2}+C_{\varepsilon}\int_{\mathbb{R}^3}
n^{p+4q-3\alpha -4}\abs{ \nabla c}^{2} }.
\]
Let $r_{1}=\frac{6(\alpha-2q+3)}{14\alpha-22q+27-3p}$. Due to
\eqref{sept27-50}, H\"{o}lder inequality and Sobolev embedding, we
have
\[
\int_{0}^{T}\int_{\mathbb{R}^3} n^{p+4q-3\alpha -4}\abs{ \nabla
c}^{2}  \leq \int_{0}^{T}\norm{n^{p+4q-3\alpha
-4}}_{\frac{\alpha-2q+3}{^{p+4q-3\alpha -4}}}\norm{\abs{ \nabla c
}^2}_{\frac{\alpha-2q+3}{^{4\alpha -6q+7-p}}}
\]
\[
\leq \int_{0}^{T}\norm{n}_{\alpha-2q+3}^{p+4q-3\alpha
-4}\norm{\nabla^2 c}_{r_1}^{2} \leq C\int_{0}^{T}\norm{\nabla^2
c}_{r_1}^{2},
\]
where we used that $\norm{n}_{L^{\infty}_tL^{\alpha-2q+3}_x}<C$
proved in Theorem \ref{weak-2}. From maximal regularity for heat
equation and results in Lemma \ref{velocity-sept12-10}, we have
\[
\int_{0}^{T}\norm{\nabla^2 c}_{r_1}^{2} \leq
C\int_{0}^{T}\bke{\norm{n}_{r_1}^{2}+\norm{u\cdot \nabla c}_{r_1}^2}
+ \norm{\nabla^2 c_0}_{r_1}^{2}
\]
\[
\leq C\int_{0}^{T}\bke{\norm{n}_{r_1}^{2}+
\norm{u}_{6}^{2}\norm{\nabla c}_{\frac{6r_1}{6-r_1}}^{2}} +
\norm{\nabla^2 c_0}_{r_1}^{2} \leq
C\int_{0}^{T}\bke{\norm{n}_{r_1}^{2}+\norm{\nabla^2
c}_{\frac{6r_1}{6+r_1}}^{2}} + \norm{\nabla^2 c_0}_{r_1}^{2}.
\]
\[
\leq
C\int_{0}^{T}\bke{\norm{n}_{r_1}^{2}+\norm{n}_{\frac{6r_1}{6+r_1}}^{2}+\|u\|_{6}^{2}\norm{\nabla
c}_{r_1}^{2}} + C\bke{\norm{\nabla^2 c_0}_{r_1}^{2}+\norm{\nabla^2
c_0}_{\frac{6r_1}{6+r_1}}^{2}}
\]
\begin{equation}\label{sept27-65}
\leq
C\int_{0}^{T}\bke{\norm{n}_{r_1}^{2}+\norm{n}_{1}^{2}+\norm{u}_{6}^{2}\norm{\nabla
c}_{r_1}^{2}}+ C:=\text{{\Romannumeral 1}}_{1}
\end{equation}
If $r_1\le \alpha - 2q + 3$, then \eqref{sept27-65} is bounded, due
to the result of Theorem \ref{weak-2}, by
\begin{equation}\label{sept23-30-2}
\text{{\Romannumeral 1}}_{1}\le C\int_{0}^{T}\norm{\nabla
c}_{r_1}^{2}+ C(1+T).
\end{equation}
On the other hand, in case that $r_1>\alpha - 2q + 3$, we can see
that $r_{1} < 3(p+\alpha)$, since $\text{max}\left\{ \alpha - 2q +
3, \ 3\alpha-4q+4 \right\} < p < 4\alpha -6q+7$, and thus
\eqref{sept27-65} is estimated as
\[
\text{{\Romannumeral 1}}_1 \leq C\int_{0}^{T}\bke{\norm{n}_{\alpha
-2q+3}^{2\theta_{1}}\norm{n}_{3(p+\alpha)}^{2(1-\theta_{1})}
+\norm{\nabla c}_{r_1}^{2}}+C(1+T)
\]
\begin{equation}\label{sept27-70}
\leq C\int_{0}^{T}\bke{\norm{\nabla
n^{\frac{p+\alpha}{2}}}_{2}^{\delta_{1}}+ \norm{\nabla
c}_{r_1}^{2}}+C(1+T),
\end{equation}
where
\[
2(1-\theta_{1})=\frac{(p+\alpha)(3p-14\alpha+22q-21)}{3p+2\alpha+2q-3},
\ \delta_{1}=\frac{2\bke{3p-14\alpha+22q-21}}{3p+2\alpha+2q-3}.
\]
Here we used that $\alpha - 2q + 3 < r_{1} < 3(p+\alpha)$ and
$r_1\ge \frac{6r_1}{r_1+6}$. We note that $\delta_1<2$, since
$\frac{10q-9}{8} <\alpha < \frac{3q-1}{6}$.

Next we estimate the term $\text{{\Romannumeral 2}}$. Let
$r_{2}:=p+2q-\alpha-1$. Using H\"older, Young's and maximal
regularity for heat equation, and following similar computations as
in \eqref{sept27-65}, $\text{{\Romannumeral 2}}$ is estimated as
follows:
\[
\text{{\Romannumeral 2}}
\leq C\int_{0}^{T}\bke{\norm{n}_{r_{2}}^{r_{2}}+\norm{\nabla^2
c}_{r_{2}}^{r_{2}}}\leq C\int_{0}^{T}\bke{\norm{n}_{r_{2}}^{r_{2}} +
\norm{u\cdot \nabla c}_{r_{2}}^{r_{2}}}+C\norm{\nabla^2
c_0}_{r_2}^{r_2}
\]
\[
\leq
C\int_{0}^{T}\left(\|n\|_{r_2}^{r_2}+\|n\|_{\frac{6r_2}{6+r_2}}^{r_2}+\|\nabla
c\|_{r_2}^{r_2}\right) +C\bke{\norm{\nabla^2
c_0}_{r_2}^{r_2}+\norm{\nabla^2 c_0}_{\frac{6r_1}{6+r_1}}^{2}}.
\]
We note that
\[
\int_{0}^{T}\left(\|n\|_{r_2}^{r_2}+\|n\|_{\frac{6r_2}{6+r_2}}^{r_2}+\|\nabla
c\|_{r_2}^{r_2}\right)\le
\int_{0}^{T}\left(C\|n\|_{r_2}^{r_2}+C\|n\|_{1}^{r_2}+\|\nabla
c\|_{r_2}^{r_2}\right)
\]
\[
\leq C\int_{0}^{T}\left(\|n\|_{\alpha
-2q+3}^{r_2\theta_{3}}\|n\|_{3(p+\alpha)}^{r_2(1-\theta_{3})}
+\|\nabla c\|_{2}^{r_{2}\theta_{4}}\|\nabla
c\|_{6}^{r_{2}(1-\theta_{4})}\right)+CT
\]
\[
\leq C\int_{0}^{T}\left(\|n\|_{\alpha
-2q+3}^{r_2\theta_{3}}\|n\|_{3(p+\alpha)}^{r_2(1-\theta_{3})}
+\|\nabla c\|_{2}^{r_{2}\theta_{4}}\|\nabla^2
c\|_{2}^{r_{2}(1-\theta_{4})}\right)+CT
\]
\[
\leq C\int_{0}^{T}\left(\|\nabla
n^{\frac{p+\alpha}{2}}\|_{2}^{\delta_{3}}+ \|\nabla
c\|_{2}^{r_{2}\theta_{4}}\|\nabla^2 c\|_{2}^{\delta_{4}}\right),
\]
where
\[
\delta_{3}=\frac{2r_{2}(1-\theta_{3})}{p+\alpha}=\frac{6(p-2\alpha+4q-4)}{3p+2\alpha+2q-4},\qquad
\delta_4=r_{2}(1-\theta_{4})=\frac{3}{2}(p+2q-\alpha-3).
\]
Here we used that $\alpha -2q+3 < r_{2}< 3(p+ \alpha)$ and
$2<r_{2}<6$ and we observe that $\delta_3<2$ and $\delta_4<2$, since
$\frac{10q-9}{8} <\alpha < \frac{3q-1}{6}$ and $\text{max}\left\{
\alpha - 2q + 3, \ 3\alpha-4q+4 \right\} < p < 4\alpha -6q+7$.
Combining estimates of $\text{{\Romannumeral 1}}$ and
$\text{{\Romannumeral 2}}$, we obtain
\begin{equation*}
\sup_{0 \leq t \leq T}\norm{n(t)}_{p}^{p}+C\int_{0}^{T}\norm{\nabla
n^{\frac{p+\alpha}{2}}}_{2}^{2}  \leq C
\bke{\int_{0}^{T}\norm{\nabla c}_{r_{1}}^{2} +1}.
\end{equation*}
We can see also that  $2<r_{1}<6$ as long as  $\frac{11\alpha
-16q+18}{3} < p < 4\alpha -6q+7$, which is valid, since
$\frac{11\alpha -16q+18}{3}> 3\alpha-4q+4$, in case that
$\frac{10q-9}{8} <\alpha < \frac{3q-1}{6}$. Therefore, for any $p$
with $\text{max}\left\{ \alpha - 2q + 3, \ 3\alpha-4q+4 \right\} < p
< 4\alpha -6q+7$ we obtain
\begin{equation}\label{sept27-75}
n\in L^{\infty}(0, T;L^{p}(\mathbb{R}^3)),\qquad  \nabla
n^{\frac{p+\alpha}{2}} \in L^{2}(0, T;L^{2}(\mathbb{R}^3))
\end{equation}
Let $p_{0}=\frac{6-3\alpha}{4}$. We can see that $1 < p_{0} <
4\alpha -6q+7$ via $\frac{10q-9}{8}< \alpha < \frac{3q-1}{3}$, and
thus it is evident from \eqref{sept27-75} that
\begin{equation}\label{Lp0-1}
n \in L^{\infty}(0, T;L^{p_{0}}(\mathbb{R}^{3})).
\end{equation}
Next, we will show that $n \in L^{\infty}(0,
T;L^{p}(\mathbb{R}^{3}))$ for any $p_0<p<\infty$. Similarly as
before, multiplying equation $\eqref{eq:Chemotaxis}_1$ with
$n^{p-1}$ and using Gronwall inequality, we get \eqref{sept27-20},
namely,
\[
\sup_{0 \leq t \leq T}\norm{n}_{p}^{p} \leq \exp \left\{
Cp^2\int_{0}^{T}\norm{\nabla^2 c}_{\frac{6p}{2p+3\alpha+6-6q}}^{2}
\right\}\int_{0}^{T}\norm{\nabla^2
c}_{\frac{6p}{2p+3\alpha+6-6q}}^{2} + \norm{n_0}_{p}^{p}.
\]
From maximal regularity for heat equation, we have
\begin{equation*}
\begin{split}
\int_{0}^{T}\norm{\Delta c}_{\frac{6p}{2p+3\alpha+6-6q}}^{2} & \leq
C\int_{0}^{T}\bke{\norm{n}_{\frac{6p}{2p+3\alpha+6-6q}}^{2}+\norm{u\cdot
\nabla
c}_{\frac{6p}{2p+3\alpha+6-6q}}^{2}} + \int_{0}^{T}\norm{\Delta c_0}_{\frac{6p}{2p+3\alpha+6-6q}}^{2}\\
& = C(\text{{\Romannumeral 3}}+\text{{\Romannumeral 4}}) +
\int_{0}^{T}\norm{\Delta c}_{\frac{6p}{2p+3\alpha+6-6q}}^{2}.
\end{split}
\end{equation*}
For $p>\alpha -2q+3$, we have
$p_{0}<\frac{6p}{2p+3\alpha+6-6q}<3p_{0}+3\alpha$ and thus, the
terms $\text{\Romannumeral 3}$ and $\text{\Romannumeral 4}$ are
estimated as exactly the same as \eqref{sept29-10} and
\eqref{sept27-40}. Hence, we skip its details. We finally conclude
the boundedness of $L^\infty$-norm of $n$. Indeed, since $n \in
L^{\infty}(0,T;L^{p}(\mathbb{R}^3))$ for all $1 \leq p < \infty$, we
can see that $c_{t}, \ \nabla^c, \ u_{t}$ and $\nabla^{2}u$ belong
to $L^{p}((0,T)\times \mathbb{R}^3)$ for all $p<\infty$ and
therefore, we also note, due to parabolic embedding, that $\nabla c
\in L^{\infty}((0,T)\times \mathbb{R}^3)$. Following similar
procedure as \eqref{L-infiny-1} in Theorem \ref{bdd weak-2}, we have
\begin{equation*}
\frac{d}{dt}\int_{\mathbb{R}^3}\abs{n}^{p}~dx \leq
Cp^2\norm{n}_{p(1-\delta)}^{p(1-\delta)}.
\end{equation*}
Due to exactly same computations as in \eqref{oct03-10}, we obtain
\begin{equation*}
\norm{n(t)}_{\infty} \leq 1+\norm{n(0)}_{\infty}.
\end{equation*}
H\"older continuity can be verified similarly as in Theorem \ref{bdd
weak-1}, and thus the details are skipped. This completes the proof.
\end{thm1.9}

\section{Appendix}

\subsection{Proof of Local energy estimates, Propositions~\ref{P:EE} and ~\ref{P:LEE}}\label{SS:ProofEnergy}
\indent \\
\\
\textbf{Proof of Proposition~\ref{P:EE}}: For a nonnegative bounded
weak solution $n$, set up the test functions
\[
\varphi_{\pm} =  \pm 2 (n - \mu_{\pm} \pm k)_{\pm} \zeta^2,
\]
where $\zeta$ is a piecewise linear cutoff function vanishing on the
parabolic boundary of $Q_{\rho}$. We calculate first that
\[\begin{split}
I &= \iint_{Q_\rho} n_t \varphi_{\pm} \,dx\,dt = \iint_{Q_\rho} \frac{\partial}{\partial t} \left[(n- \mu_{\pm} \pm k)_{\pm}^2\right]\zeta^2 \,dx\,dt \\
  &=  \int_{K_{\rho}\times\{t_1\}} (n- \mu_{\pm} \pm k)_{\pm}^2 \zeta^2 \,dx
      - \int_{K_{\rho}\times\{t_0\}} (n- \mu_{\pm} \pm k)_{\pm}^2 \zeta^2 \,dx\\
  &\quad - 2 \iint_{Q_\rho}(n- \mu_{\pm} \pm k)_{\pm}^2 \zeta \zeta_{t} \,dx\,dt.
\end{split}\]

Now we consider the following integral quantities:
\[
II = \iint_{Q_{\rho}} \nabla n^{1+\alpha} \nabla \varphi_{\pm}
\,dx\,dt = II_1 + II_2
\]
where
\[
II_{1} = 2(1+\alpha)\iint_{Q_{\rho}} \, n^{\alpha} |\nabla
(n-\mu_{\pm} \pm k)_{\pm}|^2 \, \zeta^2 \,dx\,dt,
\]
and (for any $\epsilon_0 > 0$ by the Cauchy-Schwartz inequality)
\[\begin{split}
II_{2}
&= 2(1+\alpha)\iint_{Q_{\rho}} \, n^{\alpha} \nabla (n-\mu_{\pm}\pm k)_{\pm}\,  (n-\mu_{\pm}\pm k)_{\pm}\, \zeta \nabla \zeta \,dx\,dt \\
&\leq 2\epsilon_0 (1+\alpha)\iint_{Q_{\rho}} \, n^{\alpha} |\nabla (n-\mu_{\pm}\pm k)_{\pm}|^2 \, \zeta^2 \,dx\,dt \\
&\quad + 8 \epsilon^{-1}_{0}(1+\alpha) \iint_{Q_{\rho}}  \,
n^{\alpha} (n-\mu_{\pm}\pm k)^{2}_{\pm} \, |\nabla \zeta|^2 \,dx\,dt
= II_{21} + II_{22}.
\end{split}\]
The first term on the right hand side is absorbed to $II_1$ by
choosing $\epsilon_0 = 1/2$.

Now we consider integral terms carrying the lower order term,
\[
III = \iint_{Q_{\rho}} \, \nabla ( B \, n) \, \varphi_{\pm} \,dx\,dt
= III_1 + III_2
\]
where
\[
III_1 = \iint_{Q_{\rho}} \, \nabla B \, n \, \varphi_{\pm} \,dx\,dt,
\]
and (by taking integration by parts with respect to the space
variable)
\[\begin{split}
III_2
&= \iint_{Q_{\rho}} \, B \, \nabla n \, \varphi_{\pm} \,dx\,dt = \iint_{Q_{\rho}} \, B \, \left[ \nabla (n-\mu_{\pm}\pm k)^{2}_{\pm}\right] \zeta^2 \,dx\,dt \\
&= - \iint_{Q_{\rho}} \, \nabla B \,  (n-\mu_{\pm}\pm k)^{2}_{\pm}
\zeta^2 \,dx\,dt -2 \iint_{Q_{\rho}} \, B \,  (n-\mu_{\pm}\pm
k)^{2}_{\pm} \zeta \nabla \zeta \,dx\,dt\\
& = III_{21} + III_{22}.
\end{split}\]
Then by applying the Cauchy-Schwartz inequality, we have
\[\begin{split}
III_{22}
&\leq \iint_{Q_{\rho}} \, k^{\alpha}  (n-\mu_{\pm}\pm k)^{2}_{\pm} |\nabla \zeta|^2 \,dx\,dt \\
&\quad + \iint_{Q_{\rho}} \, |B|^2 k^{-\alpha} \,  (n-\mu_{\pm}\pm
k)^{2}_{\pm} \zeta^2 \,dx\,dt = III_{221} + III_{222},
\end{split}\]
where the first term, $III_{221}$, and $II_{22}$ are collected
together in \eqref{EE} (the third term on the right-hand-side).
Notice that
$ (n - \mu_{\pm} \pm k)_{\pm} \leq k.$
From condition \eqref{B}, we compute
\[\begin{split}
III_{222} &\leq k^{2-\alpha} \|B\|^{2}_{2 \hat{q}_1, 2\hat{q}_2} \left[
\int_{t_0}^{t_1} \left[A^{\pm}_{k,\rho}(t)\right]^{\frac{q_2}{q_1}}
\,dt\right]^{\frac{2(1+\kappa)}{q_2}}, \\
III_1 + III_{21} &\leq \left(2k\mu_{+} + k^2\right)\|\nabla
B\|_{\hat{q}_1, \hat{q}_2} \left[ \int_{t_0}^{t_1}
\left[A^{\pm}_{k,\rho}(t)\right]^{\frac{q_2}{q_1}}
\,dt\right]^{\frac{2(1+\kappa)}{q_2}}.
\end{split}\]
\smallskip
\indent
\\
\\
\textbf{Proof of Proposition~\ref{P:LEE}} :
Due to the setting of function $\Psi_{\pm}$ in \eqref{Psi}, we compute that
\[
\Psi'_{\pm} = \frac{\pm 1}{(1+\delta)k - (n-\mu_{\pm} \pm k)_{\pm }}
\ \text{ and } \ \Psi''_{\pm} = \left( \Psi'_{\pm}\right)^2.
\]

For a nonnegative solution $n$, set up the test function
$$ \varphi_{\pm} = 2 \Psi_{\pm}(n) \Psi'_{\pm}(n) \zeta^2 . $$
Then we observe that
\[\begin{split}
I &= \int_{t_0}^{t_1}\int_{K_{\rho}} n_t \varphi_{\pm} \,dx\,dt = \int_{t_0}^{t_1} \int_{K_{\rho}} \frac{d}{dt}\left[ \Psi^{2}_{\pm}(n) \zeta^2 \right] \,dx\,dt \\
&= \int_{K_{\rho} \times \{t_1\}} \Psi^{2}_{\pm}(n) \zeta^2 \,dx -
\int_{K_{\rho} \times \{t_0\}} \Psi^{2}_{\pm}(n) \zeta^2 \,dx.
\end{split}\]

Using the properties of $\Psi_{\pm}$, note that
$$ \nabla \varphi = 2 \nabla n\left( 1 + \Psi_{\pm} (n)\right) (\Psi'_{\pm}(n))^2 \zeta^2 + 4 \Psi_{\pm}(n)\Psi'_{\pm}(n) \zeta \nabla \zeta. $$

Then we calculate various integral quantities:
$$ \int_{t_0}^{t_1}\int_{K_{\rho}} \nabla n^{1+\alpha} \nabla \varphi_{\pm} \,dx\,dt = II_1 + II_2$$
where
\[\begin{split}
II_1 &= 2\int_{t_0}^{t_1}\int_{K_{\rho}} \nabla n^{1+\alpha} \nabla n\left( 1 + \Psi_{\pm} (n)\right) (\Psi'_{\pm}(n))^2 \zeta^2  \,dx\,dt \\
&= 2(1+\alpha)\int_{t_0}^{t_1}\int_{K_{\rho}} n^{\alpha} |\nabla
n|^2 \left( 1 + \Psi_{\pm} (n)\right) (\Psi'_{\pm}(n))^2 \zeta^2
\,dx\,dt,
\end{split}\]
and
\[\begin{split}
II_2 &= 4\int_{t_0}^{t_1}\int_{K_{\rho}} \nabla n^{1+\alpha} \Psi_{\pm}(n)\Psi'_{\pm}(n) \zeta \nabla \zeta \,dx\,dt \\
&\leq  \epsilon_0 (1+\alpha)\int_{t_0}^{t_1}\int_{K_{\rho}} n^{\alpha} |\nabla n|^2 |\Psi'_{\pm}(n)|^2 \zeta^2  \,dx\,dt \\
&\quad + 16 \epsilon_{0}^{-1} (1+\alpha)\int_{t_0}^{t_1}\int_{K_{\rho}} n^{\alpha} \Psi^{2}_{\pm}(n) |\nabla \zeta|^2  \,dx\,dt\\
&= II_{21} + II_{22}
\end{split}\]
applying the Cauchy-Schwartz inequality with any $\epsilon_0 > 0$.
By fixing $\epsilon_0 = 2$, the first integral term is absorbed by
$II_1$.

Next, we handle the integral quantity carrying the lower order term:
\begin{equation}\label{LEE01}
\int_{t_0}^{t_1}\int_{K_{\rho}} \nabla (B n) \varphi \,dx\,dt =
III_1 + III_2
\end{equation}
using $\nabla B$. The integral \eqref{LEE01} produces two terms
\[
III_1 = \int_{t_0}^{t_1}\int_{K_{\rho}} \nabla B n \varphi \,dx\,dt
\leq 2\int_{t_0}^{t_1}\int_{K_{\rho}} |\nabla B| n \Psi_{\pm}
|\Psi'_{\pm}|\zeta^2 \,dx\,dt
\]
and
\[\begin{split}
III_2 &= \int_{t_0}^{t_1}\int_{K_{\rho}} B \nabla n \varphi \,dx\,dt
= \int_{t_0}^{t_1}\int_{K_{\rho}} B  \left[\nabla \Psi^{2}_{\pm}\right] \zeta^2 \,dx\,dt \\
&= - \int_{t_0}^{t_1}\int_{K_{\rho}} \nabla B  \Psi^{2}_{\pm} \zeta^2 \,dx\,dt - 2\int_{t_0}^{t_1}\int_{K_{\rho}} B \Psi^{2}_{\pm} \zeta \nabla \zeta \,dx\,dt \\
&= III_{21} + III_{22}
\end{split}\]
using the integration by parts. Then by applying the Cauchy-Schwartz
inequality, we have
\[\begin{split}
III_{22}
&\leq 2\int_{t_0}^{t_1}\int_{K_{\rho}} k^{\alpha} \Psi^{2}_{\pm} |\nabla \zeta|^2 \,dx\,dt \\
&\quad + 2\int_{t_0}^{t_1}\int_{K_{\rho}} |B|^2 k^{-\alpha}
\Psi^{2}_{\pm} \zeta^2 \,dx\,dt = III_{221} + III_{222}
\end{split}\]
where $III_{221}$ is bounded by $II_{22}$.

From that $(n-\mu_{\pm}\pm k)_{\pm} \leq k$ and \eqref{B}, we have
\[\begin{split}
III_{222} &\leq 2 k^{-\alpha}\left(\ln \frac{1}{\delta}\right)^2
\|B\|^{2}_{2 \hat{q}_1, 2\hat{q}_2} \left[ \int_{t_0}^{t_1}
\left[A^{\pm}_{k,\rho}(t)\right]^{\frac{q_2}{q_1}}
\,dt\right]^{\frac{2(1+\kappa)}{q_2}},\\
III_1 + III_{221} &\leq \left( (\ln \frac{1}{\delta})^2 +
\frac{\mu_{+}\ln \frac{1}{\delta}}{\delta k}\right) \|\nabla
B\|_{\hat{q}_1, \hat{q}_2} \left[ \int_{t_0}^{t_1}
\left[A^{\pm}_{k,\rho}(t)\right]^{\frac{q_2}{q_1}}
\,dt\right]^{\frac{2(1+\kappa)}{q_2}}.
\end{split}\]
\\
\section*{Acknowledgements}
Sukjung Hwang's work is supported by NRF-2015R1A5A1009350. Kyungkeun
Kang's work is supported by NRF-2014R1A2A1A11051161 and
NRF-2015R1A5A1009350. Jaewoo Kim's work is supported by
NRF-2015R1A5A1009350.

\end{document}